\documentclass[12pt,a4paper,oneside]{amsart}

\usepackage{amsfonts, amsmath, amssymb, amsthm, hyperref}
\usepackage{anysize}

\usepackage[T2A]{fontenc}
\usepackage[cp1251]{inputenc}
\usepackage[english,russian]{babel}
\newtheorem*{theorem*}{Theorem}

\theoremstyle{definition}

\theoremstyle{remark}

\numberwithin{equation}{section}

\renewcommand{\epsilon}{\varepsilon}
\renewcommand{\phi}{\varphi}
\renewcommand{\kappa}{\varkappa}

\begin{document}

\title{Turing machine interaction problem}
\author{Marsel~Matdinov}
\email{mlmatdinov@gmail.com}
\address{Marsel~Matdinov, Moscow Institute of Physics and Technology, 141701 Moscow region, Russia}

\begin{abstract}
The article introduces some ideas for solving special cases of the following problem, proposed in a somewhat generalized form by Marcus Hutter in 2000. Given two Turing machines $A$ and $C$, it is required to build a Turing machine $B$, such that after interacting of $A$ and $B$ on a shared tape for a fixed number of iterations, the machine $C$ outputs 1 on the communication protocol of $A$ and $B$. Details in the introduction.
\end{abstract}

\maketitle

\sloppy

%\Large

\tableofcontents

\newpage

\section{Введение.}

Если представить задачи, описанные мной в ~\cite{old_text}, в обобщенном виде, получится следующая задача. На вход дается машина Тьюринга $N$ и числа $k$ и $r$. Требуется эффективно построить эффективно работающую машину Тьюринга $М$, удовлетворяющую следующему свойству. Предположим, что вначале на основную ленту выписывается случайная строка $a_1$ длины $r$. Машина $М$ обрабатывает эту строку на других лентах (ей предоставляются другие, доступные ей ленты) и дописывает справа строку $b_1$ длины $r$. Затем справа на основную ленту дописывается случайная строка $a_2$ длины $r$. Машина $М$ обрабатывает все уже выписанные на основную ленту строки и дописывает справа строку $b_2$ длины $r$. Так продолжается до тех пор, пока не будут выписаны строки $a_k$ и $b_k$. После этого машина $N$ обрабатывает все выписанные на основную ленту строки и выдает результат - 1 или 0. Требование на машину $М$ в том, чтобы машина $N$ гарантированно выдавала 1 в тех случаях, когда этого можно добиться. (Фиксирование длин строк $a_i$ и $b_i$ происходит скорее из особенностей моего подхода к решению, в дальнейшем от одинаковости длин можно и нужно отказываться.) Хочется научиться строить алгоритм, решающий задачу для как можно более широкого класса машин $N$ (класса входных данных, если быть точнее), как можно более близкого к классу машин $N$, для которых построение машины $М$ под силy человеку за разумное время. 

Статья является продолжением рассуждений из статьи ~\cite{old_text}, все формулировки задач и определения терминов остаются в силе.

Задача представляется актуальной поскольку многие задачи дискретной математики и theoretical computer science формулируются на языке указанной задачи. 

Задача обобщается AIXI формализмом, предложенным Маркусом Хаттером ~\cite{AIXI}. Если быть точнее, той задачей, которая перед AIXI-агентом стоит. Сам этот агент в статье рассмотрен не будет. В конце введения к прошлой статье, написана формулировка AIXI задачи.

Если вкратце, строки $\{ A_i \}$ могут выбираться не случайно, а выписываться другой машиной Тьюринга. Таким образом, задача искомой нами машины Тьюринга в том, чтобы правильно взаимодействовать с другой машиной Тьюринга.

Так же формулировка задачи очень напоминает определение полиномиальной иерархии ~\cite{alternation}.

Готовый алгоритм пока не представлен, речь об идеях. Так же не представлен, а может быть и никогда не будет представлен конкретный четко очерченный подкласс задач указанного вида, на которые разрабатываемый алгоритм нацелен. В общем случае, если не ограничивать длины строк константой $r$, задача, конечно, алгоритмически неразрешима. Даже когда мы ограничиваем длину строк и время работы машины $N$, нам едва ли удастся находить решение эффективно всегда: скажем, некоторые NP-полные задачи можно легко сформулировать как экземпляры упомянутой в прошлом тексте большой задачи о выписывании строки. Но это все не так плохо: подкласс задач указанного класса, определяемый задачами, которые решают люди, тоже далек от того, чтобы быть четко очерченным и возможно не столь широк. Думаю, что в нашем случае, попытка четко определить класс решаемых задач так, чтобы принадлежность задачи этому классу можно было бы выяснить "сразу"{}, с требованием гарантированности нахождения решения, скажем, за полиномиальное время, нас бы сильно ограничивали. Нужно отталкиваться от интуиции самого алгоритма, пусть пока и не достроенного.

\section{Другой взгляд на пространство признаков.}

До этого я писал о признаках как об объектах, каждый из которых задается машиной Тьюринга и еще несколькими натуральными числами. Можно поступить иначе.

Для определенности, будем говорить о большой задаче о выписывании строки.

Можно решать эту задачу следующим образом. Есть схема, приходящая из машины Тьюринга, принимающая 2 строки $A$ и $X$. Наша цель - научиться по строке $A$ строить строку $X$, чтобы $M(A,X) = 1$. Рассмотрим прямоугольный паралеллепипед в трехмерной решетке, на узлах которого определены несколько предикатов. Верхнюю из сторон паралеллепипеда мы будем интерпретировать как историю работы $M$ (это прямоугольник, значения предикатов на узлах которого соответствуют значениям соответствующих гейтов в истории работы $M$). Весь паралеллепипед мы будем интерпретировать как схему, вычисляющую историю работы $M$ снизу вверх (каждое значение предиката гейта вычисляется по значениям предикатов в нескольких гейтах непосредственно под данным гейтом). Должно быть выполнено ограничение, что значения предикатов, соответствующие входной строке $A$ не меняются в процессе работы схемы-параллелепипеда: каждое значение предиката этой строки равно одному и тому же значению на любой высоте. Нижняя сторона паралеллепипеда (вход схемы) устроена следующим образом. Гейты, соответствующие строке $A$ заполняются строкой $A$ (для которой мы хотим получить строку $X$). Все остальные гейты заполняются нулями. В качестве выхода мы хотим получить корректную историю работы схемы $M$, которая в итоге выдает единицу. Первая входная строка этой полученной истории, как мы уже договорились, должна быть равна $A$, а вторая - некоторой строке $X$, которую мы и выдадим в качестве ответа. Таким образом, наша цель - предъявить способ вычисления предикатов на гейте параллелепипеда по предикатам на всех гейтах, расположенных непосредственно под ним так, чтобы полученная схема корректно работала для любой строки $A$.

Итак, вначале у нас есть множество гейтов паралеллепипеда. Предположим, у нас есть некоторый очень большой универсум, из которого мы можем доставать гейты и добавлять в наше множество (назовем его $\Omega$). Так же, у нас есть множество предикатов $\Delta_1$, которое мы так же можем расширять, при желании вводя новые предикаты. Стоит изначально завести несколько базовых предикатов, которые мы будем использовать всегда, таких как унарный предикат "являться гейтом паралеллепипеда"{}, бинарный предикат, равный единице, когда два гейта во-первых, принадлежат паралеллепипеду, а во-вторых, имеют одинаковую координату (такой предикат есть для каждой из трех координат), унарные предикаты, равные единице, когда гейт принадлежит паралеллепипеду и заданная координата равна своему минимальному или максимальному значению (соответственно, таких предикатов 6), бинарный предикат равенства, те предикаты, которые присутствуют на паралеллепипеде изначально и значения которых эта трехмерная схема выдает в качестве ответа.

Так же, нам понадобится множество функций $\Delta_2$, сопоставляющих каждому набору гейтов соответствующей арности некоторый гейт. Такой гейт - образ функций может быть добавлен в наше множество из универсума по нашему запросу. Если образу функции надлежит совпадать с одним из уже добавленных в $\Omega$ гейтов, мы можем сделать единицей значение на них предиката равенства, а так же, возможно, склеить их друг с другом. Вполне возможно, что в итоге без функций можно будет обойтись, оставив лишь предикаты. 

Назовем ситуацией набор значений всех имеющихся на данный момент предикатов на всех возможных наборах гейтов (размера соответствующей арности) текущего множества гейтов. Мы будем поддерживать модель, генерирующую ситуации и, тем самым, задающую распределение вероятостей на ситуациях. Какой именно будет эта модель? Будем считать, что это та же самая модель, что обычно используется в генеративно-состязательных сетях ~\cite{GAN}: глубокая нейронная сеть, принимающая на вход seed из нескольких независимых нормально распределенных случайных величин, выдающая на выходе по числу из промежутка [0,1] для каждого из предикатов и каждого набора гейтов количеством равным арности этого предиката. Это число интерпретируется как вероятность того, что данный предикат на данном наборе гейтов выдаст 1. Когда эти вероятности выданы нейросетью, все бинарные значения предикатов выбираются независимо с соответствующими вероятностями.

И мы будем изменять веса данной нейронной сети так, чтобы соответствующее распределение все в большей степени обладало некоторыми свойствами, о которых речь пойдет ниже.

Конечно, возможен и другой выбор порождающей модели (например, были соображения, что для наших целей могут хорошо подойти машины Больцмана). О желаемых свойствах этой модели поговорим чуть позже.

Признаки, вычисляемые машинами Тьюринга и история вычисления которых представляет из себя решетку определенной размерности стоит заменить на изначально не структурированное множество гейтов как раз для того, чтобы "уйти от решетки". Если бы мы мыслили "старый"{} признак как нечто, последовательно вычисляемое вдоль некоторой решетки, то в процессе изменения законов, по которым каждый следующий слой гейтов этой решетки вычисляется по предыдущему, "расположенность"{} гейтов истории вычисления в некоторой решетке нас бы несколько сковывала. В некотором смысле, отказ от решетки в пользу изначально не структурированного множества гейтов является выходом в объемлющее пространство (к слову, можно насчитать довольно много примеров выхода в объемлющее пространство в том или ином виде, возникших при развитии нашей системы).

При этом, если мы захотим, чтобы некоторое множество гейтов таки образовывало решетку, нам достаточно сделать следующую вещь. Назовем расширенным предикатом пропозициональную формулу от некоторого множества предикатов из $\Delta_1$ от установленных элементов некоторой оболочки набора гейтов (количество гейтов набора равно арности расширенного предиката), где под оболочкой множества гейтов имеется в виду упорядоченное множество самих гейтов, некоторое множество гейтов - образов функций из $\Delta_2$ от этих самых, исходных, гейтов, так же, образов функций из $\Delta_2$ от исходных гейтов и образов функций из $\Delta_2$ от исходных гейтов, и так далее. (Возможно, я допустил вольность речи определяя расширенный предикат словами, а не формулой, но формула бы получилась громоздкой.) Так вот, чтобы добиться решетки, нам достаточно потребовать равенство единице некоторых расширенных предикатов от всех возможных наборов гейтов в количестве, равном арности данного предиката. Каких именно предикатов, приводить не буду, чтобы не загромождать текст, это нам не понадобится. Алгоритм действительно допускает рассмотрение расширенного предиката от всех возможных наборов гейтов и форсирование его тождественного равенства единице всюду (постепенное изменение весов порождающей модели так, чтобы он был равен единице в итоге). Решетка как-бы формируется, "возникает из воздуха".

Итак, мы ввели "улучшенную версию"{} пространства признаков, тем самым введя в рассмотрение новый объект и хочется называть его иначе. Поэтому давайте будем называть его "пространство гейтов".

\section{Внешняя нейронная сеть.}

Рассмотрим апгрейд нашего алгоритма, заключающийся в том, что мы будем делать еще и вот что. Рассмотрим текущий набор гейтов с текущим распределением значений предикатов в них. И рассмотрим нейронную сеть, которую мы будем называть внешней, принимающую на вход все эти значения (ситуацию), а так же, еще несколько (небольшое число) унарных предикатов, каждый из которых равен единице ровно в одном из гейтов и равен нулю во всех остальных (мы будем называть такие предикаты правильными), и выдающую вещественное число. Мы будем генерировать ситуации и значения правильных предикатов и изменять веса нейронной сети, используя backpropagation таким образом, чтобы выдаваемые ей значения, эксклюзивно (сейчас поясню) коррелировали друг с другом. И остановимся в тот момент, когда они начнут эксклюзивно коррелировать друг с другом довольно хорошо (есть некая вольность в том, когда именно нужно остановиться). Под эксклюзивностью подразумевается то, что значения коррелируют, но когда мы рассматриваем значения нейросети на случайных ситуациях из равномерного распределения на ситуациях, то никакой корреляции уже не наблюдается или она будет слабее. После этого, мы начнем изменять с помощью backpropagation веса уже не внешней, а порождающей нейронной сети так, чтобы значения внешней нейронной сети коррелировали друг с другом еще сильнее (есть более правильная версия этой идеи, заключающаяся в том, что мы начинаем изменять веса внешней сети, так, чтобы значения эксклюзивно коррелировали, в некоторый момент начинаем подключать и веса порождающей сети, далее, веса внешней сети мы меняем все меньше, а веса порождающей сети мы меняем все больше - с целью, чтобы значения внешней сети коррелировали), и, опять же, в некоторый момент останавливаемся. Имеется в виду, конечно, что у порождающей нейросети есть еще и другие цели, которых она хочет достичь и сохранить достигнутое - изменение ее весов лишь немного подстраивается под еще и вот эту новую внешнюю нейросеть, старые цели (например, другие внешние нейросети) не перестают быть интересны. Возможно, в определение внешней нейросети стоит добавить, что мы копируем ее входные данные (ситуацию и набор правильных предикатов) и вставляем ее в каждый слой внешней нейросети (тем самым, каждый ее узел зависит еще и от некоторых значений входной ситуации и правильных предикатов).

Это является обобщением того, что я писал в прошлом препринте о том, как мы форсируем равенство единице значения признака, на разных аргументах коррелирующего с единицей. Поясню на примере. Допустим, нам дана задача где входными данными являются координаты центров и длины сторон $n$ квадратов на клеточной решетке, тем самым, задающие сами эти квадраты, и нужно с этими квадратами что-то сделать. Допустим, в процессе проверки решения наша машина Тьюринга рисует на решетке эти квадраты и оперирует еще одним унарным предикатом $P$ на узлах этой решетки. В какой-то момент процесса изменения весов порождающей нейросети может оказаться, что значения этого предиката в правых верхних вершинах всех квадратов коррелирует с единицей. Как может быть устроена внешняя нейросеть? Предполагаю, что нейросети скорее способны на то, чтобы научиться по входному правильному предикату, в котором гейт - единичка проецируется в $i$-ю ячейку входной строки, выдавать значение $P$ в правой верхней вершине $i$-го квадрата. Форсирование равенства значения этой нейросети единице как раз и ожидаемо реализует форсирование равенства значения $P$ единице.

Конструкция внешней нейросети напоминает конструкцию определенного в более ранних секциях признака. Мы и здесь (в случае с квадратами) могли бы сказать, что приближение к единице значений $P$ в правых верхних вершинах квадратов равносильно приближению к единице значения соответствующего признака (легко придумать признак, который по числу $i$ выдает значение $P$ в правой верхней вершине $i$-го квадрата), но такой признак нужно еще сформировать.

В общем случае, в случае, когда мы копируем ситуацию и правильные предикаты во все слои внешней сети, можно условно выделить некоторый промежуточный слой внешней сети и посмотреть на значения тех ее узлов, которые изменяются и их соотношение с узлами ситуации и правильными предикатами. Они будут образовывать некоторую, назовем это так, конструкцию. То, что значения внешней сети стали эксклюзивно коррелировать, означает, что эксклюзивно коррелируют значения некоторого признака таких конструкций (не путать с признаками, определенными в более ранних секциях). Когда мы изменяем веса порождающей модели так, чтобы значения внешней сети коррелировали еще сильнее, этот признак данных конструкций начинает коррелировать еще сильнее. Таким образом, мы, в каком-то смысле делаем похожие конструкции еще более похожими.

Развитие темы внешней нейронной сети есть в разделе "Еще мысли".

\section{Вещественнозначные гейты и лампочки.}

Чтобы конструкция внешней нейронной сети лучше вписалась в конструкцию пространства гейтов, полезно наряду с гейтами и лампочками, принимающими значение 0 или 1, ввести в рассмотрение гейты и лампочки, принимающие любые возможные вещественные значения. С ними мы будем обращаться точно так же: порождающее распределение будет генерировать на них уже вещественные значения. Теперь, внешняя нейронная сеть задается добавлением в пространство гейтов набора вещественнозначных гейтов и форсирования, для каждого нейрона, равенства соответствующим константам лампочек, олицетворяющих разности суммы произведений значений гейтов, подаваемых на вход нейрона на соответствующие веса (веса тоже задаются гейтами) и значения гейта, соответствующего элементу нейрона, передаваемого на вход нелинейной функции; а так же, форсирование, для каждого нейрона, равенства нулю лампочки, олицетворяющей разность значения гейта, соответствующего выходному элементу нейрона и нелинейной функции от значения гейта, соответствующего тому элементу нейрона, к которому нелинейная функция должна применяться.

Я только что написал о вещественнозначных лампочках, не определив их. В бинарном случае, каждой лампочке соответствует набор значений соответствующих этой лампочке гейтов и если гейты и правда принимают такие значения, то лампочка загорается. В вещественнозначном случае, правильней говорить не "загорается"{}, а "принимает вещественнозначное значение": каждому возможному набору вещественных значений соответствующих вещественнозначной лампочке гейтов соответствует вещественное значение лампочки; распределения на значениях лампочек можно точно так же изменять подобными бинарному случаю методами.  

Только что было описано, как реализовать внешнюю сеть, не зависящую ни от одного из правильных предикатов. Чтобы добавить один или несколько правильных предикатов, следует говорить о вещественнозначных предикатах, принимающих не только гейт внешней нейронной сети, но и один или несколько других гейтов (нас будут интересовать значения таких предикатов не всюду, а лишь когда эти один или несколько других гейтов лежат в ситуации - множестве, так сказать, основных гейтов; скажем, если один из этих гейтов лежит во внутреннем слое внешней нейронной сети, нас значение предиката интересовать не будет). Для каждого допустимого набора значений правильных предикатов, значение нейросети вычисляется абсолютно аналогичным образом, с участием значений указанных вещественнозначных предикатов на наборах гейтов, состоящих из одного гейта с одного из слоев нейросети, а так же, нескольких гейтов, в количестве равном числу правильных предикатов, находящихся ровно там, где предикаты принимают единичное значение. Но можно поступить проще (если не возникнет никаких проблем) и задать правильные предикаты просто добавлением нескольких новых предикатов и, добавлением соответствующих ограничений добиться того (далее по тексту будет понятно, что имеется в виду), что распределение будет выдавать в каждом из этих новых предикатах всегда (или почти всегда) ровно одну единичку, причем так, чтобы получаемый таким образом набор значений правильных предикатов принимал все (или почти все) допустимые наборы значений и делал это довольно равномерно и независимо от значений остальных предикатов.

Думаю, что понятие коциклического многочлена можно определить для случая вещественнозначных гейтов по аналогии, без особых проблем. Кстати, если говорить о коциклическом многочлене для случая бинарных гейтов, то его, как сумму значений лампочек с определенными коэффициентами, можно легко вычислить по ситуации агрегатом из нескольких вещественнозначных гейтов так, что значение этого коциклического многочлена станет просто значением некоторого гейта.

Скорее всего, в итоговой версии алгоритма останутся только вещественнозначные гейты, поддерживать два типа гейтов - не очень элегантно.

\section{Форсирование симметрии.}

В этой подсекции я укажу на связь между казалось бы разными элементами алгоритма, которая их объединяет, сводя "под одну крышу"{}, и мы понимаем, что на самом деле, делая эти апиори разные вещи, мы на самом деле делаем одно и то же. Речь пойдет, прежде всего, о форсировании значения коциклического многочлена значению, которое он должен принимать и форсировании значения внешней нейронной сети одному и тому же значению всюду.

Утверждается, что в обоих случаях мы приближаем значения соответственных значений предикатов семейства одинаковых (или хотя-бы похожих) конструкций к их среднему значению, либо на самой задаче, либо на близкой задаче - что такое близкая задача я скоро скажу. (По сути, я уже писал о такой возможности в прошлой статье, когда описывал семейства аналогичных конструкций, на соответственные элементы которых навешиваются метки.) Будет видно, что первое - частный случай второго.

Когда мы приближаем и без того несколько близкие значения внешней нейронной сети, элементами конструкций данного семейства конструкций являются лампочки, задействованные в вычислении выходного значения сети: для каждого seed-а и каждого набора значений правильных предикатов, значения всех задействованных лампочек, в том числе и веса нейросети образуют такую конструкцию. Их похожесть заключается в одинаковости для всех случаев весов нейросети, выполнении одних и тех же соотношений для значений лампочек, учавствующих в одном нейроне (это, опять же, напрямую следует из одинаковости весов) и корреляции значений лампочек, в которых записано выходное значение сети. Элементами, значения которых мы приближаем, в данном случае, друг к другу, то есть к усредненному значению на самой задаче, являются значения лампочек, в которых записано выходное значение сети. (Сейчас был озвучен второй вариант того, что можно понимать под конструкциями семейства в случае внешних сетей, первый был упомянут в конце секции "Внешняя нейронная сеть".)

Когда мы приближаем значение коциклического многочлена к значению, которое он должен принимать, для того, чтобы можно было говорить о приближении каких-то значений к соответственным значениям аналогичных конструкций, следует рассмотреть ослабленную задачу. Поясню, что это значит, и зачем это нужно. Дело в том, что когда вопрос встает о том, чтобы двигать значения предикатов в сторону усредненного значения этого предиката по соответственным элементам некоторого семейства аналогичных конструкций, вообще говоря, можно рассматривать эти самые аналогичные конструкции возникающие при рассмотрении не только данной задачи, но и конструкции, возникающие при рассмотрении близких задач, - наша задача является задачей удовлетворения некоторого семейства ограничений на лампочки, но можно рассматривать и пытаться решать задачу, в которой некоторые из этих ограничений сняты, задачу, где добавлены некоторые новые ограничения, а возможно, и задачу, в которой и сняты и добавлены некоторые ограничения. 

Так вот, в случае приближения значения коциклического многочлена к нужному, нашу задачу (давайте, для простоты, говорить о малой задаче о выписывании строки) можно рассматреть как задачу выбора значений лампочек так, чтобы, во-первых, те значения лампочек, которые защищены, были равны нулю (это условие задачи: лампочки, соответствующие наборам значений гейтов, запретным с точки зрения предписаний, задающих схему, данную в задаче - правил, по которым значения гейтов следующего слоя схемы должны вычисляться по значениям гейтов предыдущего - должны обнуляться), а во-вторых, чтобы значения лампочек образовывали согласованную ситуацию (что так же можно записать как набор пропозициональных ограничений на значения лампочек). И ограничения первого вида можно просто отбросить, оставив лишь ограничения для согласованности. Решая эту ослабленную задачу, можно выделить эту самую конструкцию - коциклический многочлен (если он реализован схемой из вещественнозначных гейтов, вычисляющих значение многочлена по ситуации, то конструкцией будет эта самая схема, а значение, указанное в элементе, выдающем значение многочлена, уже в нашей задаче - не ослабленной - и будет значением, которое приближается к усредненному значению для большинства конструкций семейства в ослабленной задаче, то есть, по сути, к значению, которое и должен принимать данный коциклический многочлен). На самом деле этот самый коциклический многочлен может быть "зашит"{} в показатель пустой конструкции - определяемый ниже объект, и тогда, при нашем варианте выбора нужных определений, наше семейство аналогичных конструкций - это семейство пустых конструкций и вся информация зашита в показателе (ниже будет понятно, о чем речь).

Итак, мы приближаем значения соответственных элементов похожих конструкций в нашей задаче к усредненному значению соответственных значений таких же конструкций, либо на самой нашей задаче, либо на близкой задаче. Но очевидно, что сама задача близка сама к себе - множество требуемых ограничений просто совпадает, поэтому, первое - частный случай второго.

Решая близкую задачу, нужно стараться, чтобы полученное распределение на решениях близкой задачи было как можно более общим, как можно более широким. Что значит последнее, будет написано ниже, в других секциях.

В этот же класс попадает возможность переносить информацию с одних конструкций на другие, аналогичные ей, даже в случае, когда этих конструкций мало (возможна, и не редка, ситуация, когда стоит переносить значения с одной конструкции на другую, аналогичную ей - то есть, когда их две). Имеется в виду, что нам удалось найти некоторое, возможно маленькое и никак не структурированное множество похожих конструкций - множеств одинакового размера, между элементами которых есть соответствие, такое, что значения некоторых соответственных лампочек на этих множествах всегда равны друг другу или хотя-бы коррелируют (можно, так же, например, говорить об одинаковости постоянных пропозициональных связей между лампочками этих множеств, для соответственных наборов лампочек, но наличие связи тоже можно задавать как предикат), и мы форсируем равенство друг другу и некоторых других соответственных значений (возможно других) лампочек. На самом деле, в общем случае, понятие конструкции стоит понимать в более стохастичном смысле, см. секции "Внешняя нейронная сеть"{}, "Еще мысли". Скорее, чтобы задать конструкцию, нужно задать несколько дополнительных предикатов и для каждого из них и для каждого набора лампочек в размере соответствующей арности задать вещественное число от 0 до 1 - степень принадлежности набора лампочек "в аспекте этого предиката"{} конструкции.

Многие рассмотренные до сих пор конструкции, сами по себе, как наборы лампочек, были постоянными, не менялись с изменением seed. На самом деле, мы будем относиться к конструкции, в основном, как к чему-то, что живет в конкретной ситуации, как это было с первым вариантом определения конструкции в случае внешней нейросети. В таком случае, при разных значениях seed, мы получаем разные семейства конструкций.

То, что происходит во всех описанных случаях можно объединить фразой "форсирование симметрии". 

Кстати, к форсированию симметрии можно отнести и нечто, что можно назвать форсированием определенности. То, о чем я писал в прошлом препринте - актуальными будут те признаки, которые имеют много связей с другими актуальными признаками, - наличие связи - то есть наличие некоторого каким-то образом вычисляемого предиката, который оказался равен константе. Равенство предиката константе и есть определенность. Так же определенностью можно было бы назвать, например, равенство константе предиката $P(x)$ во всех случаях, когда предикат $Q(y)$ равен единице, где $x$, $y$ - конкретные гейты. Элементами семейства одинаковых конструкций в данном случае будут гейты $x$ и $y$ при всех значениях seed-а порождающей нейросети, при которых $Q(y) = 1$. И на таких конструкциях мы форсируем равенство $P(x)$ некоторой константе.

Тем самым, симетрия и определенность - тоже представляются вещами одной природы. Человек, решая задачу, видя, что ему не удается наложить на свое решение ограничение в виде новой симметрии или определенности эффективным образом, часто начинает накладывать хотя-бы какую-то симметрию или определенность, которые не приносят никакой непосредственной пользы априори. И это, тем не менее, часто помогает решить задачу.

Возвращаясь к большой задаче о выписывании строки, можно представлять нашу систему следующим образом. Есть вышеописанный паралеллепипед, в котором есть гейты и лампочки, к которому, по нашему желанию могут добавляться "внешние"{} гейты и лампочки и есть какая-то порождающая модель, задающая распределение на значениях этих гейтов и лампочек. Мы добавим еще два типа объектов. Чтобы говорить о способе по которому гейт паралеллепипеда вычисляется по гейтам, расположенным непосредственно под ним, удобно говорить о "предписаниях"{} - специальных лампочках, которые "предписывают"{} гейту принимать определенное значение, когда некоторый набор из гейтов, расположенных непосредственно под ним принимает определенный набор значений. Это значит, что когда значение этого "верхнего"{} гейта не соответствует предписываемому ему значению "нижних"{} гейтов и эти нижние гейты принимают ровно те значения, которым соответствует данное предписание, а само предписание равно единице (активно), то возникает конфликт между значениями гейтов и предписанием. А именно, специально форсируем равенство нулю той лампочки, которая соответствует указанным значениям гейтов, конфликтующих с предписанием и значению 1 самого предписания. Другими словами, мы обеспечиваем правильное соотношение значений предписаний и гейтов и лампочек не-предписаний.

Кроме форсирования правильности этого соотношения, мы форсируем еще 2 вещи. Во-первых, это форсирование равенства нулю лампочек из верхнего слоя паралеллепипеда, соответствующих запрещаемых входной машиной Тьюринга наборам значений соседних гейтов. Во-вторых, это форсирование выполнения соотношений на предписания, обеспечивающих, что все предписания непротиворечивы и что каждой конфигурации гейтов, расположенных под каждым гейтом, соответствует хотя-бы одно предписание, предписывающее этому гейту принять некоторое значение. (Для последнего нужно построить маленькие схемы для всех групп предписаний, способных предписывать что-либо данному гейту, равные единице в случае, если для каждого выбора значений гейтов, расположенных под рассматриваемым гейтом, хотя-бы одно из этих предписаний действительно "активизируется"{}, форсировав соответствующие соотношения на соседние элементы этих схем - задав предписания для этих схем - и форсировав равенство этих схем единице.)

Лампочки-предписания ничем не отличались бы от обычных лампочек, если бы не одна вещь. А именно, еще мы требуем, чтобы порождающая модель в конце своего изменения выдавала для каждого из предписаний одно и тоже, постоянное значение - 0 или 1, независимое от входной строки. И мы изменяем наше распределение на значениях таким образом, чтобы его ограничение на предписания было как можно "уже"{}, как можно более сконцентрировано на одном конкретном наборе значений. Как добиться этого форсирования, пока не до конца понятно. Были идеи работать с оценкой величины, которую я называю узостью распределения. Узость распределения - это сумма по всем возможным наборам значений лампочек квадрата вероятности того, что распределение попадет именно в этот набор значений. Видно, что узость максимизируется, когда распределение сконцентрировано на одном конкретном наборе значений. Поэтому, чтобы добиться сконцентрированности распределения на каком-то одном наборе значений, целесообразно увеличивать узость. Но работать с узостью непосредственно вряд ли получится, поскольку всевозможных наборов значений много. Но были мысли о том, как ее можно оценивать.

Будет правильным обобщить понятие предписания до понятия "твердой"{} лампочки. Твердые лампочки - это лампочки, точно так же, как и предписания, накладывающие ограничения на конфигурации из не твердых лампочек, и для которых мы точно так же обеспечиваем правильное соотношение их значений и значений не твердых лампочек, форсируя обнуление соответствующих элементов точно так же, как указано чуть выше для случая предписаний. И для твердых лампочек мы точно так же форсируем увеличение узости распределения на их значениях - для нас так же важно, чтобы твердая лампочка в итоге либо всегда загоралась, либо всегда не загоралась (поэтому я их и назвал твердыми).

Например, если мы хотим "вырастить"{} признак на ситуации, "как на основании"{}, в стиле прошлой статьи, то с помощью твердых лампочек можно добиться жесткого вычисления гейта истории вычисления признака по гейтам, расположенным непосредственно "под"{} ним (слово "под"{} взято в кавычки, поскольку ориентация объектов в пространстве уже определена тем, как расположен наш базовый паралеллепипед). Так же, очень удобно форсировать с помощью твердых лампочек постоянные жесткие связи между значениями лампочек. 

Таким образом, предписания - это частный случай твердых лампочек.

Еще один тип объектов, которые мы вводим - это показатели. Каждый показатель принимает на вход ситуацию, "нарисованную"{} на ней конструкцию и вычисляет функцию из некоторого класса функций. Порождающая модель не генерирует значение показателя, показатель по входным данным вычисляется жестко. Этим классом функций может быть класс многочленов от значений лампочек ситуации и конструкции (например, многочлены степени не больше 7) или, скажем класс нейронных сетей от этих же значений. Много функций, которые хотелось бы выразить (в частности, коциклические многочлены), выражаются с помощью многочленов, но нейросети обладают кажется еще большей выразительностью, чем многочлены, поэтому пока для определенности будем считать, что указанный класс функций - это нейросети.

Алгоритм занимается тем, что он постоянно пытается изменить значения показателей так, чтобы как можно сильнее форсировать симметрию. Другими словами, он изменяет веса порождающей нейронной сети таким образом, чтобы показатели изменилясь в нужную сторону. В какую сторону надлежит изменяться конкретному показателю, определяется тем, в какую группу соответственных элементов семейств аналогичных конструкций этот показатель входит. 

Форсируя симметрию, мы одновременно решаем саму задачу и решаем близкую задачу. И мы переносим информацию с близкой задачи на нашу задачу (а может быть, и обратно): значения показателя конструкций (показатели зависят еще и от самих ситуаций, но мы это опускаем чтобы не загромождать текст, поэтому называем их показателями конструкций) в распределении для задачи приближаются к значениям показателя конструкций в распределении для близкой задачи. Если пытаться реализовать эту интуицию формально, возможно стоит уменьшать сумму квадратов разностей значения показателя конструкции для задачи и значения показателя конструкции для близкой задачи по всем наблюдаемым парам из конструкции в задаче и конструкции в близкой задаче. (Или, что то же самое, уменьшать сумму квадратов отклонений значений показателя в задаче от среднего арифметического значений показателя на близкой задаче.)

Если мы выделили семейство аналогичных конструкций для распределения на нашей задаче и семейство аналогичных им и друг другу конструкций для распределения на близкой задаче, и выделили для всех этих конструкций семейство соответственных элементов (мы выделили один показатель, соответственными элементами будут экземпляры этого показателя для разных конструкций), то данные два семейства аналогичных конструкций определяют, как изменяется каждый конкретный элемент для нашей задачи (не для близкой), определяют силу, с которым нам бы хотелось, чтобы он изменялся (которую мы так же будем назвать "аргументом"). Для каждого такого элемента, эта сила, так сказать, распространяется вниз, по backpropagation, к весам порождающей нейросети. Пораждающая сеть меняет свои веса, "учитывая пожелания"{} всех показателей, которые мы рассматриваем, причем такое пожелание учитывается тем сильнее, чем сильнее указанная выше сила.

В частности, порождающая нейросеть изменяется так, чтобы делать менее дискриминированными коциклические многочлены: если коциклический многочлен присутствует в списке показателей, которые мы учитываем, то информация об этом многочлене спускается вниз к весам порождающей нейросети. Но, кроме "вылечивания"{} коциклических многочленов, мы форсируем усиление и многих других симметрий.

В случае внешней нейросети, близкая задача будет совпадать с самой задачей и решаться она будет одновременно с самой задачей. Тем самым, значения показателя конструкций для задачи - это суть значения показателя конструкций для близкой задачи, таким образом, приближаться они будут друг к другу (иными словами, к их среднему). В частности, если большинство значений этих показателей близки к некоторой константе, а меньшинство - далеки, возникает некоторая сила, притягивающая значение для меньшинства к значению для большинства - к указанной константе.

Сделаю отступление, что на самом деле, не исключено, что чего-то можно добиться и без форсирования симметрии, скажем, точно так же генерируя некоторое распределение на значениях предикатов и предписаний дополненного параллелепипеда, увеличивая узость для распределения, ограниченного на предписания и, подобно тому, как мы уменьшаем loss function при обучении нейросети, уменьшая "функцию проблемы": функция проблемы - это просто суммарное количество, опять же, во-первых, нарушений, несоответствий между значениями предписаний и значениями окружающих их гейтов, во-вторых нарушений запретов в верхнем слое паралеллепипеда, заданных условием задачи (включая запрет ответу схемы быть нулем), в-третьих, нарушений запретов на соседние элементы упомянутых выше схем, проверяющих, что хотя-бы одно предписание активизировалось (включая запрет ответу такой схемы быть нулем) и, в-четвертых, конфликтующих пар предписаний. Поскольку, когда порождающая нейросеть сработала, значения лампочек еще не определены, для них есть лишь вероятности (с которыми, сразу после этого, сами бинарные значения выбираются независимо), то можно говорить лишь об ожидаемой функции проблемы. И эту ожидаемую функцию проблемы можно спустить вниз по backpropagation, по слоям порождающей нейросети и изменять веса последней так, чтобы она ожидаемо уменьшалась. Не исключаю, что и этим можно чего-то добиться в сочетании с какими-то другими идеями.

Та сила, с которой мы изменяем показатели, соответствующая некоторому семейству аналогичных конструкций на задаче и на близкой задаче тем больше, чем в большей степени данное семейство аналогичных конструкций является "непринужденно найденным". О том, что значит "непринужденно найденный"{} и как это понятие определять, я буду говорить ниже, в других секциях.

Сразу скажу, что причина, по которой перед тем, как усиливать корреляцию в ответе, выдаваемом внешней нейронной сетью, для нас важно, чтобы ответ уже коррелировал, как раз в том, что конструкции, соответствующие внешней нейросети будут считаться найденными непринужденно именно в этом случае.

Разумеется, для одного и того же семейства аналогичных конструкций одни показатели более предпочтительны для рассмотрения, чем другие. Скажем, предпочтительны простые показатели. Например, если конструкции семейства представляют из себя просто упорядоченные наборы из нескольких гейтов и, таким образом, между парами конструкций есть естественное соответствие и если это соответствие в значительной степени "уважает"{} много разных связей в поднаборах этих гейтов, что и делает эти конструкции аналогичными, то одним из наиболее простых является показатель, просто выдающий значение некоторого предиката в $i$-м по счету гейте. Такой показатель идеально соответствует интуиции "соответственные элементы семейства конструкций".

Возможно, опять же, с большей силой стоит изменять те семейства соответственных элементов, для которых уже прослеживается корреляция со значениями аналогичных элементов на близкой задаче.

Превращение хаотичного набора гейтов в структурированную однородную решетку (например, задающую признак, в смысле прошлой статьи) тоже может стать результатом форсирования симметрии. Это экземплифицирует тот факт, что форсирование симметрии можно воспринимать как уменьшение сложности.

Проявлением форсирования симметрии в мышлении человека, помимо всего прочего, можно считать применение теорем в математике. Что есть применение теоремы при решении математической задачи? Это суть нахождение в новой ситуации старой, уже известной конструкции и переносом некоторых данных с других ситуаций - тех, в которых данная конструкция так же присутствует, на эту ситуацию.

Позволю себе немножко порассуждать о нашей физической реальности, высказав даже не гипотезу, а допущение, которое имеет некоторую вероятность оказаться верным, пусть и не большую. Наша вселенная обладает некоторым простым свойством, которое я называю экстраполируемостью. А именно, во вселенной можно часто пронаблюдать семейства объектов (предметов, моментов времени, ситуаций, событий, элементов человеческой психики, и т.д.) с одинаковыми или похожими свойствами. Более того, если мы наблюдаем семейство объектов, которые, при этом, обладают неким объединяющим их свойством $P$, и при этом выясняется, что все наблюдаемые объекты обладают еще и свойством $Q$, то можно законно предположить, что и следующий встретившийся нам объект, обладающий свойством $P$ будет обладать еще и свойством $Q$. Но почему это так? И вообще, почему подобные семейства объектов встречаются нам столь часто? Например, если взять длинную строку, в каждой ячейке которой, независимо от остальных ячеек, генерируется символ 0 или 1, с вероятностью $1/2$, то такой объект похожим свойством обладать не будет: если мы выберем случайно набор ячеек этой строки (более общо - набор значений признака, в смысле, определенным в прошлом тексте, на этой строке и некоторых унарных аргументах; здесь нужно еще потребовать, что признаку не должно быть свойственно слишком часто обращаться в 1 либо слишком часто обращаться в 0) и символы в них (ячейках, унарных аргументах, соответственно) окажутся равными 1, то это совершенно не говорит ничего о том, какой символ будет записан в следующей случайно выбранной нами ячейке (унарном аргументе). Далеко не все распределения экстраполируемы. (Хотя распределения, рассматриваемые людьми - возникающими при решении задач реального мира - конечно, часто бывают.) Кстати, интересно строго определить и поизучать экстраполируемое распределение, возможно даже стоит пытаться обучаться экстраполируемому распределению (возможно, это поможет).

Возникает вопрос, что же сделало наш мир экстраполируемым? Можно попытаться объяснить очень богатую экстраполируемость во вселенной и, в частности, в жизни на земле, изотропностью и однородностью пространства и времени и выполнением некоторых, в некоторой степени известных нам базовых законов физики, по очереди объясняя каждую следующую закономерность из предыдущих, более "низкоуровневых". Но если провести мысленный эксперимент, мысленно создав вселенную в которой выполняется этот набор уравнений на движение частиц - как-то ее инициировать и запустить движение частиц, так, чтобы уравнения выполнялись - и рассмотреть одну из сформировавшихся планет, на которой образовалась некоторая форма жизни (если такая планета вообще будет), то будет ли эта жизнь столь же экстраполируема как наша или все будет гораздо более хаотично?

Что если экстраполируемость сама непосредственно форсируется некоторым физическим законом? Думаю, такой закон можно было бы интерпретировать по-разному - форсирование экстраполируемости, форсирование выполнения бритвы Оккама, форсирование уменьшения сложности, наконец, форсирование симметрии. (Интуиция подсказывает, что форсируя одно, форсироваться будет и все остальное из списка, поскольку эти названия говорят примерно об одном и том же.) Я думаю, что такое не исключено, что симметрия во вселенной форсируется каким-то статистическим способом, в чем то похожем на то, что происходит в описываемом заготовке алгоритма (это предложение более внятно выражает это гипотетическое свойство вселенной, чем моя пока не удавшаяся попытка строго определить экстраполируемость и мое предположение о ее форсировании).

В том числе, форсироваться могут и уже известные нам физические законы и уравнения. Каждому такому закону может соответствовать свой признак. Речь о предполагаемо существующих аналогах определенных в тексте признаков для непрерывного мира. На самом деле, думаю, есть много способов определять признаки для непрерывного мира. Например, аналогом правильного предиката может служить точка, выбранная в пространстве или момент, выбранный во времени. А чтобы определить непрерывный признак, может быть можно как-нибудь определить непрерывную машину Тьюринга - аналог дискретной. У меня были мысли о том, как ее можно определять, но сейчас речь о другом. Могу сказать ключевую идею, сделав тем самым отступление. У нас есть пространство, на котором что-то происходит - задан некоторый набор гладких функций, так же у нас есть несколько головок - точек на этом пространстве. Мы располагаем некоторым множеством евклидовых пространств некоторых размерностей, на которых расположены точки, которые мы называем состояниями, которые, перемещаясь, на этих евклидовых пространствах что-то "рисуют"{} (можно сказать, что каждая точка, перемещаясь, рисует свою кривую непрерывно меняющимся "цветом"{} и множество "цветов"{}, которыми она что-то рисуют, так же образует евклидово пространство). Каждое состояние (и головка) гладким образом изменяется по некоторому закону: его производная - это некоторая фиксированная гладкая функция самого этого состояния (положения головки), некоторых других состояний и их производных (возможно, с производными больших порядков тоже стоит работать) и того, что на некоторых из евклидовых пространств нарисовано (когда есть зависимость от нарисованного, гладкость данной функции, по сути, являющейся функцией в том числе и от гладких кривых, следует определить естественным образом). То, что рисует точка-состояние на своем евклидовом пространстве в данный момент определяется аналогично. С помощью таких непрерывных машин Тьюринга можно пытаться записывать утверждения типа "точка находится внутри футбольного мяча"{} или "точка в данный момент времени находится внутри человека, который сейчас нагнется, возьмет палку и кинет ее вперед". (Последнее - для случая, когда мы работаем с пространственно-временным континуумом, а не просто с трехмерным пространством.)

Так вот, речь о том, что, скажем, некоторый физический закон может быть равносилен тождественному обращению в 0 некоторого непрерывного признака, на нашем мире и любых значениях входящих в его определение аналогов правильных предикатов. Если симметрия действительно форсируется, то вполне может оказаться, что этот непрерывный признак все же не тождественный нуль, но очень близок к тому, чтобы быть тождественным нулем и описываемый гипотетический закон форсирования симметрии изо всех сил приближает его к тому, чтобы он был тождественным нулем. 

Но сделать его строго равным нулю тождественно не выходит, поскольку есть и другие нарождающиеся симметрии, которые нужно тоже форсировать и которые могут быть далеки до того, чтобы войти в учебник физики на правах физического закона, например, закон "каждое утро такой-то голубь садится на крышу такого-то сарая"{} тоже может форсироваться и в некоторой мере заставлять голубя каждое утро садиться на этот сарай снова и снова (ожидаемо, в маленькой мере: ожидать того, что голубь на следующее утро опять сядет на крышу сарая следует не из соображения аргумента, соответствующего указанному ряду похожих конструкций и, соответственно, показателей, а скорее из того, что мир под многовековым влиянием форсирования симметрии стал сильно экстраполируемым - симметрия влечет другую симметрию: если в объекте удалось выделить большое количество параметризованных семейств одинаковых конструкций, то высока вероятность, что в нем найдется и много других параметризованных семейств одинаковых конструкций). И эти другие нарождающиеся симметрии могут давать свой вклад в общий процесс изменения мира, пусть иногда и очень слабый. Другими словами, грубо говоря, физические уравнения могут выполняться не вполне строго, оставляя происходящему в мире некоторую свободу и эта свобода может использоваться для форсирования других симметрий, в том числе, окружающих нас в повседневной жизни (другие симметрии немного "тянут на себя"). Не исключено, что этой идеей можно объяснить некоторые жизненные закономерности.

Если все это так, то не исключено, что некоторые, скорее всего очень "масштабные"{}, вещи, происходящие в мире могут "сломать"{} физические законы - с точки зрения форсирования симметрии и одного из аспектов ее идеологии, который можно описать поговоркой "где тонко, там и рвется"{}, который будет обсуждаться позже, если какой-то физический закон станет этим "тонким местом"{}, оставлять закон выполняющимся станет нецелесообразно и он перестанет выполняться, впоследствии возможно заменившись другими законами.

Мы видели, что форсирование симметрии - это частный случай устранения противоречия, устранения конфликта близких задач, так что не исключено, что правильней скорее предположить, что в природе форсируется устранение противоречия.

Даже если в мире действительно что-то похожее происходит, то пока эта идея далека до реализации и проверки на экспериментах, просто мне кажется, что такую возможность интересно иметь в виду. 

Не исключено, что и в мышлении человека постоянно форсируется симметрия. Человек очень любит определенность, которая, как мы видели, тоже является некоторой формой симметрии. Это проявляется в том, как он любит вносить определенность в моральные установки (что-то хорошо всегда/что-то плохо всегда), определенность в принимаемых законах ("закон един для всех"{}, "закон применим во всех таких ситуациях"{} и т.д.), само понятие информации включает в себя некоторую определенность. Человеку очень удобно на что-то опереться, подобно тому, как коциклический многочлен "опирается"{} на защищенную лампочку. 

Еще по этому поводу можно указать на некоторую связь, которая, конечно, может быть просто совпадением. То форсирование еще большей корреляции объектов, которые и так коррелируют, что мы видели выше (внешняя нейронная сеть, трактовка коциклического многочлена в аспекте форсирования симметрии - на решениях близкой задачи коциклический многочлен коррелирует с собой максимально сильно), может иметь отношение к правилу Хебба для нервной системы: нейроны, загорающиеся часто одновременно, усиливают связь между собой и, тем самым, начинают загораться одновременно еще чаще.

Закончу эту секцию тем, что по моим наблюдениям, человек начинает искать причину всегда либо у наличия симметрии, либо у нарушения симметрии (в нашей интерпретации понятия симметрии), это указывает на возможную связь с причинностью. Есть так же не проверенная гипотеза о том, что сама по себе причина нарушения симметрии является глубокой проблемой, а причина симметрии является зоной согласованности (об обоих объектах речь будет идти ниже; называя причину, человек как-бы задает зону согласованности или глубокую проблему). Но в последнем с большой вероятностью могу ошибаться.

\section{Обратный поиск.}

Отдельно хочется выделить такую возможность. Наша цель - выбрать предписания в паралеллепипеде таким образом, чтобы на верхнем его слое выполнялись некоторые соотношения между значениями гейтов. Что если попробовать "защитить"{} эти самые соотношения, а именно, для каждого из соотношений соорудить коциклический многочлен, доказывающий, что это соотношение выполняется, при этом, "воспользовавшись"{} тем, что выполняются некоторые предписания или, более общо, загораются некоторые твердые лампочки? Имеется в виду следующее. Изначально мы, конечно, не знаем, какие из твердых лампочек у нас загорятся, а какие - нет. Для каждой из лампочек мы будем поддерживать вещественнозначную величину, - "намерение"{}, отображающую, в какой степени мы намерены защитить эту лампочку. (Можно понимать слово "защитить"{} чуть более общо: не только добиться равенства лампочки нулю, но и добиться равенства лампочки единице, - в первом случае, за фразой "воспользоваться этой лампочкой"{} кроется то, что некоторый коциклический многочлен имеет один из своих положительных коэффициентов в самой лампочке, а во втором случае, эта фраза подразумевает, что коциклический многочлен имеет своим положительным коэффициентом лампочку - отрицание указанной лампочки. Если мы лампочкой воспользовались, это значит, что мы объявили ее значение, в первом случае, нулем, а во втором случае - единицей. Разумеется, воспользоваться лампочкой можно далеко не всегда.) 

Изначально мы намерены защитить лишь некоторые лампочки из верхнего слоя паралеллепипеда. Еще один вещественнозначный параметр, который мы будем поддерживать - это "пригодность"{} или "пригодность для защиты". Абсолютная пригодность для защиты изначально есть только у предписаний паралеллепипеда и некоторых твердых лампочек (грубо говоря, если мы пронумеруем гейты вне паралеллепипеда натуральными числами, то мы можем потребовать, чтобы загорались твердые лампочки, предписывающие гейтам с большими номерами вычисляться по гейтам с меньшими номерами, если эти предписания не могут противоречить друг другу, - по значениям гейтов паралеллепипеда мы всегда восстановим значения и этой "внешней"{} последовательности гейтов непротиворечивым образом, но когда речь заходит о наложении каких-то других дополнительных ограничений, связывающих эти гейты, не всегда ясно, почему мы вправе требовать этих ограничений априори). Причем пригодность для защиты предписания полностью пропадает, когда мы решаем защитить другое, противоречащее ему предписание.

Предположим, мы намерены защитить некоторую лампочку. Чтобы полностью защитить эту лампочку, нам достаточно найти коциклический многочлен, доказывающий, что эта лампочка принимает нужное значение. При этом у этого многочлена будет несколько положительных коэффициентов и мы заинтересованы, чтобы эти коэффициенты стояли на нулевых лампочках. Если такая лампочка с положительным коэффициентом, скажем, является отрицанием предписания, такого, что никакое из противоречащих ему предписаний защищено не было, мы имеем право защитить эту лампочку сразу. Но в общем случае, чтобы защитить лампочку, которой мы воспользовались, нужно приложить некоторые усилия - сконструировать другой коциклический многочлен, доказывающий, что эта лампочка действительно принимает нужное значение. Сконструировав его, мы опять чем-то воспользуемся. Возможно, то, чем мы пользуемся, удастся все защитить сразу, а возможно, - нет и те лампочки, которые не удается защитить сразу, нужно, опять же, защищать с помощью своих коциклических многочленов.

Если нам удастся защитить все лампочки из верхнего слоя паралеллепипеда, которые нужно защитить, доказательства равенства их значений должным числам будут представлять из себя деревья из коциклических многочленов, почти то же самое, что "обобщенная клетка с подвешиваниями и дополнительным построением"{}, о которой шла речь в начале прошлого препринта. Коциклические многочлены этих деревьев могут и будут пользоваться выполнением некоторых предписаний или других вещей, которыми они "вправе"{} воспользоваться - изначально мы не знаем, каковым будет значение предписания, мы это понимаем лишь когда мы решили воспользоваться равенством предписания тому или иному значению, ради "пропускания"{} коциклического многочлена. И если предписания, зажигание которых предполагают построенные доказательства и правда не противоречат друг другу, то задача решена.

Строя коциклический многочлен в надежде что-то защитить, следует выбирать его так, чтобы, во-первых, он воспользовался как можно меньшим числом лампочек, а во-вторых, чтобы те лампочки, которыми он все-таки воспользовался, были как можно более пригодны для защиты. Возможно, стоит искать такой многочлен поиском в пространстве коциклических многочленов подобным тому, что использовался в прошлом препринте, это сводит к некоторой оптимизационной задаче.

Пригодность лампочки для защиты следует определять, основываясь на мысленных попытках ее защитить. Например, если существует коциклический многочлен, доказывающий, что лампочка примет нужное значение, и который при этом воспользуется всего двумя значениями других лампочек, одно из которых можно защитить легко, а второе можно вообще защитить сразу "бесплатно"{}, то скорее всего исходную лампочку тоже будет защитить легко и ее стоит рассматривать как пригодную для защиты. Если же со временем ситуация изменится и одна из двух лампочек станет менее пригодной для защиты, то и исходная лампочка может стать менее пригодной для защиты.

Пригодность лампочки для защиты зависит от пригодностей других лампочек и от намерений (о намерениях я писал чуть выше, зависимость пригодностей от намерений скоро станет понятна). Намерения же, в свою очередь, могут зависеть от пригодностей (нам вряд ли стоит быть намеренными защищать лампочку, непригодную для защиты). Намерения и пригодности зависят друг от друга и изменяются в едином динамическом процессе. 

В этот же процесс стоит включить, собственно, и распределение значений на самих лампочках, с которым мы работаем вне этой подсекции. Например, предположим, что мы твердо намерены защитить несколько лампочек. И перед нами встает выбор, защитить ли еще одну лампочку, отличную от этих лампочек. И мы видим, что если нам и правда удастся защитить все эти лампочки, включая новую, то, с учетом того, что все эти лампочки примут соответствующие значения, у нас образуется очень напряженный коциклический многочлен. Это будет означать, что намереваться защищать эту новую лампочку, возможно, не стоит (вот нам и пример зависимости пригодности от намерений, - намерение защитить те лампочки, о которых идет речь, кроме новой, делают новую лампочку непригодной для защиты). 

Скажу отдельно еще и такую мысль, что при построении коциклического многочлена, помимо того, чтобы значения лампочек, которыми мы пользуемся, были как можно более пригодны для защиты и их было как можно меньше, очень приветствуется, чтобы эти значения лампочек как можно сильнее коррелировали друг с другом. Это нужно для того, чтобы заставляя эти лампочки принимать такие значения, мы лишались примерно одних и тех же свобод, о чем речь - станет понятно дальше.

Возможно, помимо таких величин, как пригодность и намерение, так же будет полезна такая величина, как "желаемость". Желаемыми мы называем такие лампочки, защита которых была бы очень полезна для нас. При этом, намерение защитить желаемые лампочки может быть слабым, например, потому, что мы можем понимать, что защитить такую лампочку будет очень сложно. Это еще одна величина, которая, наряду с другими, может быть включена в процесс динамического изменения, в котором все эти величины влияют друг на друга. Возможно, все эти величины стоит реализовывать как распределения, а не как просто набор вещественных величин, по одной для каждого гейта.

\section{Абстракции и непринужденность.}

Напомню, в большой задаче о выписывании строки у нас есть схема, приходящая из машины Тьюринга, принимающая 2 строки $A$ и $X$. Наша цель - научиться по строке $A$ строить строку $X$, чтобы $M(A,X) = 1$. Мы рассматриваем прямоугольный паралеллепипед в трехмерной решетке, на узлах которого определены несколько предикатов. Верхнюю из сторон паралеллепипеда мы интерпретируем как историю работы $M$. Весь паралеллепипед мы интерпретируем как схему, вычисляющую историю работы $M$ снизу вверх (каждое значение предиката гейта вычисляется по значениям предикатов в нескольких гейтах непосредственно под данным гейтом). В качестве выхода мы хотим получить корректную историю работы схемы $M$, которая в итоге выдает единицу. 

Точно так же можно ставить себе задачу так. Когда мы имеем дело с решеточным паралеллепипедом, мы говорим о произведении $[0,...,m] \times [0,...,k] \times [0,...,n]$, где $n$ - высота паралеллепипеда, а $m$, $k$ - параметры схемы. Но мы можем захотеть видеть в верхнем сечении схемы, помимо гейтов схемы какой-то дополнительный набор признаков $S$ с теми выполненными пропозициональными связями между ними, которые должны выполняться. Таким образом, мы имеем дело со схемой, представляющей из себя произведение $S \times [0,...,n]$, нижний уровень которой заполнен нулями всюду кроме входной строки $A$ и мы хотим выбрать правило вычисления гейтов в схеме так, чтобы она выдавала набор признаков, связанных между собой нужными пропозициональными связями и, конечно же, чтобы признак, соответствующий ответу, выдаваемому $M$ был равен единице. Добавим схему $R$, вычисляющую свое значение по верхнему слою паралеллепипеда и выдающую 1 как раз в этом случае (все признаки верхнего слоя связаны нужными пропозициональными связями и ответ, выдаваемый $M$, равен единице) и только в этом случае.

Таким образом, наша цель в том, чтобы $R$ выдала 1.

В общем случае рассматриваемые схемы, "соединяющие"{} вход с выходом, конечно же, не обязаны быть чем-то напоминающим прямоугольную решетку или вот так разбиваться на $n+1$ одинаковый слой.

Предположим, что нам удалось построить такую схему $Q_i$ в ряде случаев $D_1$, $D_2$, ..., $D_r$ и каждый из этих случаев определяется некоторым набором признаков $P_i$ строки $A$, обращающихся в единицу.

Мы сейчас захотим построить некоторый объект и от построения этого объекта мы захотим некоторого свойства. Объектом является схема $Q$ с аналогичным входом и выходом, а так же набор коциклических многочленов $Z_i$ - доказательств того, что в случае $D_i$, $R(Q) = R(Q_i)$. 

Остановимся на том, что значит "построить $Z_i$". Создадим новый набор гейтов и ограничений. К каждому набору гейтов, образующему лампочку прикреплен набор гейтов, олицетворяющих числа от $-N$ до $N$, где $N$ - достаточно большое натуральное число. Если мы работаем над $\mathbb{Q}$, то все коэффициенты коциклического многочлена можно привести к общему знаменателю, тем самым, сделав их целыми. Так что, чтобы задать $Z_i$, можно при каждой рассматриваемой лампочке сопоставить единицу гейту, олицетворяющему число - коэффициент многочлена при данной лампочке после приведения к общему знаменателю и сопоставить ноль всем остальным гейтам, олицетворяющим числа. Чтобы проверять, что такое сопоставление - это то, что нам нужно, следует добавить и проследить выполнение признаков, проверяющих, во-первых, что при каждой лампочке есть ровно одна единица, во-вторых, что незащищенным лампочкам сопоставлены неположительные числа, в-третьих инвариантность значения многочлена при заменах значений гейтов с нуля на единицу и, в-четвертых, положительность этого значения на произвольном наборе значений гейтов. Таким образом, построение доказательства ничем не отличается от построения выполняющего набора некоторой задачи SAT. Если не нравится то, что $N$ нужно выбирать достаточно большим, можно вспомнить о том, что у нас есть еще и вещественнозначные гейты и применяя аналогичный рассматриваемому в данном тексте способ заставить выполненяться указаные ограничения для вещественнозначного случая. Последний вариант мне нравится больше, но в дальнейшем мы пока будем подразумевать первый вариант, поскольку в вещественнозначном случае есть некоторые нюансы, если быть точнее, забегая вперед, штрафовать нужно не только за увеличение узости, далее будет понятно о чем речь.

Итак, для каждого $i$ мы строим доказательство $Z_i$ того, что в случае $D_i$, $R(Q) = R(Q_i)$, а так же, саму схему $Q$. Мы предположили, что в случае $D_i$, $Q_i$ выдает то, что нужно. Возможно, у нас в руках даже есть доказательства того, что $Q_i$ выдает то, что нужно, в случае $D_i$, тем самым, для каждого $i$, всю схему $Q_i$ стоит добавить в рассмотрение при поиске нужных нам объектов, поскольку в ней содержится много лампочек, которые можно защитить, а значит, ими можно будет воспользоваться. Имеется в виду, что эта схема будет актуальной для поиска $Z_i$. Заметим, что $Z_i$ может использовать, как защищенные, признаки $P_i$. По желанию, конечно, можно добавить и другие признаки, помимо перечисленных.

Вначале у нас есть только "остов"{} схемы $Q$, набор гейтов, без каких-либо запрещенных конструкций - предписаний, как $Q$ должна вычислять следующие гейты по предыдущим. Найти эти предписания - наша задача. Когда мы строим $Z_i$, мы вынуждены создавать такие предписания: когда мы хотим поставить положительное значение коэффициента $Z_i$ на некоторую лампочку, мы должны сделать такую лампочку защищенной, то есть, потребовать выполнения некоторого предписания (или потребовать защищености какой-то другой лампочки, не относящейся к схеме $Q$, впоследствии защитив эту лампочку). Если мы хотим построить $Z_i$, мы должны делать это так, чтобы они были "согласованы"{} на $Q$, а именно, чтобы среди индуцированных таким образом предписаний $Q$ не было двух противоречащих друг другу - предписывающих, в некотором случае, одному и тому же гейту быть равным нулю и единице соответственно. Если мы это сделаем, схема $Q$ будет построена, так как она определяется своими предписаниями. 

Можно на это смотреть так: мы для каждого из случаев $D_i$ построили доказательство того, что схема, которую мы строим, будет работать правильно в случае $D_i$, причем постарались, чтобы доказательства были согласованы друг с другом. После этого, сами доказательства можно выкинуть, а предписания, задающие нужную нам схему останутся. И построенная схема будет работать уж точно на объединении всех случаев $D_i$. Но вот вопрос - будет ли она работать и в остальных случаях?

Мое предположение в том, что эта схема будет работать и в остальных случаях, если мы строим $Q$ и $Z_i$ правильным образом, как описано ниже.

Если в двух словах, нам нужно стараться как можно сильнее форсировать симметрию, при этом как можно медленнее увеличивая узость. 

Значения предписаний для схемы $Q$ и значения элементов $Z_i$ будут сформированы путем изменения совместного распределения всех этих значений, генерируемого нейронной сетью (опять же, либо какой-то другой порождающей моделью), узость которого мы заставляем в итоге быть максимально возможной. Такой совместный набор значений мы так же будем называть ситуацией.

Как можно догадаться из позапрошлого абзаца, нам понадобится величина $\hat{S}$, являющаяся оценкой другой величины $S$, которая, в свою очередь, олицетворяет количество симметрии, которая в данный момент присутствует в распределении (величина, похожая на величину, обратную колмогоровской сложности). 

Пока я не определяю величину $S$ и не пишу, о том, как ее оценивать, вычисляя $\hat{S}$. Попытка определить и оценить $S$ и связанные с этим объекты и построения, а так же их дальнейшее использование - отдельная тема, которой я собираюсь заниматься дальше, какие-то соображения на этот счет можно прочесть в секции "Еще одна идея".

Скажу лишь то, что $S$ является мерой того, насколько много в текущем распределении семейств похожих конструкций. Причем для каждого из этих семейств важно то, как именно эти конструкции расположены относительно друг друга. Нужно, чтобы они были близки к тому, чтобы быть "параметризованными"{}, как в это происходит в случае с внешней нейронной сетью или коциклическим многочленом. Об этой "параметризации"{} конструкций я поговорю позже.

Кстати, если в итоге все действительно получится правильно определить, то это открывает дорогу для определения колмогоровской сложности для непрерывных объектов (скажем, колмогоровская сложность многообразия).

Итак, предположим, что мы поддерживаем некоторую величину $\hat{S}$, оценивающую количество симметрии в распределении, которую мы хотим сделать как можно больше. 

Казалось бы, - можно просто двигаться наверх по градиенту $\hat{S}$. Но есть что-то еще. Дело в том, что распределения с одинаковым $\hat{S}$, но разной узостью для нас не равноправны. Давайте для удобства говорить о "широте"{} распределения - величине, противоположной узости, как по смыслу, так и по значению ("широта"{} = -"узость"). Я утверждаю, что иметь в руках распределение большей широты - ожидаемо лучше, чем иметь в руках распределение меньшей широты. Поэтому я предлагаю, кроме того, чтобы "награждать"{} за повышение оценки симметрии распределения, так же и "штрафовать"{} за уменьшение широты.

Что будет, если не штрафовать за уменьшение широты? Заметим, что в общем случае, повышать симметрию за счет уменьшения широты гораздо легче, чем за счет ее увеличения (скажем, в распределениях наибольшей возможной широты вообще симметрии очень мало). И наш алгоритм будет склонен к тому, чтобы скорее уменьшать широту, чем ее увеличивать. Не накладывая штрафа за уменьшение широты мы, по сути, разрешаем бесполезные уменьшения широты. Чем это чревато?

Мы хотим, чтобы в конечном итоге наше распределение сузилось до распределения, выдающего правильные решения: ситуации без противоречия, и вообще, с высокой симметрией (напомню, отсутствие противоречия - правильность значений, выдаваемых коциклическими многочленами - это тоже проявление симметричности). Этих правильных решений в общем случае вообще не так много: предполагаю, что частой является ситуация, когда непротиворечивых ситуаций высокой симметрии не больше 20, а остальные непротиворечивые ситуации уже несколько проигрывают им в симметричности. Так вот, совершая бесполезные уменьшения широты, мы теряем эти самые правильные решения. И если за уменьшение широты не штрафовать, мы довольно быстро их все растеряем - вероятность того, что распределение выдаст правильное решение станет очень низкой и нам придется вновь расширять распределение, чтобы эту вероятность восстановить.

Вторая причина в том, что даже если ограничение на ситуацию, задаваемое бесполезным сужением, не приведет к потере всех правильных решений и мы потом сможем сузить результирующее распределение до искомого, то на этом пути нам будет мешать это "искажение"{}, вносимое в $\hat{S}$ ненужным сужением.

Итак, имеет смысл штрафовать за уменьшение широты. Естественным представляется решение, в котором мы движемся по градиенту суммы $\hat{S}$ и функции распределения $\rho$, являющейся монотонно возрастающей функцией широты. (О деталях мы пока не говорим, поэтому здесь, опять же, конкретную функцию не указываем.)

Иными словами, мы увеличиваем симметрию, как можно медленнее увеличивая узость. Для краткости, величину 
$\hat{S} + \rho$ мы будем называть дополненной симметричностью. (Вообще, напряженность коциклических многочленов и, более общо, маленькая симметричность - это, так сказать, потенциальное уменьшение широты: нам придется уменьшить широту, чтобы избавиться от напряженного коциклического многочлена.)

То, что я написал, как-бы подразумевает, что мы поддерживаем величину $\hat{S} + \rho$ и изменяем веса порождающей нейронной сети, используя backpropagation, таким образом, чтобы $\hat{S} + \rho$ увеличивалась как можно сильнее. Но на самом деле, мы можем использовать разные уловки - "действия"{} над распределением, которые ожидаемо приведут к увеличению дополненной симметричности не непосредственно сейчас, а впоследствии. 

Например, рассмотрим такую ситуацию. Есть гейт $u$, про который мы знаем, что если форсировать его равенство единице, то это приведет к трем приятным последствиям: гейты $t_1$, $t_2$, $t_3$ станут равны значениям $a_1$, $a_2$, $a_3$ соответственно. И известно, что это хорошо скажется на симметричности распределения. Поясню оба момента. Первая вещь может означать, что у нас имеются коциклические многочлены, доказывающие - или хотя-бы дающие аргументы, указывающие на то, что если $u$ примет значение 1, то $t_i$ примет значение $a_i$. Второе может означать, что у нас есть доказательства или хотя-бы аргументы, указывающие на то, что если мы хотим, чтобы в ситуациях, даваемых распределением выполнялись те, целевые, ограничения, которые в них должны выполняться, то $t_i$ примет значение $a_i$.

Другими словами, если мы форсируем равенство $u$ единице, то $t_1$, $t_2$ и $t_3$ ожидаемо вскоре будут реже принимать те значения, которые они не должны принимать. А это и означает позитивный эффект для симметричности: скажем, в ситуациях, в которых $t_1$ принимает значение, отличное от $a_1$, соответствующий коциклический многочлен (дающий аргумент, указывающий на то, что так не должно быть), будет более напряжен; и выкидывая большое число ситуаций, в которых коциклический многочлен более напряжен, мы ожидаемо повышаем симметричность. (Опять же, согласованность - это хорошо.)

Мы выяснили, что форсировать равенство $u$ единице - это хорошо. Предположим теперь, что есть еще гейт $v$, для которого точно так же верно, что форсирование равенства его значения единице приведет к трем приятным последствиям: гейты $t_4$, $t_5$, $t_6$ примут значения $a_4$, $a_5$, $a_6$ и есть аргументы, указывающие как на то, что они их таки примут, так и на то, что так и должно быть. Тем самым, мы так же знаем, что форсирование равенства $v$ единице - полезная вещь.

Предположим так же, что у нас есть гейт $w$, для которого аналогичным образом верно, что если мы форсируем его равенство нулю, то $t_i$ станет равным $a_i$ для всех $i = 1, 2, ..., 6$. Так же предположим, что в указанных трех случаях - форсирования $u$, $v$ или $w$ соответствующему значению, мы лишаемся свободы (широты) примерно в одинаковой степени: во всех трех случаях широта уменьшается примерно на одинаковую величину.

Сравним теперь две возможности: первая состоит в том, что мы форсируем одновременно равенство значений $u$ и $v$ единице, а вторая - в том, что мы форсируем равенство значения $w$ нулю. Как с чисто интуитивной точки зрения, так и в соответствии с логикой алгоритма, вторая возможность представляется ожидаемо более желаемой.

Что нам известно - это то, что как в первом, так и во втором случае наступает эффект в виде 6 одних и тех же приятных последствий. Так же нам известно, что в первом случае мы лишаемся широты (или свободы) ожидаемо сильнее, чем во втором (априори, когда про $u$, $v$, $w$ мы больше ничего не знаем). Первое говорит о том, что эффект, произведенный на симетричность в обоих случаях ожидаемо одинаков (на этот раз, за фразой "ожидаемо одинаков"{} кроется скорее равенство матожиданий, чем примерное равенство). Второе говорит о том, что в первом случае мы потеряем ожидаемо больше широты, чем во втором. Это все и объясняет предпочтительность второго варианта.

Если, скажем, обе операции представляют из себя выбрасывание из распределения части ситуаций и распределение вероятностей между остальными ситуациями в тех же пропорциях, то нам известно, что мы, судя по всему выбрасываем, помимо всего прочего, в ожидаемо примерно одинаковых количествах в этих двух случаях, ситуации, на которых отклоняются от нормы указанные 6 коциклических многочленов (и которые, кстати говоря, по этой причине не будут "правильным решением"); но общее количество выброшенных ситуаций в первом случае ожидаемо больше, чем во втором. Результирующее распределение во втором случае ожидаемо шире, при ожидаемо той же симметричности (вся информация об изменении симметричности приходит из 6 указанных приятных последствий).

Сейчас я экземплифицировал свою идею, которая несколько раз, в разных или похожих контекстах, у меня появилась при работе над алгоритмом, поэтому ее я считаю важной и решил посвятить ей следующую подсекцию. Можно ее сформулировать двумя способами.

Во-первых, можно так сказать: следует действовать так, чтобы при фиксированном эффекте, лишаться свободы наименее сильно, а при фиксированной потере свободы обретать наибольший эффект.

Вторая формулировка такая же, но она оперирует не словом "свобода"{}, а словом "усилие". В самом деле, в вышеописанном примере форсирование равенства значения гейта константе можно воспринимать как усилие. И в первом из описанных в примере случаев мы применяем ожидаемо больше усилий, чем во втором. И правило звучит так: следует действовать так, чтобы при фиксированных усилиях добиваться максимального эффекта, а для достижения фиксированного эффекта прилагать минимальное число усилий.

Скажу еще о том, что описанный алгоритм с оптимизацией величины $\hat{S} + \rho$ следует применять и в случае общего алгоритма (нахождения предписаний схемы-паралеллепипеда) - не только нахождения $Q$ и $Z_i$. Добавление слагаемого $\rho$ в оптимизируемую функцию и соображение о том, что стоит делать шаги, приводящие к улучшению в более отдаленном будущем, чем непосредственно сейчас, - все это применимо и в случае общего алгоритма.

Можно проиллюстрировать указанную идею на старой версии общего алгоритма, когда наше распределение на лампочках было независимым (намерение выбирать независимое распределение уже оставлено, но идею проиллюстрировать можно). В прошлом препринте я писал, что когда вероятность единицы при лампочке была равна, скажем, $1/2$, и у нас было очень много аргументов от разных коциклических многочленов в пользу того, чтобы сделать вероятность единицы равной 1 и очень мало аргументов в пользу того, чтобы сделать ее нулем, можно смело начинать приближать эту вероятность к единице. Делая так, мы, очевидно, наше распределение сужаем. Скажем, если мы, в итоге, сделали вероятность строго равной единице, то в гипперкубе всевозможных расстановок нулей и единиц в лампочки, мы половину возможностей просто отсекаем, а в другой половине вероятности расстановок удваиваем; если делаем вероятность почти равной единице, все несколько мягче. Потеря широты окупается тем, что удается много коциклических многочленов вылечить: эффект ожидаемо превосходит усилия.  

Вообще, при работе алгоритма, ему стоит запоминать моменты, когда ему удалось особенно удачно лишиться свободы, чтобы если впоследствии будет много неудач, то можно было вернуться в один из таких моментов - с удачным значением величины $\hat{S} + \rho$ и пойти в другом направлении - уже не совсем по градиенту.

Должен признать, пока я не знаю с абсолютной определенностью, каким именно свойством должна обладать искомая ситуация для того, чтобы полученная схема $Q$ была обобщаемой - есть 3 довольно близких свойства - кандидата на эту роль. Когда мы описывали алгоритм, мы как-бы оправдывали его тем, что в итоге хотим получить ситуацию максимально возможной симметрии. Это первый кандидат. Второй кандидат - это ситуация, в которой высокой симметрией обладает только часть этой ситуации, а именно, сама схема $Q$, а то, как устроены $Z_i$, - нам не важно. Третий кандидат (менее вероятный) - это результат указанного процесса сам по себе: возможно, для того, чтобы быть хорошо обобщаемой, схеме $Q$ стоит не столько обладать высокой симметрией, сколько быть полученной, руководствуясь принципом "меньше усилий, больше эффекта".

Вполне возможно, что под "усилиями"{} стоит понимать не только уменьшение широты, но и меру того, на сколько большими по модулю оказываются или собираются оказаться вещественные значения вещественозначных гейтов. В вещественнозначном случае можно штрафовать еще и за это. В случае, когда мы представляем коэффициенты коциклических многочленов как единицы в некоторых бинарных гейтах соответствующих числам от $-N$ до $N$, здесь такой штраф возможно тоже уместен: всем бинарным гейтам можно и, вообще говоря, более правильно приписать некоторую априорную вероятность единицы, в общем случае, отличную от $1/2$, и гейтам, соответствующим большим коэффициентам коциклических многочленов, приписать меньшую априорную вероятность единицы (она и вообще у всех гейтов, соответствующих большим коэффициентам коциклических многочленов будет маленькой). Как априорную вероятность единицы приписывать в общем случае, пока предстоит выяснить. В общем случае, можно указать "априорное распределение"{}, не обязательно независимое, приписывающее априорную вероятность всем возможным наборам значений, принимаемым гейтами. И тогда усилие - это мера отклонения от этого априорного распределения: чем сильнее, в результате некоторого акта распределение отклонится от априорного, тем больше соответствующее усилие. В случае равномерного априорного распределения, определение усилия реализуется отклонением суммы квадратов вероятностей всевозможных значений, в общем случае, его можно реализовать чем-то другим, предполагаю, это не сложно сделать правильным и в меру естественным образом (возможно, евклидово расстояние между векторами, составленными из вероятностей ситуаций годится).

\section{Минимум усилий, максимум эффекта.}

Когда мы решаем жизненные задачи, решение разных проблем часто не мешает друг другу и решать их можно как-бы по-отдельности. В математике несколько чаще происходит так, что решая одну проблему мы создаем другую проблему, решая другую, мы создаем третью, решая третью, мы, например, возвращаемся к первой. И общее решение требует некоторой изобретательности. Проблемы как-бы "сцепляются"{} и решение такой "сцепленной"{} проблемы - то есть нечто, что требует одновременного достижения многих разных целей вызывает у нас чувство красоты.

Слово "проблема"{} очень близко по смыслу к понятию напряженного коциклического многочлена. Именно напряженный коциклический многочлен или, более общо, нарушение симметрии является прообразом житейского понятия "проблема"{} в нашей ситуации.

Сейчас мы поговорим об идее, близкой к понятиям переобучения и регуляризации вот в каком смысле. Вообще, задачу, представляющую собой вот такое сочетание проблем, можно решать по-разному. Можно решить каждую из таких проблем по-отдельности, применяя свои усилия к каждой проблеме отдельно, а можно решать задачу более "интересно"{}, достигая по несколько целей одновременно: одни и те же усилия приводят к достижению сразу нескольких целей. Нас во многих случаях будут интересовать именно такие, "интересные"{}, решения.

Алгоритм должен всеми силами находить и запоминать моменты, когда мало усилий приносят много результата. Когда некоторые вещи происходят "сами собой"{}, без дополнительных усилий, когда нет такого, что мы каждую цель прорабатываем отдельно, лишаясь свобод для каждой из них в отдельности и в итоге получается дорого.

Можно назвать близкое к данной идее понятие "переиспользования". Полезно, когда одни и те же части строимой нами конструкции служат разным целям (например, когда защищая много разных лампочек мы воспользовались одной и той же лампочкой - когда те множества лампочек, которыми мы воспользовались при защите некоторого набора лампочек сильно пересекаются, мы лишаемся меньшей свободы, чем если бы они пересекались слабо), или когда части конструкции, предназначенные для одного, являются по совместительству составными частями чего-то другого.

Решения, в которых много всего переиспользовано и многое получено с помощью совместного достижения целей, лучше обобщается. Чем больше выражен принцип "мало усилий, много эффекта"{} при решении проблемы, тем выше обобщаемость этого решения.

Алгоритм, в каком-то смысле, идет по пути наименьшего сопротивления: если, скажем, для того, чтобы решить некоторую проблему, перед нами встает выбор, устанавливать единицу в значение одного гейта или устанавливать единицу в значение трех гейтов, причем все 4 гейта равноправны в смысле свободы которой мы лишаемся при установлении их в единицу, то мы выбираем установить единицу в значение одного гейта. В психологии есть принцип, согласно которому человек так же всегда идет по пути наименьшего сопротивления. Возможно, с этим всем оно тоже согласуется.

В прошлом препринте я дал неполное определение зоны согласованности - зоной согласованности был назван объект, который который не идеально соответствует желаемой интуиции. Сейчас хочется остановиться на этом чуть более подробно. На самом деле, зона согласованности - это множество гейтов, на котором реализуется этот самый принцип, о том, что минимальным числом усилий мы добиваемся наилучшего результата. 

Предположим, например, что мы выбрали ограничение распределения на некоторое множество гейтов, после этого, у нас однозначно восстановились значения на некотором другом подмножестве гейтов, как функции значений на первом множестве гейтов (восстановились применением коциклических многочленов или форсированием симметрии, в общем случае), затем мы восстановили что-то еще и так далее. В какой-то момент мы вернулись к некоторым гейтам из первого множества гейтов. Если выяснилось, что вновь предсказанные значения этих гейтов полностью соответствуют изначально выбранным, это хорошо: нам не нужно прилагать дополнительные усилия, чтобы удовлетворить последние предсказания, эффект их удовлетворения уже был достигнут усилием выбора ограничения распределения на первом множестве гейтов.

Зона согласованности - это подмножество пространства гейтов, набор ограничений на них (это могут быть не только пропозициональные ограничения, а например какие-то связи, даваемые коциклическими многочленами, частично проходящими "снаружи"{} рассматриваемого подмножества пространства гейтов), которые удалось в значительной степени выполнить, применив небольшое число усилий. Зона согласованности возникает, например, когда нам удается намеренно выполнить несколько нужных с точки зрения форсирования симметрии ограничений, после чего многие другие ограничения, которые должны выполняться с точки зрения форсирования симметрии (и не только те, которые должны: появление какой-то новой симметрии как всегда так же приветствуется), оказываются выполненными сами собой.

Мы выбрали именно такое название для зоны согласованности, потому что в ней выполнение многих разных ограничений и предсказаний действительно происходит согласованным образом.

Поскольку прослеживается связь между получаемыми согласно описываемому принципу решениями и регуляризацией (притянутые за уши решения не так интересны), такие решения мы будем называть оправданными регуляризацией, а так же, непринужденными (или полученными непринужденно), соответственно, сам принцип - принципом непринужденности.

Можно пронаблюдать корреляцию переиспользования и непринужденности с воспринимаемой человеком красотой во многих областях жизни. 

Игра футболиста кажется красивой, когда он очень умело пользуется своей, порой совсем небольшой, свободой. Название бренда или товара, совмещающее конец одного слова и начало другого смотрится тем эффектнее (так сказать), чем больше букв вошло в совпавший фрагмент. Это проявляется в стихах (есть глубокий смысл, а еще и рифмуется с другой строчкой). В математике мы видим красоту и чувствуем удовлетворение, когда выясняется, что один построенный нами из каких-то соображений объект подходит нам еще и по другим соображениям (собственно, об этом и речь в нашем алгоритме: если мы, скажем, выбрали распределение значений некоторых гейтов, чтобы они не образовывали противоречия с другой группой гейтов и при этом выясняется, что эти значения и с третьей группой гейтов согласуются (не образуют противоречия), то это значит, что мы "убили двух зайцев одним выстрелом"{} и построенное распределение на значенях этих признаков стоит запомнить и вообще порадоваться). Мы говорим - "емкая фраза"...

По моему мнению, человек чувствует наибольшее удовлетворение от своей деятельности, когда его поведение оправдано регуляризацией, и обратно: человека злит, когда он теряет много свободы из-за какой-то мелочи.

Если человек внутри себя поддерживает что-то наподобие распределения на значениях признаков реального мира, то я думаю, в каком-то виде, у него в голове происходит форсирование симметрии и когда он узнает какой-то новый факт, он подстраивает под него свое распределение, стараясь исправить в своих воззрениях как можно меньше (бритва Оккама, шествие по пути наименьшего сопротивления, - я думаю, именно эти вещи всплыли сейчас у вас в голове).

Когда с нами случается что-то плохое, мы спрашиваем себя - а чему же научил нас опыт? Мы пытаемся найти в своем поведении во всех плохих эпизодах одинаковую конструкцию, исправив которую, неприятностей можно было бы избежать. То есть, наша цель - внести в наше поведение как можно меньшее изменение, которое приведет к исправлению от как можно больших проблем (последнее кажется справедливо и вне контекста нахождения одинаковой конструкции в неприятном опыте, больше об этой интуиции написано в секции "Еще об абстракциях"). Это тоже проявление принципа "меньше усилий, больше результата": если бы мы для каждой неприятности внесли отдельное изменение в поведение и эти изменения были бы никак не связаны друг с другом, это вряд ли обобщилось бы на другие проблемные ситуации. А когда конструкция общая, это избавит нас и от многих проблем в будущем.

На мой взгляд, не исключено, что принцип непринужденности и форсирование симметрии - это всё, что определяет задачи в математике, которые интересуют человека и которые он успешно решает. Соответственно, класс задач, на которые нацелен наш алгоритм можно неформально определить как задачи, при попытке решения которых описываемым способом, возникает много симметрии. Разумеется, это не четко очерченный класс.

\section{Еще об абстракциях.}

Поговорим об аналоге арного признака из пространства признаков в случае пространства гейтов (первым кандидатом на роль такого аналога является, разумеется, гейт, но сейчас речь пойдет немного о другом). Предположим, что есть множество гейтов $M$ (от значений гейтов из которого и будет вычисляться "признак"), есть множество правильных предикатов, определенных на $M$, в количестве, равном арности нашего "признака"{}, и есть множество $N \supset M$ - прообраз истории вычисления признака в старом его определении, которое мы будем называть историей признака, а так же, 0-арный функциональный символ, указывающий всегда на один и тот же элемент в $N$ - гейт, в котором записано значение нашего признака. Еще должны быть ФИКСИРОВАННЫЕ предписания - аналоги предписаний, вычисляющих следующие гейты истории вычисления "старого"{} признака по предыдущим. Заметим, что они уже не обязательно вычисляют нечто следующее по чему-то предыдущему, они просто задают какие-то ограничения на значения гейтов $N$ - множество $N$ не обязано быть упорядоченым, будем понимать признак более широко. Как это работает - очевидно из только что приведенного описания. Заметим, что описанная конструкция не предполагает однозначного детерминированного вычисления значения признака по значениям предикатов на $M$, значение - лишь то, что выдает распределение в соответствующем гейте. Давайте называть такой объект современным признаком над множеством $M$.

Предположим, что у нас имеются современные признаки $D_1$, $D_2$, ..., $D_k$ над множеством $M$ и им соответствуют истории $N_1$, $N_2$, ..., $N_k$; $N_i \supset M$. Для простоты предположим, что $N_i \cap N_j = M$, при $i \neq j$. Предположим, так же, что признаки были выбраны так, что при текущем распределении все $D_i$ желаемы: для каждого из $D_i$, форсирование повсеместного равенства его единице приводит к положительному эффекту. По аналогии с тем, что происходило в подсекции "абстракции и непринужденность"{}, попробуем найти современный признак $Q$ над $M$ с историей $N$, пересекающейся с каждым из $N_i$ лишь по $M$, и набор коциклических многочленов $Z_i$, где $Z_i$ - доказательство того, что из равенства $Q$ единице всюду следует равенство $D_i$ единице всюду. Опять же, мы превращаем возможные коэффициенты доказательств в гейты, рассматриваем совместное распределение значений на полученных гейтах и предписаниях для $N$, добиваемся того, чтобы это распределение было в итоге максимально узким и полученная таким образом ситуация удовлетворяла нужным нам свойствам. 

Утверждается, что если такая ситуация была получена непринужденно, то велик шанс, что наложение данного ограничения ($Q = 1$) даст дополнительную выгоду в плане симметризации, помимо того, что $D_i = 1$, - чем более непринужденно она получена, тем ожидаемо больше этот дополнительный эффект. Например, из равенства $Q$ единице может следовать еще несколько желаемых современных признаков.

Если, опять же, проводить аналогию со старой версией алгоритма, когда все вероятности были независимы, и мы строим указанные доказательства или, более общо, аргументы, опирающиеся на то, что современный признак $Q$ равен единице, и говорящие в пользу того, что желаемые современные признаки ${D_i}$ равны единице, то, при изменении вероятностей единицы при предписаниях $Q$, усиливать вероятности наиболее охотно и сильно стоит при тех предписаниях, которые позволяют как можно больше защитить из того, что нужно защитить для того, чтобы нарождающиеся искомые аргументы были сильными: например, при таких предписаниях, при отрицании которых уже сформированы положительные коэффициенты для многих из этих аргументов (или все идет к тому, что эти коэффициенты будут положительными) - если при отрицании предписания стоит ноль, то оно защищено и нарождающиеся положительные коэффициенты не будут ослаблять указанные аргументы.

Этот пример и пример из секции "абстракции и непринужденность"{} объединяет то, что мы оба раза на что-то "нажали"{} и получили хорошую прибавку к симметричности - сразу же или, хотя-бы, в перспективе. 

А именно, во втором случае - мы "нажали"{} на предписания, образующие $Q$ (форсировали равенство их значений соответствующим константам) и в каждом из случаев $D_i$ у нас есть доказательство того, что $R(Q) = R(Q_i)$ и того, что $R(Q_i) = 1$, то есть, по сути, того, что $R(Q) = 1$. Вот, что происходит: у нас есть паралеллепипед, в котором мы хотим расставить предписания. Еще у нас есть "нарост"{} сверху на этот паралеллепипед - схема $R$. Мы изначально нажимаем на образующие предписания схемы $R$ и на гейт, в котором записано значение $R$ (устанавливаем в единицу гейт - значение $R$). И мы как-бы создаем копию паралеллепипеда, совмещая его верхнее и нижнее основания с верхним и нижним основаниями исходного паралеллепипеда. 

После этого мы нажимаем на значения предписаний в копии паралеллепипеда, в соответствии с $Q$. Последние представляют из себя новые, дополнительные ограничения на ситуацию. Давайте, для отчетливости, создадим копию схемы $R$, крепимую к паралеллепипеду так же и там же, как и сама схема $R$, эта копия выдает какой-то ответ в гейте $c$. До того, как мы нажали на рассматриваемое новое дополнительное ограничение, значение в гейте $c$ могло быть как-то нетривиально распределено. При этом, желаемое для него значение - единица. Когда мы форсируем рассматриваемое новое дополнительное ограничение ($Q$), в таком случае, в гейте $c$ значение становится равным единице и симметрия увеличивается. Но до нажатия, в гейте $c$ уже могло быть сформировано константное (или очень близкое к константному) значение 1. В таком случае, нажатие на $Q$ как-бы "подтверждает"{} равенство значения $c$ единице. Даже в таком случае, нажатие на $c$ так же будет для нас ценным и это важная тонкость.

Почему это важная тонкость. Предположим, при нажатии на $Q$, во многих или во всех $D_i$ уже достигнутое желаемое значение некоторых гейтов описываемым образом подтверждается. Это ожначает, что в других случаях, отличных от $\{D_i\}$, из нажатия на $Q$ может так же следовать некоторое желаемое значение некоторых гейтов, но при этом, на момент нажатия, такое желаемое значение может еще не быть достигнуто, а это значит, что нажатие на $Q$ в данном случае может принести уже реальную ощутимую пользу в виде повышения симметричности.

Верится мне, что можно определить и оценить количественную характеристику такого "нажатия"{}, выражающую то как много желаемых значений признаков будет получено, либо "подтверждено"{} при таком нажатии - своего рода "импульс"{} нажатия. Когда из нажимаемых значений следует много желаемых значений признаков (которые и сами по себе дают хороший импульс), импульс нажатия хороший, положительный; когда нажатие независимо от желаемых значений, то импульс нулевой (хотя, может быть и отрицательным, если нажатие лишает распределение широты); когда нажимаемые значения противоречат желаемым значениям, импульс отрицательный. Первым впечатлением могло бы быть, что импульс нажатия можно определять как то, как много симметрии это нажатие даст, возможно, с поправкой на то, какой широты оно нас лишает (вернее, на изменение $\rho$), но становитя ясно, что так делать нельзя, поскольку данная величина может быть хорошей и положительной и в случае описанного выше подтверждения (когда симметрия не обязательно увеличивается и широта не обязательно меняется).

В первом случае, описанном в первой части секции, мы нажимаем на предписания, образующие $Q$, нажимаем на значение ответа, выдаваемого $Q$, заставляя его быть равным единице, и из этих нажатий следует равенство признаков $D_i$ желаемым значениям (надо не забывать, что еще есть гейты, отвечающие правильным предикатам, и для каждого такого предиката стоит нажать на ограничения, заставляющие этот предикат быть правильным). Грубо говоря, здесь реализуется интуиция, что если нажатие получено непринужденно и дает хороший импульс в разные области пространства гейтов, то оно ожидаемо может дать хороший импульс и в другие области пространства гейтов (там происходило разбиение по случаям, здесь происходит разбиение по пространственному признаку - на области; размеется, в общем случае, стоит эти два типа разбиений комбинировать).

Итак, в обоих случаях мы нажимаем на что-то непринужденно найденное и это дает ожидаемо хороший импульс глобально.

На самом деле, когда мы находим внешнюю нейронную сеть, коррелирующую на распределении с нулем, и изменяем распределение так, чтобы она коррелировала с нулем еще сильнее, происходит примерно то же самое: установление весов нейронной сети в нужные значения и приближение ответа, выдаваемого ей к нулю - это нажатие, ожидаемо увеличивающее симметрию, а значит, ожидаемо дающее хороший импульс. И это нажатие переводит наше распределение в распределение, дающее ситуации, на которых наша внешняя нейросеть, как правило, близка к нулю достаточно сильно. Можно считать, что свойство ситуации "внешняя нейросеть на ней близка к нулю"{} является некоторым признаком $F$. Можно так же считать, что свойство ситуации "внешняя нейросеть на ней близка к нулю еще сильнее"{} является признаком $G$. Таким образом, изначально у нас был нажат признак $F$ и мы как-бы нажали на признак $G$.

Я писал о том, что нужно искать нейросеть таким образом, чтобы она коррелировала с единицей эксклюзивно. Требование эксклюзивности можно оправдать тем, что мы хотим, чтобы признак $F$ был как можно сильнее для того, чтобы признак $G$ так же был как можно сильнее. Зачем? Я объясняю это тем, что мы хотим добиться как можно более редкой симметрии. Усилия результирующие в получении корреляций значений соответственных гейтов аналогичных конструкций (далее мы будем называть их рядами значений), которые коррелируют редко - это лучше, чем такие же усилия, результирующие в получении корреляций рядов значений, которые коррелируют часто. Например, усилия, результирующие в получении корреляций рядов значений, которые равны одной и той же константе вообще всегда, нам вообще не нужны.

Этим же принципом можно объяснить то соображение, что стоит в первую очередь решать наиболее сложные проблемы (точнее, это соображение - частный случай данного принципа). Если коциклический многочлен $P$ принимает значение, близкое к правильному на намного меньшем числе ситуаций, чем число ситуаций, на которых принимает правильное значение многочлен $Q$, то усилия, позволяющие "вылечить"{} многочлен $P$ предпочтительнее, чем те же усилия, позволяющие вылечить многочлен $Q$.

Скорее, то, насколько редкая симметрия ценнее, чем частая, должно быть "зашито"{} в величины $S$ и $\hat{S}$.

Вообще, когда мы на что-то нажимаем и из этого следует какой-то признак, который обязан выполняться, стоит обращать внимание на то, насколько этот признак силен, - чем он сильнее, тем ожидаемо больше симметрии он даст. Например, то, что из нажатия последовал признак, тождественно равный единице, нам совершенно не интересно.

В случае нехождения внешней нейросети, наш принцип непринужденности принимает вид стандартного обстоятельства при обучении нейросетей, в соответствии с комментариями в конце секции "Абстракции и непринужденность"{} усилия трактуются как увеличение штрафа регуляризации (сумма квадратов весов и т.п., - смотря что мы выберем)

Мы изменяем веса нейросети так, чтобы она удовлетворяла некоторому свойству (эксклюзивно коррелировала с нулем). И если нейросеть удовлетворяет данному свойству на обучающей выборке и была получена по принципу "лучший результат при наименьших усилиях"{} (штраф регуляризации оказался маленьким), то, вероятно, она удовлетворяет данному свойству и вне обучающей выборки.

В данном случае все выглядит несколько по-другому, чем в двух других описанных случаях абстракций, хотя суть та же: мы ищем нажатие, дающее как можно больший импульс.

\section{Башня абстракций. От локального к глобальному.}

Здесь я напишу о том, что может в алгоритме происходить, не вдаваясь в детали. Я писал о том, как можно строить современный признак, из которого следуют несколько желаемых современных признаков. То есть, по сути, построить его предписания. Можно так же рассматривать современные признаки, и вообще, ограничения, "более высоких порядков". Признак второго порядка является признаком предписаний одного или более описанных признаков первого порядка (а возможно, не только предписаний, но и обычных предикатов - так сказать, предписаний нулевого порядка). Соответственно, его предписания мы будем называть предписаниями второго порядка. Таким образом, у нас возникает "башня"{} из признаков, каждый из которых стоит на основании одного или нескольких других признаков (имеет "соседей снизу") и имеет "соседей сверху"{} - признаков, стоящих на нем как на основании и "соседей сбоку"{} - признаков того же уровня, имеющих общих с ним соседей снизу или сверху.

От предписаний высших порядков нас тоже будет интересовать, чтобы они давали хороший импульс. Есть по крайней мере три основных "направления"{}, из которых мы будем черпать информацию для построения высокопорядковых признаков: "снизу"{}, "сбоку"{} и "сверху". 

Начнем с направления "снизу". Речь о том, что мы строим признак таким образом, чтобы он давал хороший импульс в соседей снизу. Можно этого добиться одним из уже описанных способов добиться хорошего импульса. Еще можно, например, обучаться этому признаку: допустим, мы хотим научиться строить ситуации, удовлетворяющие некоторому подмножеству нужных нам ограничений (а может быть, и с некоторыми добавленными ограничениями, речь о близких задачах). И нам удалось построить много таких ситуаций (возможно, порождающую модель, генерирующую такие ситуации). В таком случае, можно провозгласить равенство единице рассматриваемого признака на построенных ситуациях и равенство его нулю на ситуациях, которые нам не годятся (не удовлетворяют подмножеству нужных ограничений) - построить последние обычно не сложно. После этого можно построить признак так, чтобы он удовлетворял провозглашенным ограничениям, после чего использовать его уже для работы с основным распределением, форсировав его равенство единице (тем самым, получив новые ограничения, которые должны выполняться, которые можно и дальше включать в близкие задачи и получать для них новые признаки). 

Построить его можно, опять же, методом форсирования симметрии (скажем, возникает много коциклических многочленов, которые должны принимать константные значения, когда провозглашенное значение выполняется, и эти многочлены можно делать менее напряженными), а можно обучиться провозглашенному значению, используя backpropagation. Последнее применимо в случае, когда все гейты истории вычисления признака предполагаются вещественнозначными. Впрочем, как я уже говорил, в итоговой версии алгоритма, скорее всего все гейты будут вещественнозначными. Backpropagation тоже является частным случаем форсирования симметри. Если немножко разобраться, то в нахождении частных производных loss function по весам нейросети можно увидеть предельный случай нахождения импульса, даваемого маленькими изменениями этих весов.

Направление "сбоку". Здесь все то же самое, но импульс надлежит давать и в соседей сбоку. 

Направление "сверху". Предположим, мы выделили семейство близких задач, назовем их задачами первого порядка. Для каждой из близких задач мы выделили свое множество близких к ней задач (назовем их задачами второго порядка). Для каждой из задач первого порядка мы нашли свою абстракцию (назовем ее абстракцией первого порядка), дающую хороший импульс в задачи второго порядка (если нажать на эту абстракцию, мы получим для каждой или почти каждой из задач второго порядка желаемые значения некоторых признаков). Далее, рассмотрим подзадачи первого порядка, к которым добавлены ограничения в виде найденных для этих задач абстракций. Если мы найдем абстракцию, дающую импульс в каждую из этих преобразованных задач первого порядка, то это будет абстракция, полученная сбоку. Можно поступить иначе. Можно рассмотреть все полученные абстракции первого порядка естественным образом совмещенные и реализованные на одном и том же множестве гейтов (то есть, по сути, установить соответствие между гейтами каждой из абстракций с одним и тем же множеством). После этого можно обучиться абстракции более высокого этажа: признаку этих абстракций, который обращается на них в единицу, является, как всегда, простым (обладает высокой симметрией) и обращается в единицу как можно реже. Тем самым, он задает некоторое общее свойство этих абстракций. Далее, выделив семейство близких задач по-другому, и поступив точно так же, мы получим другое свойство абстракций (на этот раз множество гейтов этих абстракций совмещать с множеством гейтов первых абстракций мы не будем), аналогично получим и третье свойство абстракций для третьего семейства близких задач. Заметим, что все эти свойства задают ограничения не только на соответствующие абстракции, но и на исходное множество гейтов (есть какая-то абстракция, удовлетворяющая этому свойству). Предположим, что мы нашли еще одну абстракцию, хорошо согласующуюся с этими тремя свойствами, давая в них хороший импульс. Теперь самое время заполнить три образовавшихся "свободных окошка"{}, выбрав предписания для трех абстракций так, чтобы они обладали указанными свойствами и хорошо согласовались с указанной новой абстракцией. Когда мы будем это делать, информация, приходящая к ним из соответствующих им свойств и есть то, что мы нарекаем информацией, приходящей сверху.

В любом случае, откуда бы не приходила информация - снизу, сбоку, сверху или же от более далеких "соседей"{}, всегда нужно поступать так, чтобы был хороший импульс.

Мы понимаем близкие задачи, как задачи с близким множеством ограничений - подмножеством исходного множества ограничений, надмножеством, либо исходным множеством ограничений, из которого мы немножко ограничений выкинули и немножко добавили. Как ослабление, так и усиление требуемых ограничений может позволить нам легче что-то обнаружить. Когда речь о подмножестве, то  такие задачи могут быть не совсем "близкими"{}, поскольку выбрать ограничений можно и очень мало, но их все равно иногда очень нужно рассматривать. Одним словом, нужно уметь определять, насколько правильным является перенос рядов значений с другой задачи на исходную, чтобы понять, стоит ли ее называть близкой.

Можно например рассматривать подпространство пространства коциклических многочленов или конечное подмножество из коциклических многочленов. Этим надо пользоваться очень активно. Вообще, нахождение близких задач - это очень важный и актуальный вопрос, я думаю, что в итоговом алгоритме этот аспект будет сильно бросаться в глаза.

Правильней называть близкими не только задачи с множеством требуемых ограничений не сильно отличающимся от множества требуемых ограничений исходной задачи, но и задачи с множеством требуемых ограничений не сильно отличающимся от множества требуемых ограничений исходной задачи, дополненного уже произведенными на данный момент нажатиями. Последние даже более актуальны.

Кроме того, чтобы рассматривать близкие множества ограничений, можно рассматривать подмножества пространства гейтов. Решив задачу для каждого из подмножеств некоторого подмножества пространства гейтов, можно использовать полученные решения для задачи на самом подмножестве, далее полученное решение, вместе с другими, можно использовать для решения задачи на надмножестве этого подмножества и так далее. То же самое можно делать с подмножествами ограничений, решая задачу для этих подмножеств, затем объединяя группы таких ограничений и решая задачу для объединения, и так далее. Происходит решение задачи посредством перехода от локального к глобальному.

Рассмотрим подмножество множества гейтов. Можно задаться вопросом: каким свойством должно обладать подмножество ограничений, чтобы рассматриваемое множество гейтов и данное подмножество ограничений образовывали наиболее актуальную задачу. Я думаю, что первым условием должна быть простота: маленькие количества ограничений высокой симетричности предпочтительны, поскольку с ними легче работать. Вторым условием должно быть то, что данное множество ограничений задает как можно более сильное ограничение на маргинальное распределение на данном множестве гейтов. Думаю, должно быть и третье условие, связанное с симметричностью маргинального распределения.

Например, если подмножество состоит из всего одного гейта и некоторый коциклический многочлен восстанавливает значение на этом гейте однозначно, то подзадача, состоящая из одного данного гейта и одного данного многочлена актуальна: она очень сильно ограничивает маргинальное распределение на этом гейте.

Вообще говоря, возможно, не обязательно поддерживать единое распределение сразу для всех гейтов. Можно рассматривать семейство подмножеств гейтов и для каждого из этих подмножеств рассмотреть свое распределение, в идеале являющееся маргинальным распределением общего распределения. Поддерживая эти распределения, следует изменять их, помимо всего прочего так, чтобы они были как можно более сильно согласованы друг с другом на пересечениях. Вот вам еще один пример выхода в объемлющее пространство.

Признаю, что описанная в данном подразделе возня с "этажами" была описана не подробно, моя цель была передать идею того, как это может в итоге выглядеть. Я думаю, что если с этим подробно разобраться, мы получим более элегантный алгоритм, без явного разделение на этажи, обращающийся с пространством гейтов в большей степени как с однородной массой. (Хотя само понятие порядка предписания скорее всего сохранится.)

Я видел видео, где Йошуа Бенджио рассказывал о какой-то интуиции, где было множество из разных элементов, на картинке с ними что-то происходило, что-то менялось. Не помню подробностей, но была какая-то идея, возможно связанная как раз с переносом информации с меньших подмножеств на большие подмножества. Правильней было бы найти результат и привести ссылку, но к сожалению видео я не нашел. Я это к тому, что возможно идея перехода от локального к глобальному для computer science не нова, взять хотя-бы сверточные сети.

\section{Нахождение глубокой проблемы.}

Я писал, что житейское понятие проблемы в нашем случае представлено напряженным коциклическим многочленом. На самом деле, понятие проблемы можно расширить, обобщив понятие коциклического многочлена (мы уже обобщили его понятием ряда значений, сейчас речь пойдет почти о том же самом).

Рассмотрим два подмножества $A$ и $B$ множества ограничений, которые нужно выполнить - две близкие задачи. Предположим, нам удалось их в некоторой степени решить и найти для каждой из них наиболее широкое распределение как можно большей симметричности, а может быть, и семейство распределений. Предположим, так же, что нам удалось найти вещественнозначный признак, который принимает на множестве решений для $A$ множество значений как можно более отделимое от аналогичного множества значений на решениях для $B$. Можно рассматривать не один а несколько признаков и нас будет интересовать отделимость множеств их совместных значений на решениях для задач $A$ и $B$ в евклидовом пространстве, скажем, гиперплоскостью, или парой гиперплоскостей с как можно большим зазором, которые можно искать, например, методом опорных векторов.

Это и есть проблема. Нам будет сложно найти решение для объединения $A$ и $B$, поскольку такое решение должно быть нетипичным для задачи $A$ или для задачи $B$. Мы будем ее называть $(A,B)$-проблемой.

Понятие дискриминированного коциклического многочлена можно понимать по-разному. Во-первых, это может быть многочлен, который должен принимать 1, но на деле, как правило, принимает отрицательные значения. Будем рассуждать на примере малой задачи о выписывании строки. Множеством $A$ будет множество всех ограничений задачи. Чтобы получить множество $B$, нужно выкинуть из множества ограничений, собственно, условие задачи - предписания схемы и равенство единице ответа. Единственным существенным условием, которому должно удовлетворять решение задачи $B$, является согласованность ситуации. Тем самым, мы видим, что решения задачи $B$ можно описать полностью и коциклический многочлен равен 1 на каждом из них, по определению. Таким образом, множества значений для решений $A$ и $B$ хорошо отделимы нулем.  

Еще дискриминированный коциклический многочлен можно было бы определять примерно так. Это коциклический многочлен, который должен принимать 1, но в реальности значения которого распределены так, что окрестность сравнимая по длине с суммой модулей коэффициентов многочлена содержит мало возможных значений этого многочлена. Имеется в виду, что если, скажем, по выборке его значений мы видим, что подавляющее большинство этих значений сосредоточено между 20 и 21, а так же, между -17 и -15 (такое бывает, поскольку, напомню, распределение на значениях лампочек у нас уже не независимое), то ожидать, что он очень часто принимает значение 1 более странно, чем в случае, если, распределение похоже на нормальное со средним в единице. В случае такого определения, стоит отделять множество значений для решений $B$ (единица) от множества значений для решений $A$ не одной точкой, а двумя точками - границей указанной окрестности (нульмерной сферой).

Перейдем к определению абстракции проблемы. Зафиксируем множества ограничений $A$ и $B$. Пусть $C_i$ - множество как можно более разных ограничений при каждом $i$. Для всех допустимых значений $i$, попробуем найти как можно больше $(A \cup C_i,B \cup C_i)$-проблем, назовем их $P_{i,j}$, где $j$ пробегает некоторое множество значений, разное, для разных $i$ (в случае коциклических многочленов, это могут быть дискриминированные многочлены, которым разрешено пользоваться тем, что ограничения из $C_i$ выполняются - ставить на соответствующие этим ограничениям лампочки произвольный коэффициент). Далее, можно попытаться для каждого $i$ от всех $P_{i,j}$ абстрагироваться, точно так же, как мы делали в случае обычных абстракций, найдя свойство $Q_i$, которым все $P_{i,j}$ удовлетворяют. Попробуем теперь найти проблему $P$, удовлетворяющую каждому из свойств $Q_i$. Для каждого $i$, она ожидаемо принадлежит классу $(A \cup C_i,B \cup C_i)$ проблем, а значит, она ожидаемо может являться $(A,B)$-проблемой. В частности, если речь о коциклических многочленах, то полученный многочлен ожидаемо дискриминирован в общем случае и ему уже не подобает пользоваться никакими из ограничений $C_i$ (ему конечно можно, но это не рекомендуется). 

Мы получили проблему, так сказать, более глубокую, чем все $P_{i,j}$. Заметим, что в указанном процессе, при построении свойства, которому удовлетворяют подмножество рассматриваемых проблем, иногда можно брать в это подмножество и проблемы $P_{i,j}$ для разных значений $i$, иногда это может помочь, хотя удовлетворение проблемы такому свойству ничего не означает априори.

Если все проблемы предполагаются быть коциклическими многочленами, то $\{ Q_i \}$ определяют области на пространстве коциклических многочленов, а $P$ - точка на пересечении этих областей.

Как я уже писал, в жизни и в математике часто бывают ситуации, когда решая одну проблему (накладывая ограничение на мир/пространство гейтов) мы получаем другую, решив другую (сняв некоторые старые ограничения и наложив новые), получаем третью, решив третью можно опять получить первую, и так далее, дальше решить первую по-другому, получить четвертую, потом снова вернуться ко второй, и так далее. Возможно, так мы вертимся вокруг одной более глубокой проблемы, чем все эти менее глубокие проблемы указанной "проблематики". Найдя эту более глубокую проблему указанным способом, можно надеяться на то, что ее решение повлечет за собой одновременное решение сразу всех проблем указанной проблематики.

В общем случае, словом "проблема"{} мы иногда будем называть именно такую проблематику, когда возникает конфликт между некоторыми рядами значений и нам сложно "устаканить"{} этот конфликт, в том смысле, что попытка устранить конфликт в подмножестве этого набора рядов значений влечет конфликт в другом подмножестве набора этих рядов значений.

В нашем контексте, я бы назвал "идеей"{} наложение ограничения (нажатие), которое не сильно лишает нас свободы и дает сильную прибавку к симметричности (возможно, решая много проблем). Идея тем глубже, чем более выражено для нее это превосходство прибавки к симметричности над лишением свободы. Ожидается, что алгоритм будет представлять из себя взаимодействие идей и проблем: поиск и применение наиболее глубоких идей и нахождение и исправление наиболее глубоких проблем.

\section{Однородность алгоритма.}

В этой секции хотелось бы рассказать о некоторых возможностях и их применении к видоизменению самого искомого алгоритма (ищущего другой алгоритм). В особенности, хотелось бы более подробно поговорить о моих ассоциациях об объектах, визуально чем-то напоминающим топологические накрытия, чего я решил пока не делать в предыдущем препринте и рассказав о них тогда лишь немного, в контексте склеивания.

\subsection{Сцепление при работе с самой задачей.}

Я сначала напомню уже вцелом описанную схему, чем-то ее дополнив, возможно где-то повторюсь. Предположим, у нас имеется несколько близких задач, возможно, ограниченных каждая на свое подмножество гейтов. Далее мы будем называть их локальными задачами. На самом деле, поскольку мы хотим сделать все как можно менее дискретным, поддерживать стоит не бинарную принадлежность гейта такому подмножеству, а скорее степень принадлежности - вещественное число от 0 до 1 - результат применения сигмоида к вещественному числу, генерируемому нейросетью, которое в свою очередь будет изменяться нужным образом, в зависимости от нашего способа выбора подмножеств. С выбором ограничений для близких задач стоит поступать аналогично.     

Предположим, что все эти задачи (a) удалось в значительной степени решить с помощью нажатия на ограничения на предикаты на гейтах соответствующего множества и, возможно, некоторых других гейтах (назовем это нажатие решающим, полученное распределение мы так же назовем решающим), и (b) удалось в значительной степени решить такими же нажатиями, как и в пункте "a"{}, но несколькими разными способами, после чего мы абстрагировались и нашли общее свойство всех этих нажатий (которое и само является нажатием), получили распределение на вводимых ограничениях, с повышенной вероятностью удовлетворяющих данному свойству и с повышенной вероятностью решающих задачу (в таком случае мы так же назовем это распределение решающим, так же, как и нажатие на расматриваемое свойство). Мы, конечно, могли найти свойство свойства и так далее, с этими случаями все происходит так же, я их не описываю, чтобы не заграмождать текст.

Наша цель - найти общее для всех этих локальных задач единое нажатие, дающее хороший импульс в выбранные ограничения, решающие распределения и решающие нажатия для этих задач. Это общее нажатие должно как можно больше согласоваться с каждой из локальных задач и с каждым из локальных нажатий: вызывать новые и подтверждать уже присущие решающему распределению признаки, желаемые для решения локальной задачи и для реализации решающего нажатия. Если в результате общего нажатия решающее распределение потеряло широту, то это тоже плохо, поэтому общее нажатие не должно ему сильно противоречить. Одним словом, общее нажатие должно делать локальным задачам хорошо.

Искомое нажатие должно быть простым - обладать высокой степенью симметричности - именно в этом случае, при условии, что оно дает хороший импульс в локальные задачи, этот импульс ожидаемо будет обобщаться и на другие локальные задачи. Мы будем называть такое общее нажатие сцеплением.

Для каких наборов локальных задач стоит искать сцепление? Для таких, в которых решающие нажатия задач этих наборов дают хороший импульс в другие задачи набора. Мы будем называть такие наборы импульсно скоррелированными. Для импульсно скоррелированного набора задач найти сцепление легче, чем для набора задач, в котором решающие нажатия независимы и, тем более, чем для набора, в котором решающие нажатия противоречат друг другу. Сильная импульсная скоррелированность нажатий означает, что совершая эти нажатия мы, по существу, лишились одних и тех же свобод. Импульсная скоррелированность - понятие, противоположное понятию "противоречивость".

Первый из описанных способов искать нажатие, дающее хороший импульс, когда мы строим коциклические многочлены, являющиеся аргументами в сторону того, чтобы желаемые значения признаков были приняты, скорее всего обобщается до чего-то похожего на то, что было описано в секции "обратный поиск". Мы ищем нажатие, дающее хороший импульс, допустим, мы его частично построили, мы видим несколько лампочек, таких, что если мы нажмем и на них тоже, то результирующий импульс этого "добавления"{} будет хорошим, это делает эти лампочки желаемыми. Далее, мы видим, что есть несколько лампочек, нажатие на которые повлечет нажатие на много только что описанных желаемых лампочек. Это делает и эти новые лампочки желаемыми тоже. И так далее. Своего рода "обратное распространение желаемости". Нас, опять же, в большей степени интересует нажатие на лампочки, за которым следует принятие как можно большего числа желаемых значений другими лампочками (с учитываемым "весом"{} их желаемости). Из двух наборов лампочек, нажатие на которые влечет одинаковое число принятий другими лампочками желаемых значений, априори следует нажимать на тот набор, который меньше по размеру: нажатие на меньшее число лампочек повлечет к утрате меньшего числа свободы (а оцениваемый нами эффект одинаков). 

Так или иначе, хочется слить все найденные алгоритмы нахождения простого нажатия с хорошим импульсов в единый алгоритм. Вскоре мы поговорим о еще одном способе искать такие нажатия.

Следует добавить, что в описанной выше ситуации, когда мы ищем локальные решения из локальных соображений, а затем сцепить их в единое решение, предворяющей подготовкой может служить приложение усилий, специально направленных на то, чтобы добиться как можно большей импульсной скоррелированности - при поиске локальных решений нужно хотя-бы немножко "слушать"{} и решения других локальных задач, готовых или только нарождающихся. 

Так же хочется сказать вот о чем. В прошлом препринте я писал о том, что следует искать наиболее дискриминированные коциклические многочлены, у нас были 2 метода их поиска. В контексте не независимых распределений, применим только один из этих методов. Так же, по аналогии с этим, следует разрабатывать методы, чтобы искать наиболее сложные локальные задачи: искать подмножества гейтов и подмножества ограничений, для которых тяжело построить распределение на этих гейтах так, чтобы данные ограничения выполнялись.

Мы ползаем по пространству гейтов, строя маргинальные распределения на подмножествах этих гейтов, (строя нейросеть, генерирующую распределение на подмножестве гейтов, мы как-бы "описывая маргинальное распределение снизу"), а так же, совершая нажатия на признаки значений в этих подмножествах гейтов, как-бы "описывая соответствующие маргинальные распределения сверху".

И вот так ползая по пространству гейтов и строя описания маргинальных распределений снизу и сверху, мы находим наиболее сильные противоречия, наиболее сложно выполнимые локальные задачи. Эти противоречия заставляют корректировать текущие описания снизу и описания сверху наиболее сильным образом.

Скажем, если в одном из гейтов одного из маргинальных распределений сформировалось значение, противоречащее некоторым сформированным значениям гейтов и признаков других маргинальных распределений, то нужно где-то что-то поменять.

Если же в этом гейте сформировалось ровно то значение, которое из упомянутых значений гейтов и признаков следует, то это хорошо и мы на правильном пути (и можно совершить подкрепление, сделав некоторые значения более "робастными"{}, но об этом чуть позже).

Когда маргинальные распределения "выправлены"{} с помощью наиболее сильных противоречий описанным образом, у нас подготовлена почва для построения распределения на объединении множеств и для сцепления.

Если говорить о противоречиях, которые представляют из себя коциклические многочлены, такие противоречия могут быть недостроенными: такой многочлен может воспользоваться некоторым набором лампочек, возможно, на расширенном множестве гейтов, которые он только собирается защитить. В таком случае имеет смысл перестраивать распределение таким образом, чтобы этот набор лампочек защитить не удавалось: чтобы в нашем распределении данные лампочки не могли принимать такой набор значений.

Для этого могут понадобиться "противоречия второго порядка"{} - коциклические многочлены, которые "бьют по позициям"{} коциклических многочленов - противоречий первого порядка и, при этом, могут так же воспользоваться некоторым набором пока еще не защищенных лампочек.

Далее, могут появиться противоречия третьего порядка - коциклические многочлены, бьющие по этим самым лампочкам. И так далее. Возникают две противоборствующие силы: противоречия четного порядка и противоречия нечетного порядка, бьющие по позициям друг друга.

Каждая из сторон старается все делать так, чтобы противоречия, бьющие по позициям противника были как можно более сильными.

Чем больше коциклических многочленов с одной из сторон собираются воспользоваться некоторым значением лампочки, тем ценнее это значение для этой стороны и тем ценнее становятся для другой стороны коциклические многочлены, которые эту самую лампочку разрушают.

Разумеется, кроме того, чтобы бить по лампочкам, на которые опирается противник, стоит делать противоречия, которые создает противник более слабыми: например, если коциклический многочлен противника таков, что он принимает значения сильно меньше нуля, то для нас ценен будет такой кочиклический многочлен, который делает маленькой вероятность при лампочке, коэффициент первого коциклического многочлена при которой сильно меньше нуля.

Итак, нам нужны наиболее сильные противоречия для того, чтобы наиболее быстро и наиболее правильно сформировать общее распределение и, возможно, произвести сцепление. Этакий принцип узких мест на локальном уровне ~\cite{Shap}. 

Эту подсекцию можно подытожить, сказав, что мы ищем решение локальных задач из локальных сообразений, а затем находим сцепление, являющее собой решение общей задачи.

\subsubsection{Накрытие.}

Здесь мы поговорим о том, как можно построить сцепление "по частям". 

В упомянутой аналогии с накрытием, как роль накрывающего пространства, так и роль базы играют множества гейтов. Первое мы будем называть разверткой, второе так и будем называть базой, которая и будет играть роль сцепления.

Предположим, что мы, ползая по пространству гейтов, много раз встретили одну и ту же конструкцию. Мы это повторение специально искали, прилагая усилия именно для этого, и вот, нашли. Причем, встретили не обязательно непосредственно, она много раз именно "проявилась": конструкции предполагаются быть найденными как семейства признаков развертки с одинаковыми связями между этими признаками или, более общо, как семейства гейтов (опять же, с одинаковыми связями между этими гейтами) с множеством связей с разверткой, опосредованных другим, соединяющим их с разверткой набором гейтов. Это могут быть непринужденно найденные конструкции; конструкции, дающие хороший импульс в развертку; конструкции, представляющие из себя интерпретацию того, что происходит в развертке (все перечисленое - близкие вещи). И предположим, это произошло несколько раз: мы нашли несколько семейств одинаковых конструкций.

Что если предположить, что для каждого из семейств, все эти конструкции "проецируются"{} на один и тот же участок сцепления, на один и тот же участок базы? Конструкцию, которую представляет из себя этот участок, надлежит найти, проинтерпретировав конструкцию, соответствующую данному семейству так, чтобы полученные участки сцепления были согласованы на пересечении, что, разумеется, является необходимым их свойством.

Возможно, я выразился слишком размыто; чтобы было более понятно, как это примерно может происходить, приведу пример.

\subsubsection{Абстракция похожая на накрытие.}

Сейчас будет описан частный случай абстракции, описанной в секции "Абстракции и непринужденность". Предположим, опять же, что мы ищем алгоритм - схему, решающую какую-то задачу. Давайте рассуждать в тех же обозначениях, что и в только что упомянутой секции. Предположим, что нам удалось построить такую схему $Q_i$ в ряде случаев $D_1$, $D_2$, ..., $D_r$ и каждый из этих случаев определяется некоторым набором признаков $P_i$ строки $A$, обращающихся в единицу.

В основном, нас будет интересовать случай, когда $r = 1$, случай $r > 1$ в большей степени экзотичен. Так же, отличием будет то, что схема $Q$, которую мы хотим построить, общая для всех $D_i$, предполагается сильно меньшей по размеру, чем $Q_i$.

Предположим, что для каждого $i$ нам удалось построить такую конструкцию. На $Q_i$ выделено семейство маленьких множеств гейтов - окрестностей, плотно заполняющих $Q_i$. Будем называть их верхними окрестностями. В то же время, на базе $Q$, для каждого $i$, также выделено семейство маленьких множеств гейтов, которые мы будем называть нижними окрестностями.

Подобно тому, как множества окрестностей накрывающего пространства проецируются на одну и ту же окрестность базы, для каждого $i$, множество верхних окрестностей будет разбито на семейства и будет установлено соответствие между верхними окрестностями и нижними окрестностями так, что двум верхним окрестностям соответствует одна и та же нижняя окрестность тогда и только тогда, когда эти верхние окрестности принадлежат одному семейству.

Более того, для каждой из верхних окрестностей, построена схема, принимающая гейты этой верхней окрестности как вход и вычисляющяя значения в гейтах соответствующей ей нижней окрестности. При этом, должны выполняться следующие условия, из которых следует, что вычисляемые таким образом значения в гейтах нижней окрестности будут ровно теми, которые получатся, если запустить $Q_i$ и $Q$ на одной и той же строке, покрытой случаем $D_i$.

А именно, предположим, что выполняется следующее условие. Все верхние окрестности занумерованы и у каждой из этих окрестностей выделено небольшое множество окрестностей - "предшественников"{} меньшего номера. Должно выполняться следующее: если выбрать значения в гейтах предшественников таким образом, чтобы эти значения были согласованы на пересечениях предшественников и согласованы с предписаниями схемы $Q_i$, то, посредством предписаний $Q_i$, гейты самой окрестности могут быть восстановлены однозначно. Во-вторых, аналогичное должно выполняться для соответствующих этим верхним окрестностям нижним окрестностям: если выбрать в базе значения для гейтов нижних окрестностей, соответствующих предшественникам рассматриваемой верхней окрестности, то значения оставшейся части нижней окрестности, соответствующей рассматриваемой верхней окрестности, в соответствии с предписаниями базы, восстановятся однозначно.

Еще одно требуемое условие такое. Предположим, что у нас есть верхняя окрестность $A$ и список ее предшественников $B_1$, $B_2$, ..., $B_k$. Окрестности $A$ соответствует нижняя окрестность $A'$, а окрестностям $B_1$, $B_2$, ..., $B_k$ - нижние окрестности $B'_1$, $B'_2$, ..., $B'_k$, соответственно.

Предположим, что на $B_1 \cup B_2 \cup ... \cup B_k$ выбраны значения гейтов таким образом, чтобы они не конфликтовали с предписаниями развертки, а вычисленные по этим значениям с помощью схем значения гейтов $B'_1$, $B'_2$, ..., $B'_k$ не конфликтовали с предписаниями базы и не конфликтовали друг с другом. Тогда продолженные с этих значений по предписаниям базы значения гейтов $A'$ должны совпадать со значениями $A'$, полученными по схеме из значений гейтов $A$, полученных как продолжение по предписаниям развертки значений гейтов, заданных на $B_1 \cup B_2 \cup ... \cup B_k$. Предположим так же, что указанная коммутативность доказана с помощью коциклического многочлена для всех верхних окрестностей. Более подробно, у нас есть коциклический многочлен, доказывающий отсутствие конфликта на выходе схемы из $A$ в $A'$, которому можно пользоваться предписаниями развертки, предписаниями базы, предписаниями упомянутых схем, а так же, признаками входной строки (в том числе, признаками развертки и признаками базы), защищенными в случае $D_i$ (если эти многочлены не пользуются последними, то в таком случае справедливо равенство выходов, выдаваемых разверткой и базой вообще всегда, а не только в случае $D_i$).

Также следует добавить начальные и конечные условия: если верхняя окрестность расположена в непосредственной близости к входу, то каким-бы ни был вход, покрываемый случаем $D_i$, соответствие, устанавливаемое между указанной верхней окрестностью и соответствующей окрестностью базы, работает. Аналогично, если работает соответствие между верхней и нижней окрестностями, примыкающими к выходу, то соответствие с помощью схемы подразумевает, что ограничение значений гейтов окрестностей на выход одинаково (в пересечении).

Я умею строить коциклический многочлен, доказывающий, что при выполнении всех вышеупомянутых условий для случая $D_i$, развертка и база дают одно и то же. Это делает описанное частным случаем абстракции, описанной в секции "абстракции и непринужденность"{}, как я уже упоминал. То обстоятельство, что база предполагается сильно меньшей по размеру, чем развертка, позволяет предположить обобщаемость со случаев $D_1$, $D_2$, ..., $D_r$ на общий случай, даже при $r = 1$. В общем случае, чтобы сцепление оказалось хорошо обобщаемым, нужно, чтобы оно было простым. А маленькое количество гейтов - это как раз плюс к простоте: если про две конструкции мы не знаем ничего, кроме того, что вторая содержит меньше гейтов, чем первая, то вторая конструкция ожидаемо проще. Вообще, добавление нового гейта стоит рассматривать как усилие.

\subsubsection{Обобщенное накрытие.}

Можно говорить о возможности, когда найденные нами одинаковые конструкции проецируются не в одно и то же место, а в аналогичные места. Мы находим одинаковые конструкции не только для того, чтобы переносить информацию с них друг на друга, но и для создания конструкций, дающих хороший импульс: там, где появляется много одинаковых конструкций, там "пахнет"{} импульсом, там ожидаемо можно "извлечь"{} конструкцию небольшой сложности, дающую хороший импульс, пусть найденным одинаковым конструкциям и соответствуют разные места на общей конструкции. Такую общую конструкцию маленькой сложности найти было бы гораздо сложнее в случае, когда найденные конструкции не имели бы между собой ничего общего.

Возможна смесь "накрытия"{} и "обобщенного накрытия"{}, состоящая в том, что найденные одинаковые конструкции проецируются на частично перекрывающиеся места общей конструкции: у них есть одно и то же общее пересечение, в остальном они разные (хотя и аналогичные).

\subsection{Изменение общего алгоритма.}

Здесь мы поговорим о том, что форсирование симметрии можно применять и для обучения самого алгоритма - до этого мы пытались предъявить алгоритм, который делает что-то явно, но можно добавить еще и элемент обучения.

\subsubsection{Метапредписания.}

Алгоритм, который ищет общий алгоритм в большой задаче о выписывании строки, как мы уже много раз обсуждали, задается с помощью предписаний. Я верю в то, что общий алгоритм так же можно будет задать с помощью предписаний, которые мы назовем метапредписаниями: мы же считаем какие-то производные, что-то инкрементно изменяем, что-то с чем-то складываем - все это можно делать, совершая операции с вещественными числами, подобные тем, которые происходят при работе со значениями гейтов в дискретном случае.

Общий алгоритм предполагает отличаться от алгоритмов, которые он строит тем, что метапредписаний может быть довольно много и им, вообще говоря можно немножко не согласовываться: если одно метапредписание заставляет параметры изменяться с одной производной, а второе - немножко с другой, и эти производные хорошо согласуются, то когда мы усредним эти производные, вполне вероятно, что оба метапредписания реализуют свою программу.

Метапредписания реагируют на свои триггеры, которые указывают на то, что делать алгоритму в различных ситуациях. В алгоритме нет единого "стержня"{}, как это происходит с алгоритмами, которые строит общий алгоритм, он может реагировать на различные триггеры, которые он встречает при своей работе. Алгоритм вероятностный. Его ход может зависить, скажем, от случайных нормально распределенных величин - аналогов случайных битов в непрерывном случае (впрочем, первое сводит ко второму). Это позволяет на каждом из возможных участков хода работы алгоритма определить "естественость"{} хода работы на этом участке. Более вероятный ход работы более естественнен, чем менее вероятный.

Если где-то имеющиеся метапредписания конфликтуют, то их можно точно так же "лечить"{} с помощью форсирования симметрии, в частности, с помощью коциклических многочленов, как это и происходило в случае предписаний искомого дискретного алгоритма.

\subsubsection{Сделай везде все смотря по ситуации и пойми, что везде сделал одно и то же.}

Можно организовать общий алгоритм второго порядка, который наблюдает за работой общего алгоритма и обучает его. Делает это он следующим образом. Он может увидеть, как алгоритм, решая конкретную задачу, справляется с ней плохо - либо очень медленно, либо вообще, так сказать, метапредписаний для этого не хватает. При этом, он может понимать, как именно ее следовало бы решать: он понимает, как следовало бы решать задачу, пусть и не самым естественным образом. Способа сделать это с хорошей степенью естественности нет, но есть способ сделать это с удовлетворительной степенью естественности. Мы стараемся сделать это наиболее естественным способом из всех имеющихся. 

Проделав это для нескольких разных задач, у нас, в каждом из этих случаев, появится аналог локального решения, о котором шла речь несколько выше: набор метапредписаний (или признаков метапредписаний), которые, при их добавлении, реализуют то, что происходит при этом наиболее естественном из не столь естественных решений. И здесь описываемый алгоритм второго порядка может реализовать сцепление: выбрать некоторый простой набор общих метапредписаний, реализующих в каждом из случаев то, что реализуют соответствующие "локальные"{} метапредписания.

В данном контексте, сцепление можно описать фразой "сначала найди решение для каждого случая, а потом пойми, что в каждом из случаев ты сделал то же самое и продолжай делать так дальше". Мы, в каждом из случаев, сначала делаем что-то "руками"{}, а потом "сцепляем". Это специализация с последующей "антиспециализацией"{} (сцеплением) - переходом от частных решений к общему с тем же результатом для частных задач, что и указанные частные решения. Для того, чтобы конструкция сцепления хорошо обобщалась, она, как всегда, должна быть простой. Есть прямая аналогия с тем, что происходило вначале секции: мы ищем локальные решения из локальных соображений и сцепляем их. Конечно, аналогия с абстрагированием этих локальных задач так же уместна: импульс можно давать и в абстракции локальных решений. Так алгоритм обучается, сцепление - это акт обучения. Так он "ищет себя".

Иными словами, мы делаем ситуативно естественные шаги (шаги из локальных соображений), не обязательно являющиеся естественными глобально, возможно, абстрагируемся, выделяя свойство таких шагов в каждой из ситуаций, затем сцепляем, сделав эти не вполне естественные глобально шаги более естественными глобально, затем снова делаем другие ситуативно естественные шаги (в других ситуациях), затем снова сцепляем, так вот и развиваемся. При этом выбирать ситуации, в которых мы такие шаги будем делать, стоит самые сложные. То, что шаг не оказался естественным глобально, - это своего рода излом. Мы что-то ломаем для того, чтобы потом вновь что-то обрести.

Идейно такая же вещь происходила в секции "абстракции и непринужденость"{}, когда у нас была задача построить единый алгоритм, с одинаковыми предписаниями, для всех случаев, но мы временно "сломали"{} эту одинаковость: мы для разных случаев построили разные алгоритмы. Но после этого мы сцепили и полученное сцепление можно использовать для восстановления единого алгоритма, "скопировав"{} его в нашу схему - параллелепипед, вновь обретя одинаковость.

\subsubsection{Накрытие в случае развития общего алгоритма.}

Интуиция сцепления "накрытие"{} применима и к общему алгоритму.

Алгоритм работает с разными объектами (гейты, лампочки, ограничения и т.п., возможно появится много чего еще) и "на разных участках"{} множества объектов, с которыми он работает, которое тоже можно воспринимать как "развертку"{}, производит надлежащие действия. И здесь мы точно так же можем в этой предполагаемой развертке что-то там поменять, пошевелить, с минимальным ущербом для ее производительности так, чтобы выделились семейства одинаковых конструкций на разных участках, чтобы между разными семействами одинаковых конструкций все согласовалось (возможно, для этого и сами полученные конструкции нужно будет "преобразовать"{}, представить в другом виде) и она и в самом деле стала что-то "накрывать".

Опять же, следует прилагать усилия, чтобы такие семейства конструкций выделить: скорее всего, они будут выделены не на развертке непосредственно, а как множество гейтов связанных друг с другом и связь которых с разверткой опосредована некоторым числом гейтов - посредников.

Каждое полученное семейство - это конструкции, которые проецируются на одно и то же место сцепления. То, что происходит на каждом участке сцепления, тем самым, интерпретируется многими разными способами в терминах развертки. Происходит переиспользование: одни и те же части сцепления служат, так сказать, разным целям.

Мы форсируем одинаковость (или симметрию, иначе говоря) в структуре самого алгоритма, чтобы он лучше форсировал одинаковость (симметрию) в самой задаче.

В данном контексте, опять же, речь скорее о нахождении скорее коррелирующих конструкций, чем полностью совпадающих. Хотя при осуществлении накрытия, конечно, совпадение должно быть точным: скорее всего нужно слегка изменить найденные коррелирующие (аналогичные) конструкции, так, чтобы они стали совпадать строго, после чего осуществить накрытие.

\subsubsection{Обобщенное накрытие в случае развития общего алгоритма.}

Так же интересна интуиция сцепления "обобщенное накрытие"{} в контексте применения к общему алгоритму.

Мы же можем точно так же, наблюдать за тем, что наш алгоритм делает в разных задачах и в разных ситуациях или что ему надлежит делать (это отсылает нас в позапрошлую подсекцию) в разных (особенно, наиболее сложных) ситуациях и уловить в происходящем одинаковую конструкцию, возможно, одинаковый принцип и сделать сцепление, в котором эти конструкции спроецируются на аналогичные места. Точно так же можно сделать что-то в каждом из случаев, возможно, абстрагироваться в каждом из случаев один или несколько раз, после чего понять, что во всех этих случаях мы сделали аналогичные вещи, после чего сцепить, проведя обобщеное накрытие.

Поняв, что в нескольких разных случаях алгоритм применил один и тот же принцип, он находит, где этот принцип может быть применен еще и внедряет его в себя, так, что этот принцип начинает применяться алгоритмом повсеместно. Возможно, можно грубо охарактеризовать это повсеместным и согласованным его применением "на локальном уровне".

Я предполагаю, что при некоторых условиях, если предоставить алгоритму возможность развиваться, начиная с более раннего этапа, чем это собираемся сделать мы, то примером сцепления "обобщенное накрытие"{} может стать слияние оперирования внешними нейронными сетями и оперирования коциклическими многочленами: в данном случае, одинаковой конструкцией может стать то, что мы выделяем семейства похожих конструкций и делаем их еще более похожими. После обобщенного накрытия мы просто применяем эту конструкцию повсеместно - иными словами, форсируем симметрию.

Так же примером этого явления может стать нахождение абстракции методами, похожими на те, что описаны в секции "абстракции и непринужденность"{} и усиление корреляции значений, выдаваемых нейронными сетями: здесь общей конструкцией является то, что мы находим нажатие, дающее хороший импульс.

Мне встретилось много примеров применения этого же принципа, когда я проделывал работу по классификации актуальных признаков, о которой я еще раз упомяну ниже. Там повсеместно происходило усмотрение за, на первый взгляд, разными типами актуальных признаков одного и того же принципа и я понимал, что, на самом деле, рассматриваемые типы актуальных признаков можно "слить"{} в один тип.

Хочется довести алгоритм до того, чтобы он проделывал такие вещи самостоятельно.

Понимание алгоритмом одинаковости того, что он делает локально, на высоко абстрактном уровне, в разных ситуациях, и внедрение этого элемента в себя повсеместно предлагается сделать движущей силой, фактором развития алгоритма.

Если говорить о том, что происходит на сцеплении, то врядли найденная одинаковая конструкция присутствует там непосредственно, скорее, эта конструкция ПРОЯВЛЯЕТСЯ на локальном уровне.

Принцип можно рассматривать как проявление правила "сначала найди решение для каждого случая, а потом пойми, что в каждом из случаев ты сделал аналогичные вещи".

Возможно, решая конкретную задачу, стоит обучаться специально для решения именно этой задачи - метапредписания можно формировать эксклюзивно для решаемой в данный момент задачи, отчасти отказываясь от них в дальнейшем, при решении других задач.

Разумеется, все может происходить не только так, дискретно, что мы нашли одинаковую конструкцию в том, что в разных случаях делает общий алгоритм, В ГОТОВОМ ВИДЕ, а потом распространяем эту конструкцию всюду. Обычно мы будем находить несколько семейств чем-то ПОХОЖИХ конструкций и начинать делать их более похожими, находить еще больше напоминающих их конструкций и тоже сближать их с ними, понимать, как можно генерировать все такие конструкции, потом, опять же, форсировать одинаковость их всех (форсирование симметрии на метауровне) и только потом проводить обобщенное накрытие.

\par\bigskip
\par\bigskip
\par\bigskip

Пара слов в заключение этой секции. Как я уже писал, если два метапредписания указывают на то, что в одной и той же ситуации нужно делать две конфликтующие вещи, возникает конфликт: напряженный коциклический многочлен или, в общем случае, нарушение симметрии, вследствие чего, хотя-бы одно из этих метапредписаний предстоит исправить. Общий алгоритм надлежит видоизменять таким образом, чтобы таких противоречий было как можно меньше. И, вообще говоря, нам нужно больше, чем просто непротиворечивость, нам нужна высокая импульсная скоррелированность метапредписаний.

На некотором участке множества объектов, с которыми метапредписания работают, что-то происходит, информация из этого участка переходит на другой участок, чтобы все работало корректно, на этом другом участке далее должно нечто произойти. Метапредписания с другого участка должны как-бы подхватить эстафетную палочку и сработать надлежащим образом. И это действительно происходит. И это свойство форсируется.

Эта фраза "должно нечто произойти"{} вполне может реализовываться, скажем, коциклическим многочленом: некоторый коциклический многочлен, возможно, один из специально найденных наиболее информативных коциклических многочленов, о которых шла речь в первой подсекции этой секции, как-то проходит по первому и второму участку (возможно, еще где-то) и предписывает произойти тому, что должно произойти на втором участке - если этого не произошло, он напряжен. (Правильней было бы называть его метакоциклическим многочленом - аналог коциклического многочлена для работы с метапредписаниями, напряженность которого является количественной мерой конфликтов между метапредписаниями, а так же, конфликтов между тем, что происходит в реальности и тем, что должно происходить по внутреннему плану алгоритма той или иной долгосрочности. Впрочем, последнее наверно тоже можно оформить как конфликт метапредписаний, чуть подробнее об этом написано в секции "Еще мысли".) Высокая импульсная скоррелированность как раз и нужна для того, чтобы предписываемое этим коциклическим многочленом событие произошло.

Я не удивлюсь, если в мозге, по аналогии с нашим алгоритмом, деятельность одного и того же участка или группы нейронов в один и тот же момент времени может выполнять две или более совершенно разные задачи, то есть, его/ее работа может быть проинтерпретирована совершенно по-разному, в каком-то смысле, "мир накрывает мозг".

"Обобщенное накрытие"{} так же наводит на мысль, что в мозге может быть много симметрий. Возможно, "на локальном уровне"{} мозг всегда решает одну и ту же задачу и, развиваясь, "внедряет на локальный уровень"{} все новые и новые идеи. Впрочем, не знаю, есть ли в мозге силы, форсирующие симметричность неокортекса, в нашем смысле этого понятия, поддерживая аналогичность задач, выполняемых аналогичными его структурами (аналогичные структуры не обязательно должны присутствовать явно, они могут "проявляться"{} подобно описываемому выше), будь то филогенез или онтогенез.

\section{Еще мысли.}

Я много писал о том, что нужно искать похожие конструкции и делать их еще более похожими, но мало писал о том, как именно следует их искать. Рассмотрим внешнюю нейронную сеть. Вначале по seed-у, состоящему из нормально распределенных случайных величин, порождающая нейронная сеть генерирует ситуацию (правильней на этом этапе не преобразовывать полученные генеративной сетью вещественные значения от 0 до 1 в бинарные значения с соответствующими вероятностями, а оставить полученные вещественные значения как есть). Затем, возможно принимая дополнительный seed, представляющий из себя один или несколько правильных предикатов, внешняя нейронная сеть, работая некоторое время, до условного момента $M$, по ситуации строит то, что мы называем конструкцией. При этом, сама ситуация может присутствовать в неизменном виде на всех слоях сети, вплоть до момента $M$, помогая эту конструкцию вычислять. Тем самым, в момент $M$ мы получаем нашу ситуацию, на которой как-бы "нарисована"{} построенная конструкция, одна из конструкций семейства. 

Далее, по тому, что построено в момент $M$, другой частью сети вычисляется некоторая вещественнозначная величина, записанная в итоговом гейте $N$. Мы строим сеть так, чтобы итоговая величина в $N$ эксклюзивно коррелировала с некоторой константой, а затем изменяем веса генеративной модели так, чтобы эта корреляция усилилась.

То, что мы называем конструкцией, в данном контексте можно понимать двояко. Я называл конструкцией саму внешнюю нейронную сеть. Так же, под конструкцией можно понимать то, что я только что назвал - то, что "нарисовала"{} сеть на ситуации в момент $M$. В этом смысле, "конструкция"{} понимается нами в стохастическом смысле - это не обязательно набор конкретных гейтов с определенными связями между конкретными из них. Похожими их делает то, что многие эти конструкции в большой степени обладают одним и тем же признаком. Речь о названной корреляции в гейте $N$.

Естественным образом возникает возникает вопрос - почему бы не сделать, чтобы таких признаков было несколько? Кусок внешней сети что-то рисует в момент $M$ и дальше "сверху"{} параллельно друг другу присоединяются $k$ нейронных сетей, выдающих результат в гейтах $N_1$, $N_2$, ..., $N_k$ соответственно. Мы стремимся, чтобы получаемые нами на разных seed-ах вектора из значений $N_1$, $N_2$, ..., $N_k$ коррелировали с некоторым константным вектором. И мы впоследствии эту корреляцию усиливаем.

Можно так же сделать следующее. У нас есть часть нейросети, вычисляющей нечто до момента $M$, назовем ее стволом. Так же есть $k$ нейросетей сверху, вычисляющих $N_i$. Назовем их объединение короной. Так вот, можно рассмотреть семейство внешних сетей, у которых корона одинаковая, а ствол разный.

Нейросети семейства обучаются так, чтобы соответственные веса корон всегда оставались одинаковыми. На стволы таких ограничений нет. В остальном обучение и последующее изменение распределения происходят так же. Таким образом, эти корону и ствол можно прикреплять по-разному к разным местам пространства гейтов. Как именно прикреплять - отдельный вопрос. В данном случае под похожими конструкциями можно понимать, опять же, то, что нарисовал ствол, и похожими они будут, поскольку вектора из $N_i$ коррелируют, а можно понимать корону.

Важно понимать, что чтобы делать похожие конструкции более похожими, не обязательно делать более коррелирующим уже коррелирующий признак. Можно построить еще один или несколько признаков. Скажем, еще одну нейросеть, результат которой вычисляется в гейте $N_{k+1}$, и этот результат не обязан коррелировать для разных сетей семейства и разных seed-ов, и изменять веса генеративной сети так, чтобы он таки закоррелировал (конечно, остальные ${N_i}$ не должны переставать коррелировать, за этим нужно следить). Но я пока не думал, как такой признак строить, так, чтобы форсирование корреляции соответствовало интуиции "сделай конструкции более похожими".

Я писал о том, что когда ствол и корона сформированы, меняя генеративную модель, мы их уже не меняем. На самом деле, от этого правила можно немножко отказываться и в процессе изменения генеративной модели немножко изменять ствол и корону тоже. Например, если сближение некоторых признаков конструкций требует отдаления других признаков, то иногда актуально пошевелить ствол и корону так, чтобы этого отдаления не происходило. Рассуждая в большей общности, мы хотим каждый раз сближать конструкции, которые итак близки; и этот лучший кандидат на роль "итак близки"{}, вообще говоря может меняться (в смысле, меняться, в данном случае, будут ствол и корона).

Можно подытожить, сказав, что алгоритм скорее всего будет постоянно везде ползать, намеренно выискивая одинаковые конструкции, (в том числе, находя конструкции, которые часто встречались ему ранее) в том числе и в себе самом, с тем, чтобы сделать их еще более похожими.

Сначала можно задаваться целью найти две - три похожие конструкции (сейчас я про последние, с короной), потом, когда мы их нашли, можно прилагать усилия, чтобы найти много таких конструкций в других местах. Когда мы нашли много, можно попытаться все такие места научиться генерировать (параметризовать семейство конструкций): попробовать построить порождающую модель, генерирующую такие конструкции или семейства конструкций. Если строить параллель с жизнью, то можно пронаблюдать, что так часто бывает, что мы сначала замечаем совпадения, затем пытаемся вспомнить, где мы такую конструкцию видели еще, после чего пытаемся найти общую теорию, объясняющую возникновение этой конструкции во всех ситуациях класса ситуаций, в которых мы их видели и ожидаем увидеть.

Можно устроить, чтобы помимо ситуации и seed-а из правильных предикатов, ствол внешней сети принимал еще дополнительный вещественнозначный seed, тогда семейство похожих конструкций будет не дискретным, а непрерывным. Не существенное замечание, поскольку я не знаю, есть ли в этом смысл в нашем случае, но если дело дойдет до работы с непрерывной математикой, непрерывные семейства конструкций потребуются. Кстати, seed из правильных предикатов в итоге правильнее будет реализовывать через вещественнозначный seed: можно устроить нейросеть, которая по вещественнозначному seed-у строит предикат на гейтах, который всегда или почти всегда обращается в единицу ровно в одном гейте, "распределяясь"{} по этим гейтам примерно равномерно и примерно независимо от всего остального. И, разумеется, "правильность"{} - это тоже некоторая частность, от которой мы в итоге тоже скорее всего уйдем - можно генерировать предикаты, которые принимают сразу по много единиц и разное их число, в разных случаях.

Обобщением форсирования корреляции признака конструкции и константы является форсирование корреляции признака конструкции и другого признака конструкции. Когда форсировать равенство константе не получается, можно это ограничение равенства именно константе "сломать"{} и требовать равенства некоторой функции конструкции, зависящей от некоторого множества признаков этой конструкции. Потом мы поймем, что и равенство признака функции от указанного множества признаков форсировать не получается и мы решаем расширить это множество признаков. Потом, возможно, это множество признаков удастся сузить, и так далее. Равенство константе соответствует пустому множеству признаков. 

Так может происходить, когда мы встретили несколько конструкций семейства, в которых некоторый признак равен единице; мы хотели уже провозгласить тождественное равенство единице этого признака, потом мы встретили конструкцию, в которой этот признак не может быть равен единице и должен быть равен нулю; после этого, мы смотрим на другие такие конструкции и пытаемся понять - когда же этот признак равен единице, а когда нулю - мы ищем критерий. Этот критерий как раз и есть функция от как можно меньшего числа как можно более простых других признаков конструкции. (Думаю, этот акт "ломания"{} константы можно обобщить до "антисцепления"{} - акта, обратного сцеплению, когда мы заменяем конструкцию из гейтов на эквивалентную конструкцию большего размера - скажем, на развертку, в случае "антинакрытия"{} - и в результирующей конструкции что-то меняем - производим изменение, которое было бы невозможно в исходной конструкции. Собственно, возвращаясь к однородности самого алгоритма, описанного в прошлой секции, "сломать"{} нечто в самом алгоритме, решающем задачу, мы можем вообще говоря и на локальном уровне: намеренно найдя наиболее серьезную проблему, некоторую несовместимость в том, что делает алгоритм, (как мы это умеем) мы "ломаем"{} что-то на локальном уровне, отказываемся от некоторых признаков, так сказать, "локального алгоритма"{}, чтобы впоследствии этот самый локальный алгоритм сформировался более правильным образом; конечно же, при этом, одинаковость этих самых локальных алгоритмов должна форсироваться.)

Думаю, стоит отвлечься, сказав нечто о гипотезе, которую я сформулировал в конце прошлого препринта. Тогда шла речь о том, что если нам удастся построить некоторый объект, называемый коциклическим абсолютно сходящимся рядом, то это доказало бы невозможность одновременного выбора истинностных значений для всех утверждений из математики. Задумка была в том, чтобы попробовать проверить возможность того, что наш математический (а можно предположить, что и физический) мир лишь локально непротиворечив, а глобально непротиворечив быть не может и построение указанного объекта доказало бы это (существование аналогичного объекта в случае физики означало бы, что значения признаков нашей вселенной нельзя одновременно выбрать все сразу так, чтобы нигде не возникало противоречия, но можно выбрать часть: скажем, та часть вселенной, которую видит наблюдатель плюс та, которая находится в его памяти непротиворечива, но наблюдатель не сможет объять своим взором всю вселенную, в любом случае, так, чтобы противоречия нигде не было, что делало бы нашу вселенную довольно забавной - в каком-то смысле, она не могла бы существовать одновременно полностью, а лишь представляется частично, в каком-то виде виде, для наблюдателя - что-то вроде треугольника Пенроуза).

Но, к сожалению, в прошлый раз я не заметил, что из существования такой конструкции следовала бы и локальная противоречивость: поскольку ряд абсолютно сходящийся, мы можем выбрать конечное число лампочек суммарного веса сильно приближающегося к сумме всех весов лампочек. Если бы признаки, учавствующие в этих лампочках можно было выбрать непротиворечиво, то сумма коэффициентов, ограниченных на эти лампочки была бы не больше нуля - все положительные лампочки обнуляются; поэтому, при произвольном выборе значений всех остальных признаков, весов остальных лампочек не хватило бы, чтобы "добрать"{} до единицы - сумма этих весов мала. Поэтому, если мы хотим что-то подобное доказать, то конструкцию нужно модифицировать. Возможно, как-то удастся отказаться от абсолютной сходимости (основная сложность в не сходяжемся абсолютно ряде в том, что мы не можем проверять инвариантность суммы локально - неверно, что мы сможем дойти от одного набора значений до любого другого за конечное число единичных замен значений).

Если говорить о такого рода нескромных предположениях, то скорее, из этой описанной мной идеологии о том, что при выстраивании мировоззрения относительно конкретной задачи, нам полезно иногда что-то сломать (какую-то одинаковость, какую-то симметрию) с тем, чтобы впоследствии что-то обрести (то, что я говорил о "антисцеплении"{} и писал в подсекции "сделай везде все смотря по ситуации и пойми, что везде сделал одно и то же"), можно предположить, что возможно когда нибудь мы что-то аналогичным образом "сломаем"{} и в том, каким образом мы думаем о математике или о мире вцелом.

Может быть, действительно, абсолютная точность и строгость математических рассуждений является иллюзорной и на самом деле вообще все мягко, ряды значений не могут быть равны константе абсолютно строго и современная математика - это просто очень "сильная и глубокая, в нашем смысле, позиция". Может быть, когда-нибудь нам удастся объединить разные части математики, каждая из которых является непротиворечивой сама по себе, в одну большую часть и уже на этой большой части увидеть противоречие (коциклический многочлен?). При этом, само по себе предположение о противоречивости математики конечно не вносит ничего нового в науку, учитывая количество энтузиазма насчет того, чтобы это противоречие найти :) Возможно, кстати, что-то подобное описываемому происходит при фрейдовской рационализации. Человеку указывают на проблему (например, напряженный коциклический многочлен), подвергающую сомнению, скажем, моральность его поступков, он сообщает нечто - выполнение некоторого признака этого мира (как-бы совершая нажатие) и проблема уходит (многочлен становится менее напряжен), но ему указывают на другую проблему, возникающую при выполнении указанного им признака, он находит новый аргумент (новое нажатие), решающее и другую проблему и так далее и лишь когда ему указывают на более глубокую проблему, возникающую из совместного рассмотрения всех этих "маленьких"{} проблемных ситуаций, он понимает, что делает что-то не так.

Вообще, если проводить аналогии с жизнью, то я считаю, что человек оценивает признаки реального мира с двух позиций. Во-первых, он оценивает, насколько вероятно, что этот признак окажется правдой или неправдой, во-вторых, он оценивает, насколько это хорошо или плохо, если этот признак окажется правдой. С первым все более понятно - если правдивость нескольких признаков противоречит друг другу, это просто стандартная проблема, например, напряженный коциклический многочлен, со вторым нужны пояснения. Если понаблюдать за человеческим бытом и языком, быстро бросается в глаза выделяющийся аспект постоянной оптимизации - в язык зашито повсеместное оценивание человеком того, как соотносятся вещи с его целями и ценностями. И за словом "проблема"{} в языке кроется в первую очередь смысл, указывающий на то, что "все плохо"{}, а не на то, что "это все маловероятно". Поэтому нахожу разумным выделить из всех признаков реального мира признак "все будет хорошо". Выделенный признак принимает большие значения, когда предсказания, основанные на том, что происходит сейчас, указывают на удовлетворение потребностей человека/группы людей в будущем. (Вообще говоря, правильней говорить о признаке "все будет лучше, чем если бы мы поступали иначе"{} или "ожидаемый результат наших действий близок к наилучшему возможному результату в спектре всех возможных ожидаемых результатов наших действий"{}, поскольку иногда все возможные ожидаемые результаты бывают плохими.) Конечно, то, что удовлетворит одну группу людей, может не удовлетворить другую, так же, у человека могут быть разные потребности, так же то, что может удовлетворить человека через минуту или завтра может навредить ему через год. Поэтому, скорее, выделенный признак дробится на много маленьких признаков. В частности, кстати, признак, который я бы назвал "признаком лени"{}, который принимает сильно отрицательные значения, когда человеку приходится совершать много телодвижений, так что если для достижения некоторой цели человеку приходится пойти на большие затраты энергии, это становится проблемой само по себе.

Так вот, чтобы определить проблему во втором смысле ("все плохо") в нашей модели мира, можно добавить к набору значений признаков мира, которые мы хотим объявить проблемными, этот самый выделенный признак, с приписанным для него желаемым значением (или распределением на значениях). И если в этом расширенном множестве признаков со значениями мы находим уже обычную, определяемую в данных двух текстах проблему, указывающую на то, что это маловероятно, что данные признаки примут указанные значения и при этом, выделенный признак примет желаемое для нас значение, то в этом случае мы объявляем рассматриваемый набор признаков проблемой. В упомянутом чуть выше примере с рационализатором обсуждались проблемы именно такого типа. 

Именно такой подход к построению общего алгоритма подразумевался ближе к концу секции "обобщенное накрытие"{} в случае развития общего алгоритма"{}, в комментарии в скобках о метакоциклическом многочлене. Конечной целью общего алгоритма в конкретной задаче является форсирование симметрии, но при работе над задачей он вполне может ставить перед собой какие-то цели-посредники, достижение которых может привести к достижению конечной цели впоследствии, возможно, реализованные как некие обобщенные метапредписания, не суть дела. И прообразом проблемы во втором смысле, в нашем случае, могла бы стать "метапроблема"{} (в указанном комментарии в скобках я применил фразу "метакоциклический многочлен"). Скажем, такая метапроблема возникает, если "в нескольких местах"{} алгоритм собирается совершить некоторые нажатия и все эти нажатия в совокупности конфликтуют с одной или несколькими указанными целями. Тот мой комментарий, напомню, касался того, что общий алгоритм, обучаясь, с помощью сцеплений, выстраивает свои метапредписания таким образом, чтобы подобных конфликтов было как можно меньше и если он этой согласованности добьется, часто можно "надеяться"{} на то, что с объектами, с которыми работает алгоритм, в процессе его работы будут происходить именно те вещи, которые должны с ними происходить "из других соображений". 

Возвращаясь к нашей тематике, я писал о том, что есть объекты более глубокие - состоящие в конструкциях, дающих большой импульс, которые стоит менять наиболее неохотно (робастные) - а есть менее глубокие. Так же и конструкци бывают глубокими и не глубокими, робастными и не робастными: робастные - это те, которые наиболее устойчивы к изменению генеративной модели - если мы, в соответствии с алгоритмом будем изменять веса порождающей модели, то более робастные конструкции будут более ожидаемо оставаться "похожими на себя"{}, чем не робастные. Нас интересует перенос с наиболее робастных рядов значений: если в ряде значений, соответствующему семейству конструкций, эти конструкции оказались не робастными, то перенос с этого ряда менее надежен, чем в ситуации с робастными конструкциями. Соответственно, если часть конструкций робастные, а часть - нет, чтобы получить "силу переноса"{}, нужно как-то эту величину усреднять. Аналогично, при переходе от локального к глобальному, нужно в первую очередь смотреть, как взаимодействуют наиболее глубокие, наиболее робастные части локальных решений - те, которые будут присутствовать в этих решениях с наибольшей вероятностью.

Добавлю так же пару слов о непринужденности. При нахождении одинаковых конструкций (для переноса значений друг на друга и для проведения обобщенного накрытия), для нас важно то, как именно были обнаружены эти конструкции. Именно, для нас важно, чтобы они были получены непринужденно: чтобы прилагаемые усилия при их нахождении были как можно меньше по сравнению с результатом. Так же, если мы, скажем, специально искали конструкции, у которых в каком-то конкретном месте стоит единица, то переносу с ряда значений, соответствующего этому месту силу придавать следует гораздо меньше, чем если бы эта единичка появилась у нас в этом месте сама собой, непринужденно, без дополнительных приложенных именно к этому усилий.

Вообще, что значит "заметить"? Это значит, что мы произвели какие-то действия, направленные на получение информации и получили информации больше, чем ожидали получить. Это тоже проявление непринужденности.

Почему перед тем, как усиливать корреляцию в ответе, выдаваемом внешней нейронной сетью, для нас важно, чтобы ответ уже коррелировал? Это важно как раз потому, что только в этом случае конструкции, образованные внешней нейросетью можно считать непринужденно найденной: если этой корреляции нет, то и непринужденности нет - в таком случае вся совпадающая часть конструкций ограничивается весами нейросети, ничего нового в результате наших усилий по изменению весов мы не получили.

Есть еще такой аргумент в пользу того, что нас должны интересовать именно непринужденно найденные конструкции. Во-первых, заметим, что устанавливать одинаковое значение в набор показателей, являющихся соответственными показателями семейства одинаковых конструкций, априори гораздо лучше, чем устанавливать одинаковое значение в набор случайно выбранных показателей, в том же количестве. 

Покажу это на следующем примере. Предположим, у нас есть большая квадратная решетка из нулей и единиц и на ней, обращением в единицу, при всех значениях правильных предикатов, соответствующих признаков, задана такая структура: в любых двух ячейках, обе пары соответственных координат которых имеют отличие, кратное 10, значение должно быть одинаковым. Одно дело, когда мы устанавливаем, скажем, значение 1 во все ячейки главной диагонали: это соответственные элементы одинаковых конструкций - одинаковыми конструкциями являются сами эти ячейки (а соответственными элементами - значение в этих ячейках) - думаю, что нейросеть способна по правильному предикату нарисовать единичку в единственном месте, находящемся на главной диагонали и на той же вертикали, что единичка правильного предиката (или одинаковыми конструкциями могут быть вертикальные отрезки, идущие вниз, до конца, начиная от диагональных ячеек, а соответственными элементами - значение в их верхних точках). В таком случае, мы лишаемся свободы не сильно: у нас остается $2^{90}$ вариантов выбрать значения во всех узлах решетки. В случае же, когда мы устанавливаем, единицу в набор случайных узлов решетки, в размере длины главной диагонали, при условии, что решетка достаточно большая, скорее всего, у нас останется единственный вариант выбрать значения в узлах: расставить единицы всюду. Мораль в том, что когда мы приравниваем соответственные элементы одинаковых конструкций, у нас есть ощутимый шанс лишиться свободы не сильно; в случае же, когда мы приравниваем случайно выбранные элементы, скорее всего, мы лишаемся свободы очень сильно.

Так вот, если бы мы не требовали непринужденности нахождения конструкций, в которых му производим перенос по соответственным элементам, можно было бы оправдать приравнивание значений в случайно выбранных гейтах: для случайно выбранных гейтов можно тоже придумать одинаковые конструкции, соответственными элементами которых эти гейты являются, но, в общем случае, это будут либо маленькие/не похожие друг на друга конструкции (и в этом случае, сила переноса будет маленькой), либо это будут притянутые за уши конструкции (не являющиеся непринужденно найденными) - например, конструкции, составленные из гейтов, соответствующих значениям признаков и "залазиющие"{} в пространство признаков куда-то очень далеко. Требование непринужденной найденности это исключает.

Сделаю отступление о генеративных моделях, с которыми мы работаем и которые генерируют распределение, с которым мы работаем. В прошлом препринте мы работали, в основном, с независимыми распределениями. Но независимые распределения слабы: скажем, они не умеют генерировать уже никакое нетривиальное распределение на значениях ячеек только что описанной решетки, на которые наложено ограничение, что в двух ячейках, отличающихся на вектор с координатами, кратными 10, значения одинаковы. Очень хочется, чтобы генеративная модель умела генерировать нетривиальные и вообще "невырожденные"{} распределения на, скажем, решетке такие, что значения, выдаваемые моделью, каждый раз удовлетворяли ограничению, задаваемому некоторым признаком, в смысле старого препринта (при любых значениях правильных предикатов, признак равен единице). Не исключено, что есть генеративные модели, которые справляются с этой задачей лучше, чем нейросети. Думаю, правильным апгрейдом в будущем может стать "залезание форсирования симметрии внутрь"{} порождающей нейросети: симметрия форсируется и на признаках весов порождающей нейросети и значений, выдаваемых ее нейронами на промежуточных слоях, что открывает дорогу для таких вещей, как equivariance. Почему нет? Гейты на промежуточных слоях порождающей нейросети - тоже полноценные участники происходящего.

Кстати, возможно, имеет смысл обучать внешнюю нейронную сеть несколько иначе: скажем, все веса на нечетных слоях будут изменяться так, чтобы корреляция как можно сильнее усиливалась, а веса на четных слоях будут изменяться так, чтобы корреляция как можно сильнее ослабевала (мы считаем, что выходной гейт на нечетном слое). 

Есть соображение, по которому если при обучении внешней сети старым способом и только что описанным, мы получили одинаковые уровни корреляции, то далее усилять корреляцию меняя порождающую сеть целесообразнее для той сети, что обучалась только что описанным способом. Если мы обучаем сеть только что описанным способом, а дальше, параллельно продолжению этого обучения, запускаем изменение порождающей сети в сторону усиления корреляции и добиваемся очень хорошей корреляции, то, поскольку есть сопротивляющиеся слои, в каком-то смысле, в итоге мы получим больше симметрии: для сетей, полученных из данной заменой весов в сопротивляющихся слоях на другие (отличающихся от первых в пределах разумного), корреляция так же ОЖИДАЕМО усилится до очень хорошей. Нужно произнести еще некоторые слова, чтобы сделать это рассуждение более убедительным, но я сейчас этого делать не стану и вообще говоря оно имеет вероятность оказаться ошибочным. Но если оно верно, то мы получаем еще один пример, когда усиливать похожесть непринужденно найденных конструкций полезнее, чем менее непринужденно найденных (корреляция в альтернирующей сети является более непринужденно полученной, поскольку в половине случаев усилия были направлены против усиления корреляции, корреляция в большей степени получилась "сама собой"{}, чем в случае старого способа обучения сети).

Так же хочу пару слов сказать о нажатиях. Я описывал нажатия как форсирование равенства конкретных гейтов конкретным значениям. На самом деле, можно точно так же организовать нажатие на распределение: скажем, у меня бывало так, что было полезно рассмотреть распределение на значениях части лампочек и задать, совместное с ним, распределение на другом множестве лампочек (это распределение может не иметь ничего общего с реальным распределением на этих лампочках), которые связаны с первым множеством лампочек естественными ограничениями и заставить первое распределение "слушаться"{} второе, а затем использовать первое распределение в контексте самой задачи.

Еще оказывается полезным такой прием. Некоторые из значений лампочек у нас фиксированы навсегда. Мы можем от этого условия отказаться и перевести данное значение из разряда фиксированных в разряд максимально желаемых: значение не фиксировано, но есть некоторая сила, которая притягивает данное значение к нужному. Это еще один пример выхода в объемлющее пространство: до этого мы "ползали"{} по многообразию распределений, на котором это значение фиксировано, после этого мы ползаем по объемлющему многообразию распределений, на котором это значение не фиксировано. В этом объемлющем многообразии есть возможность найти более короткий путь к нужному распределению. При этом, производить такую операцию можно не только с лампочками, которые должны быть фиксированы для конкретной решаемой в данный момент задачи, но и с теми, которые должны быть фиксированы вообще всегда, постоянство которых, так сказать, следует из математических законов.

\section{Актуальные нажатия.}

В прошлом препринте я писал о том, что собираюсь выяснить, какие признаки актуальны для рассмотрения и добавления в пространство признаков. Иногда я разбираю задачи со школьных олимпиад, где просят что-то доказать или построить какой-то алгоритм. Я строю коциклический многочлен, если нужно доказательство и множество предписаний, если нужен алгоритм.

В любом случае, чаще всего для нахождения доказательства/алгоритма требуются дополнительные построения. В терминах прошлого препринта, это признаки, в терминах этого препринта это добавления новых гейтов и нажатия. (Рассмотрение признака есть нажатие на вычисляющие его предписания - нажатие, "на которое мы всегда имеем право".)

Я наблюдал за тем, какие именно построения в каждом из случаев позволяют решить задачу.

Сформировалось представление, какими свойствами должны удовлетворять нужные нам дополнительные построения, чтобы быть полезными для алгоритма при решении задачи. Я предполагаю, что такие свойства можно хорошо и отчетливо сформулировать в инвариантных терминах, не зависящих от специфики задачи. Чем в большей степени выражено такое свойство у дополнительного построения, тем ожидаемо полезнее для нас это дополнительное построение.

За все время разбора задач сформировалось более 20 типов дополнительных построений, которые нам пригождались. При более детальном рассмотрении, все эти типы посливались друг с другом и осталось всего 3 типа актуальных дополнительных построений.

Во-первых, это нажатия, дающие хороший импульс. Во-вторых, это нажатия, выявляющие проблему. И, в-третьих, это нажатия, позволяющие пропустить коциклический многочлен, доказывающий невозможность одновременного принятия конкретных значений набором из лампочек (сокращенно - пропускающее нажатие).

Для первых, как пример, можно привести признак с однозначно восстанавливаемым значением (скажем, единицей). Соответственно, нажатием будет нажатие на предписания, вычисляющие признак и нажатие на единицу в значении признака. Причем, это будет тем актуальнее, чем эксклюзивнее появляется данная единица - то есть чем более сильным ограничением эта единичка является (собственно, именно в таких случаях будет большой импульс).

Во втором случае мы ищем нажатие, которое накладывает как можно меньшее ограничение на ситуацию, но позволяет обнаружить как можно более напряженный коциклический многочлен или проблему, если рассуждать более общо. Например, это может быть просто рассмотрение набора признаков, позволяющих пропустить через свои предписания напряженный коциклический многочлен. В таком случае никакого ограничения на ситуацию просто их рассмотрение не предполагает, а проблема выявляется. При этом, если бы мы произвели нажатие, накладывающее на ситуацию сильное ограничение, плохо совместимое с задачей и другими, более ранними нажатиями, то проблема бы как раз появилась, а не выявилась, это как раз плохо. 

Вообще, выявляющийся напряженный коциклический многочлен можно интерпретировать как признак и нажатие на этот признак (заставление многочлена быть равным тому, чему он должен быть равен) будет давать хороший импульс (вообще говоря не всегда, поскольку данное нажатие может сделать напряженными другие коциклические многочлены; в таком случае может быть импульс в близкую задачу с множеством ограничений, в которое эти другие многочлены не входят). В этом смысле (и не только в этом), признак выявляющий проблему подготавливает почву для нажатия, дающего хороший импульс.

Наконец, пропускающее нажатие. Опять же, предполагается, что это нажатие накладывает как можно меньше ограничений на задачу, часто это просто рассмотрение набора признаков. Оно тем актуальнее, чем более эксклюзивной является невозможность принятия упомянутыми лампочками соответствующих значений. Указанный коциклический многочлен - доказательство должен как можно более существенно использовать сложившуюся на данный момент ситуацию, пользуясь нужными лампочками. 

При этом, сами эти лампочки должны и быть в некотором смысле актуальны: хорошо, если некоторые из них уже присутствуют в пространстве гейтов (вообще говоря не обязательно все рассматриваемые лампочки на момент формирования пропускающего нажатия должны там присутствовать, некоторые могут быть добавлены вместе с нажатием), если какая-то из них являет собой значение признака, то приветствуется, чтобы этот признак был простым и т.п. На самом деле, я допускаю, что есть тип нажатия, обобщающий данный тип нажатия и формулируемый в терминах симметрии, без упоминания коциклических многочленов.

Еще пара оговорок. Во-первых, иногда принадлежность нажатия тому или иному типу актуальных нажатий не известна нам точно, а лишь предсказана по некоторым его признакам. В таком случае, мы все равно считаем нажатие актуальным. Во-вторых, если нажатие принадлежит одному из указанных типов не для самой задачи, а для близкой задачи (правильней сказать, для локальной задачи - множество гейтов задачи тоже можно ограничивать), мы так же считаем его актуальным, хотя и чуть меньше.

Вполне вероятно, работа по поиску новых типов нажатий и их слиянию друг с другом и со старыми типами продолжится. Приведу пример, как такое может происходить. Был такой тип признаков (мы знаем, что рассмотрение признака - это суть нажатие) - признак, у которого обнаруживается неожиданная корреляция с другими актуальными признаками (на манер этого препринта, а не старого - корреляция со значениями гейтов пространства гейтов). Можно назвать это в стиле нашего текста, сказав, что у него обнаруживается эксклюзивная корреляция с другими актуальными признаками - в другой ситуации такой корреляции не было бы.

Что значит "корреляция"{} с другим признаком? Это значит, что среди пар значений, которые могут принимать наш признак и этот признак, соответственно, в нашей ситуации некоторые пары значений принимаются чаще, чем они принимаются обычно, априори, а некоторые - реже. Уже в этих словах "обычно"{}, "априори"{} слышно дыхание близких задач: "обычно"{} значит, что если рассмотреть, скажем, лишь небольшую часть требуемых ограничений (или ограничений и гейтов) и решить соответствующую близкую (локальную) задачу, то данная пара значений в наиболее общем (широком) решении этой близкой задачи принимается с некоторой частотой, которую мы и называем обычной частотой. Возникновение корреляции означает, что в результате работы с некоторыми рядами значений, эта априорная частота вынуждена нарушиться. Мы видим конфликт между рядом значений, приходящим из указанной близкой задачи и некоторыми другими рядами значений, возникающими при решении задачи. Это и есть проблема: некоторые ряды значений конфликтуют. Это позволяет отнести рассматриваемый тип нажатия к нажатиям, выявляющим проблему. В таком духе все и происходит, типы нажатий сливаются друг с другом или приводят к изменениям в общем алгоритме.

\section{Еще одна идея.}

Здесь хочется поговорить еще об одной идее, которая тоже связана с одинаковостью. Я буду рассуждать в терминах уменьшения напряженности коциклических многочленов. Я надеюсь на то, что идею можно будет обобщить и на форсирование симметрии в общем случае. Что мы делаем, когда мы вылечиваем выбранный нами набор из коциклических многочленов? Порождающая сеть генерирует ситуацию, у каждого из коциклических многочленов набора на этой ситуации вычисляется отклонение от должного значения (ошибка) и эта ошибка спускается вниз по backpropagation к весам порождающей нейросети; эти веса изменяются так, чтобы указанные ошибки как можно сильнее уменьшались суммарно. Поговорим о том, какие именно многочлены вошли в набор. Что если, скажем, половина из них полностью совпадают друг с другом? Нельзя ли этот многочлен при изменении весов учитывать с коэффициентом 1, а не с коэффициентом "половина числа элементов набора"?

Идея в том, чтобы сопоставить каждому коциклическому многочлену множество, которое мы назовем определяющим множеством коциклического многочлена, являющееся подмножеством единого для всех многочленов множества - определяющего универсума.

Хочется, чтобы для пар многочленов, у которых вектора, составленные из их коэффициентов, отличаются друг от друга не сильно, определяющие множества так же отличались друг от друга не сильно: сильно пересекались и имели маленькую симметрическую разность.

Если мы этого добьемся, мы сможем учитывать многочлены следующим образом. Если рассуждать в терминах имеющегося у нас набора каким-то образом выбранных многочленов, можно построить новый набор многочленов, выбираемых из первого набора, и многочлены нового набора учитывать с равными коэффициентами. Новый набор строится следующим образом. Мы итеративно выбираем случайную точку определяющего универсума и из всех многочленов набора, определяющее множество которых проходит через эту точку, мы выбираем многочлен случайным образом. На самом деле, среди этих многочленов, выбор может происходить не равновероятным образом, некоторым из них может отдаваться предпочтение, об этом я чуть-чуть скажу позже. Итак, мы повыбирали многочленов из исходного набора указанным образом и учитываем их с одинаковым коэффициентом. Теперь можно не бояться, что половина многочленов исходного набора совпадают: для них частота присутствия в новом наборе ожидаемо мала.

На самом деле, этот прием не избавляет нас от некоторых проблем. Например, если в исходном наборе есть два непересекающихся поднабора, по 50 совпадающих многочленов в каждом поднаборе, и еще один многочлен, определяющее множество которого является объединением определяющих множеств многочленов из разных поднаборов, (давайте здесь для удобства считать, что определяющие множества двух многочленов из одного поднабора строго совпадают, когда они совпадают лишь примерно, легко провести аналогичные рассуждения) то в таком случае, из трех фигурирующих многочленов - многочлен из первого поднабора, многочлен из второго поднабора и рассматриваемый дополнительный многочлен, третий будет учтен гораздо меньше, чем первые два. Можно справляться с этой проблемой, при выборе каждого следующего многочлена набора, находя многочлен, определяющее множество которого покрывает не одну, а, скажем, 5 случайно выбранных точек определяющего универсума. Тогда, если мы выберем 5 точек, покрываемых определяющим множеством третьего многочлена и не покрываемого никаким из определяющих множеств первых двух, то мы уже обязаны выбирать третий многочлен, а не один из экземпляров первых двух.

Как же нам строить определяющий универсум и определяющие множества? Можно поступить так. Определяющим универсумом мы объявим вероятностное пространство seed-а из нескольких нормально распределенных независимых вещественнозначных случайных величин (с самим множеством мы непосредственно работать не будем, только с seed-ами). Принадлежность точки универсума определяющему множеству многочлена набора будет вычисляться по seed-у с помощью некоторой модели: скажем, нейросеть принимает seed и узлы ее выходного слоя соответствуют многочленам набора и значение, выдаваемое нейросетью в каждом выходном узле переводится в вероятность принадлежности соответствующей точки вероятностного пространства соответствующему определяющему множеству с помощью softmax-функции. (В секции я рассуждаю об обычных множествах, но вообще говоря, правильней делать определяющие множества нечеткими множествами, в таком случае, вероятность принадлежности стоит оставлять как есть и преобразовывать ее в бинарную принадлежность не нужно.)

Поговорим о том, как мы будем эту нейросеть обучать. Изначально у нас есть предпочтения насчет того, как должны соотноситься друг с другом пары определяющих множеств: мы хотим, чтобы множества для многочленов, вектора из коэффициентов которых близки, обладали примерно одинаковым размером и маленькой симметрической разностью. Можно задать функционал (пока не уточняю, какой именно), для двух определяющих множеств $A$, $B$, по трем числам $A \cap B$, $A \setminus B$, $B \setminus A$, определяющий, насколько множества $A$ и $B$ нас устраивают - функционал большой, если устраивают и маленький, если нет. Набрав случайную выборку из seed-ов, мы можем, при текущих значениях весов сети, указанные 3 числа оценить, а значит, оценить и функционал. Можно несложно понять, как вычислить производную этой оценки функционала по весам сети. Обучение сети будет состоять в шествии по градиенту суммы функционалов по всем парам многочленов набора.

Но правильней, на мой взгляд, поступать иначе: обучать нейросеть, принимающую, помимо вышеописанного seed-а еще и коэффициенты многочлена, принадлежность определяющему множеству которого [точки, соответствующей seed-у] мы хотим выяснить и вычисляющего единственное число - вероятность принадлежности точки этому множеству. Обучать такую сеть будем способом, аналогичным вышеописанному: инициализировать множества случайной сетью, затем набирать выборки из пар коциклических многочленов, вектора коэффициентов которых отличаются не сильно и сближать множества этих многочленов (сближать тем сильнее, чем ближе вектора коэффициентов); так же следует заставить множество нулевого многочлена быть пустым.

На самом деле, в обоих случаях может возникнуть проблема, что нейросеть постоянно приходит к тому, что выдает пустые множества для всех многочленов. Если это так, нам нужна какая-то сила, которая эти множества "расталкивает"{}, не давая им исчезать полностью. (Это нужно и по другой причине, которую я позже упомяну.)

Обсудим вопрос того, как следует выбирать многочлен среди тех, множество которых проходит через выбранную нами случайным образом точку, в данном случае (я буду проводить рассуждения для случая, когда мы и правда выбираем одну точку; на случай, когда мы случайным образом выбираем несколько точек, все переносится очевидным образом). Можно построить еще одну, вторую нейросеть, принимающую наш seed и еще один дополнительный seed, по которому она будет стараться возвращать коэффициенты многочлена, множество которого проходит через точку, соответствующую первому seed-у (вернее, для того, чтобы ответ был всегда коцикличным, следует выдавать коэффициенты разложения по базису для пространства коциклических многочленов). С помощью первой нейросети мы можем подготовить множество случайно выбранных seed-ов (первого типа), для каждого из которых выбрано множество многочленов, определяющие множества которых проходят через соответствующую этому seed-у точку. Если у нас есть предпочтения насчет того, какие именно многочлены мы, в качестве ответа, выдаваемого второй нейросетью, хотим видеть больше, нам следует набирать мультимножество: более желанным многочленам сопоставлять большие коэффициенты. Когда мультимножество составлено, можно обучать вторую сеть точно так же, как обучают GAN: вторая сеть - это генератор, так же потребуется дискриминатор, обучаясь, пытающийся отличить, для каждого из seed-ов первого типа, то, что выдает генератор на этом seed-е и случайных seed-ах второго типа, от того, что должен на этом seed-е первого типа выдавать генератор (у нас есть мультимножество); задача генератора в том, чтобы дискриминатор не справился. Когда вторая сеть обучена, можно надеяться успешно применять ее и на других seed-ах первого типа. Таким образом, итоговый набор многочленов следует выбирать, итеративно выбирая случайно оба seed-а и добавляя в набор то, что выдает на них вторая сеть.

Идею определяющих множеств можно проиллюстрировать на примере, когда где-то возникло нарушение симметрии и перед алгоритмом стоит выбор: поменять распределение на трех слабо связанных друг с другом коциклических многочленах (определяющие множества которых далеки друг от друга), переведя их из разряда сильно не дискриминированных в разряд дискриминированных, либо поменять распределение на трех коциклических многочленах, любые два из которых, как вектора из коэффициентов, почти совпадают друг с другом и определяющие множества которых почти совпадают, опять же, переведя их из разряда сильно не дискриминированных (уравновешенных) в разряд дискриминированных. Интуиция подсказывает, что априори лучше выбирать второе: меняя распределение на слабо связанных многочленах, мы нарушаем симметрию априори сильнее, чем во втором случае и поскольку мы идем по пути наименьшего сопротивления, выбор стоит делать в пользу изменения распределений на второй группе многочленов. И это согласуется с тем, как мы учитываем многочлены: поскольку определяющие множества во втором случае почти совпадают, алгоритм воспринимает вторую группу многочленов почти как один и тот же многочлен и симметрия будет нарушена скорее для него, чем для трех разных.

Хочется, чтобы понятие близости определяющих множеств многочленов так же отвечало понятию чувствительности изменения распределения на одном из многочленов к изменению распределения на другом. Скажем, если два многочлена, как вектора из коэффициентов, почти совпадают, то изменения соответствующих им распределений в процессе работы алгоритма происходит почти синхронно: изменение распределения на одном подразумевает изменение распределения на другом. И их определяющие множества почти совпадают. Вообще говоря, скорее, это потребует изменения того, как мы обучаем нейросеть, строя и видоизменяя тем самым определяющие множества: если мы будем строить определяющие множества ровно так, как описано выше, то мы скорее придем к тому, что близость определяющих множеств многочленов отвечает их физической близости - тому, в какой степени совпадают вектора их коэффициентов. Хочется, чтобы указанная чувствительность к изменениям тоже учитывалась. Скорее всего, стоит изменять положение вещей в едином динамическом процессе: изменения в чувствительности значений/распределений на значениях многочленов к изменениям друг друга влекут изменения определяющих множеств, но и обратно: изменения определяющих множеств влекут изменения в чувствительности значений/распределений на значениях многочленов к изменениям друг друга.

Если среди коэффициентов многочлена мало ненулевых коэффициентов, либо много, но большинство из этих ненулевых коэффициентов - при лампочках, которые почти не меняются (может быть, защищенных, может быть, просто робастных - тех, которые принимают почти постоянное значение и это свойство устойчиво к изменениям распределения в процессе работы алгоритма), то хочется, чтобы определяющее множество такого многочлена было маленьким. Мы будем называть такой многочлен маленьким; соответственно, многочлен, не являющийся маленьким - большим. Если множество ненулевых коэффициентов одного многочлена является подмножеством множества ненулевых коэффициентов другого многочлена, причем на этом подмножестве эти коэффициенты почти совпадают, то изменения на распределении первого априори подразумевают изменения на распределении второго, но изменения на распределении второго не всегда подразумевают изменения на распределении первого. В случае, когда изменения на распределении одного многочлена влекут изменения на распределении второго, но изменения на распределении второго не всегда влекут изменения на распределении первого, хочется, чтобы определяющее множество первого многочлена было примерно подмножеством определяющего множества второго.

Возможно, имеет смысл помимо определяющих множеств для коциклических многочленов рассматривать и обучаться определяющим множествам для ограничений коциклических многочленов на области в пространстве гейтов (имеются в виду ограничения множества коэффициентов на некоторые подмножества лампочек). И помимо форсирования свойства близости множеств для близких многочленов и пустоты множества нулевого многочлена, при обучении сети, задающей определяющие множества, в таком случае, стоит также еще и набирать выборку пар ограничений коциклических многочленов на областях, согласованных на пересечении этих областей и заставлять множество соответствующего ограничения многочлена на объединении областей быть близким к объединению множеств для указанных ограничений многочленов на областях. (Возможно, не все области - множества лампочек - стоит рассматривать и добавлять в выборку.)

Другой возможный апгрейд (пока очень плохо продуманный) в том, чтобы вспомнить, что коциклические многочлены отвечают топологическим коциклам (см. прошлый препринт) и можно обучаться множествам не только для коциклов но и для коцепей; при этом, для пар коцепей, разность которых соответствует маленькому коциклическому многочлену (в частности, это цепи с совпадающей кограницей), определяющие множества стоит сближать. И в самом деле, изменение распределения для одной из таких коцепей (вернее, для соответствующего лампочного - не обязательно коциклического - многочлена) будет очень чувствительно к изменению распределения для другой. Таким образом, если, скажем, два коцикла разбиты на несколько коцепей и между этими наборами коцепей есть биекция, такая, что множества соответственных коцепей близки, то множества для коциклов скорее тоже будут близкими.

Как я уже писал, когда мы выбираем множество многочленов, которые мы будем лечить, у нас могут быть свои предпочтения. То, что мне пока приходит на ум, что можно в этом смысле предпочесть - это наиболее дискриминированные коциклические многочлены и, еще было соображение, по котором стоит предпочитать упомянутые выше маленькие многочлены (не исключено, что потом удастся найти что-то эти два типа объединяющее и обобщающее). 

В случае, когда мы конструкцию с определяющими множествами не использовали, каждый из весов порождающей сети менялся в соответствии с суммой производных отклонений многочленов набора от нормы, по этому весу. Эта сумма, деленная на число многочленов набора, аппроксимирует интеграл этого отклонения по мере на пространстве коциклических многочленов, соответствующей распределению, из которого мы выбирали многочлены набора. В случае, когда мы конструкцию с определяющими множествами используем, вес порождающей сети меняется в соответствии с, так сказать, "частичной суммой"{} указанных производных отклонений, где влияние от многочленов, несущих схожую информацию, частично "перекрывается". Возможно, такая частичная сумма тоже аппроксимирует какую-то естественным образом определяемую математическую величину, некий частичный интеграл. 

Отсюда и можно пытаться вытаскивать определение степени симметричности $S$: есть ряды значений, у каждого - свое отклонение от нормы (в норме ряд - константа), у каждого - своя предпочтительность, а по каждому из весов, у каждого своя производная, показывающая, насколько ряд значений меняется в лучшую сторону при изменении веса (можно брать, скажем производную указанной в секции "Еще мысли"{} суммы квадратов разностей значений ряда на самой задаче и значений ряда на близкой задаче); производная $S$ по весу есть частичный интеграл этой производной по всем возможным рядам значений (если удастся построить определяющие множества для рядов значений в общем случае), с учетом их предпочтительности, саму же величину $S$ можно было бы определять как частичный интеграл отклонения от нормы, с учетом предпочтительности. Ведь "уравновешенность"{} разных рядов значений часто, так же, как и в более узком случае коциклических многочленов, говорит об отчасти одном и том же. Основная проблема здесь, на мой взгляд, не в том, что частичный интеграл еще не определен, а в том, что вообще говоря в определении неявно фигурирует понятие близких задач ("по всем возможным рядам значений") и чтобы определять $S$ строго, нужно, вообще говоря, иметь в руках конкретные распределения для близких задач, но мы их ищем в процессе. Поэтому, возможно, имеет смысл определять лишь $\hat{S}$, но не $S$. Вероятно, чем "ближе"{} близкая задача к исходной, тем предпочтительнее будут соответствующие этой задаче ряды значений. Соответственно, $\hat{S}$ может оцениваться методом, аналогичным описываемому в этой секции.

Еще одним подходом к повышению широты распределения может стать отдаление определяющих множеств друг от друга: если мы будем действовать так, чтобы два ранее близких, сильно пересекающихся, определяющих множества разных признаков стали далекими, слабо пересекающимися, то изменение распределения на одном из этих признаков больше не будет влечь изменение распределения на другом, и наоборот, в той же степени, что и раньше. 

Определяющие множества, конечно же, можно определять не только для коциклических многочленов, но и для всех признаков. Если в результате наших действий, два ранее далеких, слабо пересекающихся, определяющих множества разных признаков стали близкими, сильно пересекающимися, то изменение распределения на одном из этих признаков теперь будет влечь изменение распределения на другом, и наоборот, чаще, чем раньше. Признаки стали более чувствительны друг к другу, вполне возможно, образовались какие-то новые связи, корреляции, которых не было раньше. Мы потеряли широту. Это плохо. Поэтому, действовать следует так, чтобы определяющие множества как можно сильнее расталкивать, отдалять друг от друга.

Добавлю, что когда мы сцепляем несколько современных признаков в один современный признак, как это описано в секции "Еще об абстракциях"{}, при выборе набора признаков, которые мы собираемся сцепить, стоит выбирать те наборы признаков, определяющие множества которых близки. Аналогично, когда мы видим набор проблем, скажем, коциклических многочленов, определяющие множества которых близки, стоит задуматься о том, чтобы найти более глубокую проблему, проявлениями которой данные проблемы являются.

\section{Связи с другими работами.}

Здесь хочется указать известные мне идеи и направления в науке, имеющие отношение к идеям, изложенным в этом и предыдущем текстах. Я уже упомянал о видео с Йошуа Бенджио, источник которого к сожалению был утерян. Некоторые идеи перекликаются с GAN ~\cite{GAN}, ~\cite{2013}.

На начальном этапе развития алгоритма, последний чем-то напоминал то, что происходит в реализации алгоритма игры в Go от компании Deepmind. Теоретически, этот алгоритм от Deepmind мог на меня перед этим непосредственно повлиять, приведя к появлению этих аналогичных идей, хоть я ничего тогда о нем и не знал, кроме факта его существования. Не хочу утруждать вас подробностями того, как такое могло вообще произойти, просто так жизнь сложилась. Сейчас алгоритм уже ушел далеко в сторону от методов reinforcement learning, использованных Deepmind для игры в Go, но, тем не менее, использовал схожие идеи на начальном этапе своего развития, некоторое время назад, так что считаю должным эту работу указать ~\cite{go}.

Иногда я слышу о разных идеях в науке, в которых узнаю что-то похожее на мои идеи. Во-первых, речь о гештальт психологии, в ней есть интуиция устранения конфликта, реализуемая в нашем случае коциклическими многочленами, а иногда и форсированием симметрии. Во-вторых, речь о возможно идейно близком к нашему понятию "сцепление"{} понятии "сжатие"{} в искусственном интеллекте ~\cite{whp} ~\cite{hp}. В-третьих, речь о когда-то услышанной мной идее некого, кажется, профессора лингвистики, о том, что прием, заключенный в использовании метафоры, может быть внедрен в язык и мировоззрение человека в большей степени, чем просто как художественный прием, используемый в литературе. Возможно, речь шла о чем-то похожем на перенос по аналогичным конструкциям. Если автора идеи удастся установить, я приведу ссылку в следующих работах на тему.

Так же стоит упомянуть идеи, о которых автор узнал до начала работы над задачей. Во-первых, сюда можно отнести упомянутый в тексте психологический принцип шествия человека по пути наименьшего сопротивления. Во-вторых, это сама формулировка ~\cite{AIXI}, ~\cite{GP2007}, ~\cite{alternation}. В-третьих, это идея Александра Шаповалова о принципе узких мест ~\cite{Shap}.

\section{Дальнейшая работа.}

В дальнейшем хочется разобраться, как можно "слить"{} упомянутые подходы к нахождению нажатий небольшой сложности, дающих хороший импульс (внешняя нейронная сеть, абстракции, описанные в секциях "Абстракции и непринужденность"{}, "Еще об абстракциях"{}, сцепление "накрытие"{} и т.д.), в один подход.

Следующая идея исходит из некоторого сходства у приближения друг к другу значений соответственных элементов аналогичных конструкций и сцепления "накрытие"{}. Допустим, для простоты, что в первом случае у нас имеется семейство аналогичных конструкций, в каждой из них выделено по гейту и эти гейты являются соответствеными. И мы приближаем друг к другу значения в этих гейтах. То есть, есть свойство этих соответственных гейтов, выражаемое вещественным числом, которое мы делаем одинаковым для всех этих гейтов.

Во втором случае у нас опять имеется семейство одинаковых конструкций и гейты каждого из семейств соответственных элементов мы склеиваем друг с другом. Что если совпадение/несовпадение двух гейтов тоже будет определяться равенством/неравенством некоторых вещественнозначных свойств этих гейтов?

По этой причине является естественным гейты задавать точками в евклидовом пространстве. И тем самым, склеивание соответственных гейтов семейства одинаковых конструкций может происходить постепенно: соответствующие точки постепенно приближаются друг к другу и в итоге совпадают.

Как это можно устроить? Раньше у нас гейты, как сущности, были фиксированы и у нас была генеративная модель, генерирующая значения предикатов. Теперь нам хочется, чтобы она генерировала еще и "положения гейтов в евклидовом пространстве"{}.

Можно поступать так. У нас есть множество из $n$ независимых нормально распределенных случайных величин, множество значений которых образуют евклидово пространство $U$. Гейты будут жить в пространстве $V$ размерности $m$. И у нас есть нейросеть, имеющая $n$ входов и $m$ выходов, тем самым, задающая отображение из $U$ в $V$. Заметим, что получаемые таким образом гейты как-бы размазаны по пространству $V$. Мы приближаемся к дискретной структуре гейтов, когда отображение, задаваемое нейросетью, с вероятностью, близкой к единице, бьет в одну из маленьких окрестностей на $V$ из некоторого набора маленьких окрестностей.

Заметим, что произошел еще один выход в объемлющее пространство: в процессе работы алгоритма, вместо того, чтобы всегда сохранять "целостность"{} гейтов, мы можем отказаться от этой "целостности"{} и перемещаться по пространству состояний, в которых гейты размазаны по евклидову пространству нетривиальным образом. И даже если в итоге нам указанная целостность нужна, то в процессе поиска может (и скорее всего будет) наличествовать более короткий путь по указанному объемлющему пространству (чем по исходному), который полезно использовать.

Заметим, что помимо координат точки на пространстве $U$ указанная нейросеть может принимать на вход дополнительный seed из независимых нормально распределенных случайных величин (хотя и формально не отличимый от первого seed-а) так, что по сути мы имеем не функцию - набор размазанных гейтов, а распределение на таких функциях.

Как задать предикат на таком "множестве размазанных гейтов"{}? Такой предикат представляет из себя нейросеть, которая принимает на вход $m$ чисел, которые можно интерпретировать как координаты на $V$ и выдает на выходе набор из чисел в количестве, которое можно определить как размерность предиката. Помимо указанных $m$ чисел, указанная нейросеть так же может принимать на вход дополнительный seed из независимых нормально распределенных случайных величин, что превращает этот многомерный предикат в распределение на предикатах. Предикат так же можно очевидным образом интерпретировать как функцию на $U$.

Если мы хотим задать не унарный предикат, а, скажем, бинарный, то нейросеть должна принимать не $m$ чисел, а $2m$ чисел - координаты двух точек на $V$ (плюс дополнительный seed из независимых нормально распределенных случайных величин).

Что в таком случае является аналогом внешней нейронной сети? У нас есть пространство $U$, на точках (и на парах точек, на тройках и т.п.) которых определены функции (предикаты), так же заданы несколько наборов из $n$ чисел (координаты точкек пространства $U$, каждый такой набор - аналог правильного предиката) и мы хотим по этой функции и набору наборов чисел вычислить значение некоторого функционала (допустим, опять же, многомерное).

Погрузим пространство $U$ в пространство $W$ размерности на 1 большей. Ответ, выдаваемый аналогом внешней нейросети будет являться значением некоторой функции $f$ в некоторой точке $P$ этого пространства. Если бы интересующий нас функционал зависел лишь от значения предикатов в дискретном конечном множестве точек (и наборах точек, являющихся подмножествами этого дискретного конечного множества точек), то можно было бы рассмотреть в $W$ дискретный конечный набор точек - аналогов гейтов в схеме ($P$ - одна из них), условно разделенный на слои, и сказать, что значение $f$ в такой точке является некоторой функцией от значений $f$ в установленных для нее точках предыдущего слоя. Можно убрать из этой конструкции детерминированность (вообще, сама внешняя нейросеть у нас являлась чем-то детерминированно что-то вычисляющим, но я считаю, что от всякого рода детерминированности следует уходить, так, чтобы детерминированность становилась предельным случаем чего-то более общего), сказав, что в некоторых наборах точек $W$ из соседствующих слоев (одна точка с некоторого слоя и несколько - с предыдущего) обнуляется некоторый функционал $g$ от значений $f$ в этих точках (в случае детерминированного вычисления, актуальное значение $f$ в точке, где $f$ предполагается быть вычисленной по предыдущим слоям, минус, собственно, вычисленное по предыдущим слоям должное значение значение $f$, должно быть нулевым вектором).

Теперь это легко обобщить до вычисления функционала на значениях предикатов во всех точках (не только в дискретном наборе). А именно, нужно потребовать обнуления некоторого набора функционалов $G_i$ от наборов точек $W$ и значений $f$ в этих точках (скажем, потребовать $G_1(X_1, X_2, f(X_1), f(X_2)) = 0$ и $G_2(X_1, X_2, X_3, f(X_1), f(X_2), f(X_3)) = 0$). Требование равенства $\{ G_i \}$ нулю являются непрерывным обобщением понятия нажатия. $G_i$ может быть определен не всюду (в случае нажатия мы тоже требовали установления некоторых значений не для всех лампочек, а лишь для подмножества), а для тех точек, в которых он определен, мы можем "хотеть"{} приравнять его к нулю в разной степени, тем самым, в разных точках мы будем "тянуть"{} его в направлении нуля с разной силой. Если желание приравнять к нулю такой функционал является максимально сильным на некотором подмножестве, в том смысле, что в итоге функционал абсолютно точно должен оказаться нулевым, и нулевым вне этого подмножества, то естественным вариантом происходящего является то, что данное подмножество представляет из себя подмногообразие не столь большой размерности: иначе нам будет очень сложно добиться согласованности такого нажатия для разных наборов точек.

Можно требовать выполнение всяческих связей между предикатами. Например, если у нас есть $k$-мерный и $l$-мерный предикаты, то можно протянуть между гейтами, в которых записаны их координаты, нейросеть, принимающую $k$ входов и $l$ выходов, тем самым, вычисляющую один предикат по другому (и потребовать, чтобы вычисленные значения соответствовали действительности). Разумеется, такая связь может быть реализована не только нейросетью, вычисляющей одно по другому детерминированным образом.

Более того, можно сделать трюк, подобный вышеописанному, чтобы только что указанная связь между предикатами осуществлялась не дискретной по своей природе нейросетью, а посредством непрерывного нажатия на функцию на евклидовом пространстве некоторой размерности.

Положение гейтов в пространстве и предикаты могут и сами по себе вычисляться вот таким, более непрерывным образом.

Написаное в последних двух абзацах не очень осмысленно само по себе, оно начинает иметь смысл, если мы начнем задумываться о том, что на самом деле все происходящее можно погрузить в единый универсум, в единое евклидово пространство $\Omega$ и все соотношения между гейтами, предикатами и т.д. реализуются непрерывными нажатиями. Сами же интересующие нас объекты - гейты, предикаты и т.д. будут задаваться непрерывными отображениями из $\Omega$ в $\Omega$. Скажем, если предикат бьет в одномерную прямую, то эта прямая в $\Omega$ каким-то образом вложена. Сами эти непрерывные отображения в итоге скорее тоже будут генерироваться нейросетями. Мы, опять же, будем работать не с одиночными функциями, а с распределениями на функциях.

Рассмотрим набор $A$ из небольшого числа точек $\Omega$. Рассмотрим две близкие задачи (понятие непрерывного нажатия у нас определено, так что и близкие задачи легко определить естесственным образом). Для решения одной из них распределение на значениях некоторой функции на указанном наборе точек $A$ образует некоторое, скажем так, облако. Для второй - другое облако. Движущей силой алгоритма будет, опять же, сближение этих двух облаков друг с другом (конечно, не только этих двух - подобных пар облаков будет много): распределение для каждой из двух задач следует изменять таким образом, чтобы каждая из точек одного облака приближалась ко второму облаку, а именно, старалась минимизировать сумму квадратов расстояний до точек второго облака или, что то же самое, приближалась к центру масс второго облака. Важным с идейной точки зрения частным случаем является случай, когда две близкие задачи совпадают друг с другом. Тогда наша цель в том, чтобы облако из точек - функций распределения, примененных к $A$ сжималось.

Разумеется, в качестве $A$ можно брать совершенно любой набор точек. В частности, если мы видим, что в результате некоторого телодвижения, для многих подобных наборов точек, указанное облако сжимается или, более общо, идет в нужном направлении, то это повод для того, чтобы данное телодвижение совершить.

Я писал о том, что устанавливать одинаковое значение в набор показателей, являющихся соответственными показателями семейства одинаковых конструкций, априори гораздо лучше, чем устанавливать одинаковое значение в набор случайно выбранных показателей, в том же количестве (я приводил пример с квадратной решеткой из нулей и единиц). В начале этой же секции ("Еще мысли"{}) я приводил конструкцию со стволом и короной, в которой новую корреляцию усиливать тем целесообразнее, чем больше у нас признаков, образующих корону. Последнее тоже имеет отношение к тому, что я сейчас напишу.

Скажем, устанавливать в единицу значение некоторого предиката в некотором множестве гейтов, на котором, как на аргументах - то есть как на значениях правильных предикатов, подаваемых на вход некоторому, скажем, современному признаку, довольно много таких современных признаков - которые по совместительству еще и являются простыми - обращается в единицу, априори гораздо лучше, чем устанавливать в единицу значение этого предиката в случайном множестве гейтов того же размера - в последнем случае, при наличии других уже существующих корреляций, которые необходимо сохранить (в приведенном примере фигурируют строгие корреляции - тождественные обращения в единицу некоторых признаков, но, думаю, это обобщается и на не строгие корреляции), мы потеряем ожидаемо слишком много свободы.

Как перенести эту интуицию на непрерывный случай? Что такое подмножество пространства гейтов в непрерывном случае? Можно считать, что это некоторая область вышеупомянутого пространства $V$, в котором наши размазанные гейты обитают. Само пространство V задавалось набором из $m$ чисел - то есть значениями некоторой функции в $m$ точках пространства $\Omega$. Мы хотим выделить некоторую область, на которой некоторое множество признаков выполняется. Давайте для простоты говорить о предикатах: выделяем некоторую область, на которой некоторое множество предикатов выполняется (если говорить об аналоге чего-то, что вычисляется по всему пространству гейтов, но один гейт дан в качестве правильного предиката, то в этом случае все слегка сложнее, хотя суть та же и это усложнение ничего не меняет).

Предикат - это нечто, опять же, вычисляемое по гейту и результирующее в наборе чисел, записанных в фиксированных точках $\Omega$. В непрерывном случае, вместо обращения предиката в некоторое значение строго, правильнее говорить о попадании значения в некоторую область. Таким образом, нужное нам свойство, накладываемое на область пространства $V$ в том, что довольно большой набор (чем больше, тем лучше) предикатов на этой области попадает, каждый, в довольно маленькую область соответствующего пространства (возможно, не всегда попадает, но хотя-бы часто: речь о корреляции с некоторым значением). Будем называть это свойством соответственности.

Приветствуется, опять же, чтобы такой предикат вычислялся несложным образом (в случае, когда у нас все слилось в единую непрерывную массу - $\Omega$ - это свойство несколько теряет свой смысл - неясно, что значит, что что-то по чему-то несложно вычисляется, и вообще, что значит "вычисляется"{} - детерминированности в общем случае нет - поэтому нужно разбираться, во что это свойство выливается) и приветствуется эксклюзивность: точек пространства $V$, которые попадают в эту область должно быть мало, в частности, нам сейчас не интересны предикаты, значение которых попадает в эту область вообще всегда. Из этих двух свойств простоту возможно стоит заменить на непринужденность нахождения, при этом, часто так бывает, что предикаты, не являющиеся непринужденно найденными, так же не являются и эксклюзивными, так что возможно, у непринужденности и эксклюзивности есть что-то общее, впрочем, если пока не понимаете, как можно было бы пытаться определять непринужденность нахождения в нашем случае, не берите в голову это предложение.

Итак, у нас есть область, в большой степени обладающая указанным свойством. В более дискретном случае мы говорили о приравнивании значения предиката на множестве гейтов к единице, в непрерывном случае правильно говорить о приближении предиката на указанной области к некоторой точке или, по сути, о сближении значений такого предиката друг с другом. Таким образом, наше наблюдение, при переносе на непрерывный случай, превращается в правило, согласно которому, значения предиката на области, в большой степени удовлетворяющей указанному свойству, априори более пригодны для сближения друг с другом, чем значения предиката на области, удовлетворяющей данному свойству в меньшей степени.

Заметим, что область может быть "подвижна"{}, то есть зависеть от дополнительного seed-а из случайных величин - нам ничто не мешает распространять рассуждения и для таких областей. И, наконец, правилу подвержены не только области пространства $V$, но и любые области евклидовых пространств, определяемых значениями некоторых обрабатываемых функций в некотором наборе точек $\Omega$.

Это можно неформально переформулировать как то, что если точки близки в некоторых смыслах, то их следует более охотно делать более близкими и в других смыслах. Под это правило попадает и аналог форсирования равенства друг другу соответственных элементов аналогичных конструкций, и сближение соответственных гейтов аналогичных конструкций, происходящее в "накрытии"{}. Возьмем, к примеру, форсирование равенства друг другу соответственных элементов аналогичных конструкций.

Вернемся в более дискретный случай. Если множество гейтов определяется тем, что каждый гейт из этого множества гейтов обладает рядом свойств, выполнение каждого из которых вычисляется внешней нейронной сетью, то эти гейты являются соответственными элементами семейства похожих конструкций, состоящих из, собственно, тела каждой из этих внешних нейросетей, а так же, единички, записанной в ответе каждой из этих нейросетей. Иногда это можно интерпретировать иначе: каждая из этих нейросетей может вычислять по ситуации и исходному гейту положение некоторого другого гейта и выдавать в качестве ответа значение некоторого предиката в этом другом гейте. В таком случае, конструкцией семейства похожих конструкций будет исходный гейт и набор этих самых других гейтов, положение которых вычисляется нашими нейросетями, вместе с единицей, записанной в значении соответствующего предиката. Исходный гейт относится к такому другому гейту так же, как другой исходный гейт относится к соответствующему ему другому гейту (в более общем случае, такое отношение может выражаться не только в том, что положение одного гейта вычисляется по ситуации и положению другого, но и в том, что некоторое вычисление, принимающее на вход положения этих двух гейтов и ситуацию, выдает что-то близкое к единице; конструкции, определяемые набором подобных отношений, представляется естественным так же считать похожими и нас, как всегда, будут интересовать наиболее выразительные конструкции подобного типа, (грубо) те, что располагают как можно большим числом таких отношений при фиксированном числе задействованных гейтов; вообще, свойство соответственности, конечно, следует определять в несколько более общем виде, чем примерное совпадение детерминированным образом вычисляемых свойств).

В таком случае, данное множество гейтов подлежит тому, чтобы подумать о том, чтобы сделать значение некоторого предиката на каждом из этих гейтов, равным некоторой константе. И оно этому тем больше подлежит, чем больше у нас указанных свойств. То же самое и в случае со склеиванием друг с другом соответственных гейтов одинаковых конструкций. Вообще, склеивая их, мы делаем их похожими вообще во всех смыслах.

Это все оправдывает наше построение в непрерывном случае. 

Можно возразить, что только что указанная конструкция описывает не все семейства аналогичных конструкций, поскольку, например, соответственные гейты двух таких аналогичных конструкций могут совпадать и тогда нам не вычислить по гейту другой гейт конструкции однозначно. Эту проблему можно устранить, сказав, что этот наш основной гейт (на котором мы что-то хотим приравнять к константе), так же генерируется из некоторого источника, из некоторой точки пространства $\Omega$, из которой, на тех же правах, генерируются и остальные гейты конструкции. Тогда гейты могут совпадать и это нормально. Из некоторого источника генерируются $k$ гейтов - "другие"{} элементы конструкции и из него же генерируется и основной, $k+1$-й гейт. Это еще раз наводит нас на вопрос о том, каким свойством вычисление $k+1$-го гейта должно удовлетворять, чтобы форсирование некоторого значения в нем нас в большей степени интересовало.

Я писал о том, что вычисление $k+1$-го гейта должно быть простым. Есть еще такой момент. Предположим, что вычисление каждого из первых $k$ гейтов должно состоять из некоторого вычисления $A$, за которым следует некоторое дополнительное простое вычисление. Сравним два кандидата на роль вычисления $k+1$-го гейта: первый состоит из вычисления $A$, за которым следует некоторое простое вычисление, а второй представляет из себя вычисление такой же сложности, что и первый кандидат, но в нем не присутствует вычисление $A$. Представляется логичным, что априори первый кандидат является более подходящим. И вообще, хочется, чтобы вычисление $k+1$-го гейта имело как можно больше общего с вычислением первых $k$ гейтов, чтобы оно оперировало теми же объектами, что и указанные вычисления.

Это наводит на мысли об условной колмогоровской сложности. Напомню, условная колмогоровская сложность строки $B$ относительно строки $A$ определяется как минимальная длина описания машины Тьюринга, выдающей строку $B$ на входе $A$. Если нам удастся определить условную колмогоровскую сложность в непрерывном случае, то возможно стоит требовать для вычисления $k+1$-го гейта маленькой условной колмогоровской сложности относительно вычислений первых $k$ гейтов.

Чтобы такую вещь определять, начинать стоит с определения обычной колмогоровской сложности в непрерывном случае. Самая простая вещь, которая приходит на ум, как я уже обсуждал, это мера того, сколько в объекте симметрии. Например, это может быть частичная сумма по всем рядам значений убывающей функции дисперсии этого ряда значений (о частичной сумме шла речь в секции "Еще одна идея"{}; в более непрерывном случае нужно будет работать с частичным интегралом). Эту убывающую функцию стоит выбирать так, чтобы она была равна, скажем, единице, когда дисперсия нулевая, и постепенно устремлялась к нулю при устремлении дисперсии к бесконечности. Рядов значений у нас бесконечно много, поэтому обычная сумма функции дисперсии являлась бы бесконечной, частичная же сумма, при правильном выборе деталей ее определения, будет конечной.

Вообще говоря, то, что чем больше совпадающих признаков у семейства объектов, тем логичнее предположить совпадение и нового признака, часто проявляется и в жизни. Когда совпавших признаков мало, перенос не очень уместен, мы говорим "это разное"{}.

Итак, нам более интересны корреляции на тех рядах соответственных элементов, которые являются, так сказать, наиболее соответственными (так и будем их далее называть; то, насколько соответственны эти элементы, мы так же будем называть степенью соответственности указанной корреляции), для которых есть как можно больше "дружественных"{} рядов соответственных элементов. Понятие степени соответственности корреляции уместно и в более непрерывном случае.

Если одно нажатие приводит к появлению корреляции большей степени соответственности, а другое нажатие приводит к появлению корреляции меньшей степени соответственности, то если мы про эти нажатия больше ничего не знаем, выбирать стоит первое. Это так, потому, что, во-первых, интуитивно кажется понятным, что появление первой корреляции производит больший симметризующий эффект. Появление второй же корреляции производит больший разрушительный эффект: в первом случае приближение соответствующих значений к константе, для разных точек, как-бы подразумевает друг друга, во втором случае, такое приближение для разных точек приводит к в большей степени не связанным друг с другом последствиям, и форсирование второй корреляции вносит априори больше напряженности и хаоса.

Аналогичное можно сказать в случае, когда происходит не появление новой корреляции, а что-то похожее на парралельный перенос: если одно нажатие приводит к сдвигу значений, образующих корреляцию высокой степени соответственности на некоторой области, а второе нажатие приводит к сдвигу значений, образующих корреляцию более низкой степени соответственности на другой области (тем самым, состоящей из точек, имеющих меньше отношения друг к другу), то априори выбирать нужно первое. Сдвигать нужно "цельные куски"{}.

Пока не понятно, что нужно предпочесть в случае разрушения двух корреляций разной степени соответственности.

Так или иначе, если наш алгоритм будет отчасти устроен как оптимизация некоторого функционала, то приближение более соответственных значений к соответствующему среднему значению для близкой задачи возможно следует как-то награждать, например каким-то награждающим слагаемым, включенным в функционал (такое приближение в принципе хорошо бы награждать, но когда точки являются в высокой степени соответственными, следует награждать сильнее). Такое слагаемое может представлять из себя некий суммарный эффект по всем возможным областям. Вообще понятие "области"{} скоррее всего так же будет "смазано"{}, в том смысле, что более естественно рассматривать не подмножество $\Omega$, а функцию ("характеристическую"{}) из $\Omega$ в $[0,1]$.

Так же следует увеличивать степень соответственности ряда соответственных элементов в ответ на увеличение степени соответственности дружественных рядов соответственных элементов (и, соответственно, уменьшать в ответ на уменьшение). И это логично: если у нас есть семейство рядов соответственных элементов, составленных из соответственных элементов семейства аналогичных конструкций, то если, скажем, для большинства из этих рядов, в результате изменения ситуации, образовалось свое, дополнительное, семейство аналогичных конструкций, в котором элементы соответствующего ряда являются соответственными, то мы можем естественным образом составить новое семейство аналогичных конструкций, каждая из которых составлена из конструкции исходного семейства аналогичных конструкций, а так же, из набора из нескольких конструкцй указанных дополнительных семейств аналогичных конструкций. Наличие этого нового семейства конструкций делает ряды соответственных элементов исходного семейства конструкций более соответственными, причем тем в большей степени, чем выразительнее оказались конструкции указанных дополнительных семейств.

Нечто, вычисляющее по набору гейтов/набору точек вещественное значение, соответствует понятию предиката в логике, нечто вычисляющее по набору гейтов/точек положение некоторого гейта соответствует понятию функции в логике. Надо сказать, склеивание гейтов (приближение друг к другу многомерных значений, я бы назвал это номинальным склеиванием) в общем случае более информативно, чем приближение друг к другу одномерных значений предикатов.

Политика делать уже похожее в чем-то еще более похожим может дать еще и такой приятный эффект. Когда по точкам пространства $U$ мы генерируем точки пространства $V$ (пространства гейтов), на $V$ естественным образом будут образовываться "сгустки"{} - области, маленькие относительно своего прообраза в $U$. Нетрудно видеть, что для прообраза точки, близкой к границе сгустка, в соответствии с тактикой усиления степени соответственности для областей уже высокой степени соответственности, сил, заталкивающих образ этой точки в $V$ с окрестности границы сгустка внутрь сгустка, априори будет больше, чем сил, тянущих этот образ наружу сгустка. Тем самым, сгустки будут иметь тенденцию сгущаться еще сильнее. Таким образом, гейты будут формироваться сами собой, возможно для этого не потребуется дополнительных специально для этого вводимых сил алгоритма, направляющих изменение ситуации (а формирование гейтов нам нужно: скажем, если мы хотим найти строгое доказательство теоремы, нам в итоге нужен коциклический многочлен, оперирующий дискретным множеством отделеных друг от друга гейтов).

Размазывание гейтов по евклидову пространству так же обещает чему-то приятному появиться и в вопросе перехода от локального к глобальному.

Локальная задача требует нужного взаимного геометрического расположения частей конструкции, которую мы строим, друг к другу (взаимное расположение образов некоторых функций друг к другу, частота попадания значений пары функций на некоторой области в декартово произведение двух других областей, и т.д.). Построив несколько разных решений локальной задачи, мы можем найти общее свойство найденных решений (найти абстракцию) и использовать его при построении решений для более высокоуровневой локальной задачи, являющей собой объединение ограничений для набора из нескольких локальных задач, наряду с абстракциями, найденными для других локальных задач набора.

Или, скажем, первая часть построенной конструкции диктует второй части конструкции располагаться одним образом, а третьей части конструкции - располагаться другим образом. И когда мы вторую и третью части конструкции исходя из этого правильным образом располагаем, мы можем обнаружить, что положения второй и третьей части конструкции приятным образом согласуются друг с другом из других локальных соображений. При приближении к правильному решению, такого рода согласованность нарастает лавинообразно.

Причина лавинообразности в том, что если мы хотим построить объект, удовлетворяющий некоторому простому свойству и у этого свойства есть много следствий (а может быть, и просто коррелирующих с ним утверждений) и нам удалось построить ПРОСТОЙ объект, удовлетворяющий довольно большому подмножеству этих следствий, то этот простой объект с высокой вероятностью удовлетворяет нужному простому свойству и, в частности, удовлетворяет и остальным следствиям тоже (близко к бритве Оккама). А если и удовлетворяет не всем, то с высокой вероятностью удовлетворяет довольно большому их подмножеству, образующему некоторое сильное простое следствие из нужного простого свойства. Но если и так, то объект все равно был построен не зря, поскольку построив несколько подобных объектов мы можем абстрагироваться и извлечь пользу из полученной абстракции.

Хочется указанное свойство активно использовать, разумеется, в стохастичной его форме.

По сути, именно это происходило в начале секции "Еще об абстракциях"{}. Тогда мы тоже искали простое нажатие, из которого следовало несколько желаемых современных признаков и мы предполагали, что если мы его найдем, то из него последует и много других желаемых современных признаков. 

Какое именно подмножество следствий выбирать для того, чтобы такой объект строить? Думаю, ведущее требование в том, чтобы эти следствия как можно сильнее противоречили друг другу (речь о сильно отрицательной импульсной скоррелированности) - нужна глубокая проблема. (А глубокая проблема, кстати говоря, кажется, симметрично, будет наиболее трудно обходимой, если она будет простой.)

Кстати, из общих соображений, построение проблемы на пути к построению проблемы может помочь в построении решения. Так что если, скажем, поиск напряженных коциклических многочленов в пространстве коциклических многочленов, описанный в прошлом тексте удастся адекватно обобщить на поиск глубокой проблемы в общем случае, а затем применить аналогичный метод для нахождения препятствий к построению проблемы, то это может помочь в построении решения и найденный метод можно так же пытаться слить во что-то одно вместе с остальными предложенными методами.

Можно так же подумать о слиянии форсирования симметрии со сцеплением и принципом непринужденности. Что значит, что одно и то же действие привело к достижению многих разных целей? Это значит, что есть согласованность намерений, определенность в том, что нам стоит делать. Иными словами, может возникнуть определеность в том, что мы намерены сделать, может возникнуть определенность в том, что мы хотим сделать. Скажем, так бывает, когда аргументов двигать некоторый вещественный параметр в одну сторону гораздо больше, чем аргументов двигать его в другую сторону. Но определенность - это и есть симметрия.

Так же возможно прослеживается связь между степенью соответственности и определяющими множествами: признаки, определяющие множества которых близки, в некоторых отношениях похожи на признаки высокой степени соответственности (мы можем определять степень соответственности не только для рядов из элементов/областей, но и для пар элементов/пар точек): например, они чувствительны к изменениям значений друг друга (для соответственных элементов это верно возможно далеко не всегда, но описанная выше чувствительность может указывать на идейную связь); и близость определяющих множеств, и степень соответственности можно назвать степенью одинаковости элементов.

\setcounter{secnumdepth}{0}
\section{Список литературы}
\begingroup
\renewcommand{\section}[2]{}%

\end{document}